\documentclass{article}
\usepackage[table]{xcolor}

\usepackage[scr=rsfs]{mathalpha}

\usepackage{tikz}
\usetikzlibrary{}
\usetikzlibrary{matrix,chains,positioning,decorations.pathreplacing,arrows}
\usetikzlibrary{positioning,calc}
\usetikzlibrary{graphs,graphs.standard,quotes}

\usepackage{pgfplots}
\pgfplotsset{compat=1.17} 
\usepackage[utf8]{inputenc}
\usepackage[square,sort,comma,numbers]{natbib}
\usepackage{graphicx}
\graphicspath{ {./images/} }
\usepackage{amsmath}
\usepackage{amsfonts}
\usepackage[linesnumbered,ruled]{algorithm2e}
\usepackage{mathtools}
\usepackage{mathrsfs}
\usepackage{tabularx}
\usepackage{multirow}
\usepackage{float}
\usepackage{subcaption}
\usepackage[font={small}]{caption}
\usepackage{stackengine}
\usepackage{url}
\usepackage{booktabs}
\usepackage{makecell}

\newtheorem{remark}{Remark}

\newtheorem{definition}{Definition}[section]

\newcommand{\abs}[1]{\left|#1\right|}
\newcommand{\norm}[1]{\left\|#1\right\|}

\newcommand{\bref}[1]{(\ref{#1})}

\DeclareMathOperator*{\argmin}{argmin}

\definecolor{ultramarine}{RGB}{0,0,0}
\definecolor{reviewV2}{RGB}{0,0,0}

\title{Accelerating Algebraic Multigrid Methods via Artificial Neural Networks\thanks{P.F.A and L.D. are members of the INdAM Research group GNCS. P.F.A has been partially funded by the research projects PRIN17 (n. 201744KLJL) and PRIN 2020 (n. 20204LN5N5), funded by Italian Ministry of University and Research (MUR). L.D. has been partially funded by the research project PRIN 2020 (n. 20204LN5N5) funded by MUR.}}
\author{Paola F. Antonietti, 
Matteo Caldana,
Luca Dede'$^{\dag}$\\[0.3cm]
MOX, Dipartimento di Matematica, 
Politecnico di Milano,\\
Piazza Leonardo da Vinci 32,
20133 Milano, Italy \\[0.2cm]
$^{\dag}$ {\small Corresponding author: {\tt luca.dede@polimi.it}
}}

\date{\today}

\begin{document}

\maketitle
\section*{Abstract}
We present a novel deep learning-based algorithm to accelerate -- through the use of Artificial Neural Networks (ANNs) -- the convergence of Algebraic Multigrid (AMG) methods for the iterative solution of the linear systems of equations stemming from finite element discretizations of \textcolor{reviewV2}{P}artial \textcolor{reviewV2}{D}ifferential \textcolor{reviewV2}{E}quations \textcolor{reviewV2}{(PDE)}. 
We show that ANNs can be be successfully used to predict the {strong connection} parameter that enters in the construction of the sequence of increasingly smaller matrix problems standing at the basis of the AMG algorithm, so as to maximize the corresponding convergence factor of the AMG scheme. To demonstrate the practical capabilities of the proposed algorithm, which we call AMG-ANN, we consider the iterative solution of the algebraic system of equations stemming from finite element discretizations of \textcolor{ultramarine}{two-dimensional model problems. First, we consider \textcolor{reviewV2}{an elliptic} equation with a highly heterogeneous diffusion coefficient and then a stationary Stokes problem.}
We train (off-line) our ANN with a rich dataset and present an in-depth analysis of the effects of tuning the strong threshold parameter on the convergence factor of the resulting AMG iterative scheme.


\section{Introduction}  
In the last thirty years, there has been an increasing demand for computationally efficient methods to solve sparse linear system of equations stemming from numerical discretization of Partial Differential Equations (PDEs). 
For real-life problems, the typical size of the resulting algebraic systems makes direct or classical one-level methods impractical and hierarchical iterative solvers have been intensively developed and studied.
This paper focuses on the Algebraic Multigrid (AMG) method (\cite{xuzikatanov2017}) for the iterative solution of the symmetric and positive definite systems of equations stemming from Finite Element (FE) approximation (\cite{hughes2012finite,Quarteroni_2017,QV94}) of elliptic \textcolor{ultramarine}{PDEs}. One of the main feature of AMG is that it is a purely matrix-based approach, thus it does not make use of any geometric information and the hierarchy of operators is constructed directly from the system matrix, provided that the underlying matrix has certain properties, see e.g., \cite{Brandt,ruge,STUBEN2001281,xuzikatanov2017}. AMG methods can be advantageous whenever geometric multigrid is not a viable option, e.g., whenever the sequence of coarser meshes at the basis of geometric multigrid is not available. AMG and AMG-like approaches have been developed to solve a variety of problems in the context of PDEs-based simulations; here we mention, for example, the AMG method based on element interpolation (AMGe) for solving  the discrete equations that arise in Ritz-type finite element methods, \cite{brezinafalgout2002AMGe, hensonvassilevski2001}, Maxwell’s equations \cite{joneslee2006}, linear elasticity \cite{baker2010improving}, Navier-Stokes's equations \cite{webster1994algebraic} and multi-phase porous media \cite{buiwnagoseikuffor2018}. In  \textcolor{ultramarine}{\cite{bank2015algebraic, cleary2000robustness, baker2012scaling, li2020efficient}},  AMG methods for large-scale supercomputing architectures are presented. In the paper \cite{xuzikatanov2017} by Xu and Zikatanov, AMG methods are presented and analyzed in \textcolor{ultramarine}{a} unified framework and an abstract theory for the construction of optimal coarse space as well as quasi-optimal spaces is derived. The abstract framework of \cite{xuzikatanov2017} covers most of the existing AMG methods, such as classical and energy-minimization AMG, unsmoothed and smoothed aggregation AMG, and spectral AMGe \cite{falgoutvassilevski2004}.
AMG methods for non-standard FE approximations have been also developed, for example in the context of  \textcolor{ultramarine}{discontinuous Galerkin} methods \cite{AntoniettiMelas_2020,BastianBlattScheichl_2012,SiefertTuminaro_et_al_2014}.

The AMG method relies on a set of parameters that defines how to algebraically carry out the coarsening phase. Often their tuning is based on experience and it could be rather inefficient in certain situations. In this paper, we propose using Machine  \textcolor{ultramarine}{Learning (ML)} and Deep Learning  \textcolor{ultramarine}{(DL)} algorithms to make the choice of the AMG parameters fully automatic so as \textcolor{ultramarine}{to} improve the efficiency of the method. The approach that we propose is based on the use of Artificial Neural Networks (ANNs). Artificial Neural Networks are \textcolor{ultramarine}{ML} and  \textcolor{ultramarine}{DL} models that are nowadays widely used in several problems in image recognition, speech recognition, and natural language processing \cite{goodfellow2016deep}. The introduction of  \textcolor{ultramarine}{C}onvolutional  \textcolor{ultramarine}{N}eural  \textcolor{ultramarine}{N}etworks (CNNs) \cite{lecun1998gradient} changed modern object recognition process \cite{he2015delving,krizhevsky2012imagenet,russakovsky2015imagenet}. Today, the most advanced ANNs in image recognition are variations of CNNs: ResNet \cite{he2016deep} and SENet \cite{hu2018squeeze}.

Nowadays, \textcolor{ultramarine}{ML} and \textcolor{ultramarine}{DL} models are increasingly being used in \textcolor{ultramarine}{s}cientific \textcolor{ultramarine}{c}omputing \cite{Neittaanmaki202227}, especially for the numerical approximation of ODEs and PDEs \cite{Mishra}. For example, \textcolor{ultramarine}{p}hysics\textcolor{ultramarine}{-i}nformed \textcolor{ultramarine}{n}eural \textcolor{ultramarine}{n}etworks have been introduced to approximate the solution of PDEs as a meshless method \cite{raissi2017machine,raissi2019physics} and ANNs are employed for model order reduction of parameter-dependent PDEs \cite{fresca2021comprehensive, hesthaven2018non, regazzoni2019machine}. ANNs can also be employed to enhance the performance of algorithms and solvers used in ``classical" numerical methods for the approximation of PDEs, i.e. as accelerators for scientific computing. In this context,
we mention for example: 
the enhancement of numerical stabilization methods for the FE approximation of advection-dominated differential problems, e.g. in \cite{janssensadvancing,tassi2021amachine};
the use of ANN to optimally select artificial viscosity for \textcolor{ultramarine}{discontinuous} Galerkin methods in
\cite{discacciati2020controlling}; exploiting CNN for grid refinement in \textcolor{ultramarine}{discontinuous} Galerkin and virtual element methods in
\cite{antoniettimanuzzi, antoniettimanuzzi2};
the hybrid ML-FETI-DP algorithm combining the advantages of adaptive coarse spaces in domain decomposition methods and certain supervised \textcolor{ultramarine}{ML} techniques that have been proposed in 
\cite{Heinlein_et_al_2021}.
\textcolor{ultramarine}{Moreover, ML techniques and ANNs have already been used to optimize multigrid algorithms, see \cite{KATRUTSA2020112524,Gottschalk2021,greenfeld2019learning}.}

In this work, we make use of ANNs to improve the tuning of the strong threshold parameter that enters in the definition of AMG so as to improve its performance. In order to test the proposed approach, we consider a two-dimensional elliptic equation with a highly heterogeneous diffusion coefficient discretized by the FE method. In order to use the sparsity pattern of the underlying matrix as input of the neural network, we introduce a pooling operator. We show how an ANN-enhanced approach can effectively improve the AMG performance.  The performance of the AMG method is measured in two ways: using the approximate convergence factor and using the elapsed time. We show that these two measures are strictly correlated, this entails that we have a unique way of measuring the performance. We demonstrate that, in some test cases, the value of the strong threshold parameter commonly used in literature can be improved so as to gain efficiency with respect to both measures. In particular, we test different models to tune the hyper-parameters of the model and we report the predictions of the models with the lowest loss.
\textcolor{reviewV2}{Our model shows very accurate predictions in case of a diffusion coefficient that exibits ``simple'' patterns. More complex coefficient distributions seem to work well only if the dataset reasonably rich, thus more computational work is needed.}

The paper is structured as follows. In Section~\ref{sec:amg} we recall the basic elements of the AMG methods. In Section~\ref{sec:model_prob} we introduce the model problem and its FE discretization. In Section~\ref{sec:ann} we give a brief overview on ANNs. The results of the numerical experiments are showcased in Section~\ref{sec:numeric_res}, \textcolor{ultramarine}{namely,} we report a wide set of numerical experiments aimed at testing the \textcolor{reviewV2}{algorithm's} performance when we \textcolor{ultramarine}{vary the strong threshold parameter. Then\textcolor{reviewV2}{,} in Section~\ref{sec:ann-enhancedAMG},} we design the architecture of \textcolor{reviewV2}{the} net, we introduce the pooling operator and test the model. Finally, in Section~\ref{sec:conc} we draw some conclusions.

\section{Algebraic Multigrid Methods}\label{sec:amg}
In this section, we introduce the main ingredients of AMG methods; we refer the reader to \cite{xuzikatanov2017} for a comprehensive description. We consider the linear system of equations:
\begin{equation}
\mathrm{A}_h \mathbf{u}_h = \mathbf{f}_h
\label{eq:lin_sys_base}
\end{equation}
where, for $n \in \mathbb{N}$, $\mathrm{A}_h \in\mathbb{R}^{n\times n}$ is symmetric \textcolor{reviewV2}{and} positive definite\textcolor{reviewV2}{. Here,} and $\mathbf{u}_h$, $\mathbf{f}_h\in\mathbb{R}^{n}$. Let $\mathcal{N}_h = \{ 1,...,n \}$ be the set with the indexes of all the variables. The set $\mathcal{N}_h$ is split into two disjoint subsets $\mathcal{C}_h$ and $\mathcal{F}_h$ such that $\mathcal{N}_h=\mathcal{C}_h \cup \mathcal{F}_h$ and $\mathcal{C}_h \cap \mathcal{F}_h = \emptyset$.

Let $\mathrm{I}_H^h:\mathbb{R}^{n_H} \to \mathbb{R}^{n}$ be the interpolation operator that maps coarse level vectors into fine level vectors, and let $\mathrm{I}^H_h:\mathbb{R}^{n} \to \mathbb{R}^{n_H}$ be the restriction operator that maps fine level vectors into coarse level vectors. It is assumed that $\mathrm{I}_H^h$ can be written as:
\begin{equation}
(\mathrm{I}_H^h \mathbf{e}_H)_i = 
\left\{
\begin{matrix}
    (\mathbf{e}_H)_i & \text{ if } i \in \mathcal{C}_h, \\
    \sum_{k \in \mathcal{P}_i} w_{ij}^h(\mathbf{e}_H)_j & \text{ if } i \in \mathcal{F}_h,
\end{matrix}
\right.
\label{eq:interp_generic_form}
\end{equation}
where $\mathbf{e}_H \in \mathbb{R}^{n_H}$ is a generic vector, $\mathcal{P}_i \subset \mathcal{C}_h$, for all $i \in \mathcal{F}_h$ is called a set of interpolatory variables for $i$ and $w_{ij}^h$ is a set of weights.
\textcolor{ultramarine}{One way to define $w_{ij}^h$ is the following. \textcolor{reviewV2}{We define the direct neighborhood of a point $i$ as} 
$$
\text{Neigh}(i)=\{j \neq i: (\mathrm{A}_h)_{ij} \neq 0\}.
$$
For $a \in \mathbb R$, \textcolor{reviewV2}{we define its positive and negative part as} $a^+=\max \{0,a\}$ and $a^-=\min \{0,-a\}$\textcolor{reviewV2}{, respectively}. \textcolor{reviewV2}{Similarly, we split $\mathcal P_i$ into two sets} 
$$
\mathcal{P}_i^+= \mathcal{P}_i \cap \{j \neq i: (\mathrm{A}_h)_{ij} > 0\},\quad
\mathcal{P}_i^-= \mathcal{P}_i \cap \{j \neq i: (\mathrm{A}_h)_{ij} < 0\}.
$$
Moreover, \textcolor{reviewV2}{\cite{trottenberg2001multigrid} shows that the following identity holds
\begin{equation}
 (\mathrm A)_{ii}(\mathbf e_h)_i + \alpha_i \sum_{k \in \mathcal P_i} (\mathrm A)_{ik}^-(\mathbf e_h)_k  + \beta_i \sum_{k \in \mathcal P_i} (\mathrm A)_{ik}^+(\mathbf e_h)_k = 0,
 \label{eq:amg_interpolation}
\end{equation}
where the coefficients are given by}
$$\alpha_i = \frac{\sum_{j\in \text{Neigh}(i)} (\mathrm{A}_h)_{ij}^-}{\sum_{j\in \mathcal{P}_i} (\mathrm{A}_h)_{ij}^-}, \quad \beta_i = \frac{\sum_{j\in \text{Neigh}(i)} (\mathrm{A}_h)_{ij}^+}{\sum_{j\in \mathcal{P}_i} (\mathrm{A}_h)_{ij}^+}.$$
Then, \textcolor{reviewV2}{the weights are defined as}
\begin{equation*}
w_{ik}^h = 
\begin{dcases}
-\alpha_i (\mathrm{A}_h)_{ik}/(\mathrm{A}_h)_{ii} \quad& k\in \mathcal{P}_i^-, \\
-\beta_i (\mathrm{A}_h)_{ik}/(\mathrm{A}_h)_{ii} & k\in \mathcal{P}_i^+, \\
0 & \text{otherwise}.
\end{dcases}
\end{equation*}
} \textcolor{reviewV2}{S}ince $\mathrm{A}_h$ is symmetric it is also assumed that:
\begin{equation}
\label{eq:operator_adj}
\mathrm{I}_H^h = (\mathrm{I}_h^H)^\top.
\end{equation}
Then, the coarse-level AMG matrix is defined as
$\mathrm{A}_H = \mathrm{I}_h^H \mathrm{A}_h \mathrm{I}_H^h \in \mathbb R^{n_H \times n_H}$. One of the key ingredients of the AMG method consists in the \textcolor{reviewV2}{definition} of the interpolation operator $I_H^h$ previously described. The classical coarsening algorithm prescribes to maintain at the coarse level all the strong connections that are defined through a parameter $\theta$, called the strong threshold parameter. Its rigorous definition is given in the following~\cite{Brannick}.

\begin{definition}
\label{def:strong_conn}
\textcolor{reviewV2}{Let $\mathrm A_h \in \mathbb R^{n\times n, }$.} Given a threshold \textcolor{reviewV2}{parameter} $0 < \theta \leq 1$, the variable $i$ strongly depends on the variable $j$ if
\begin{equation*}
- (\mathrm{A}_h)_{ij} \geq \theta \, \max_{k \neq i} \, \{ - (\mathrm{A}_h)_{ik} \}\textcolor{reviewV2}{, \qquad i,j=1, ..., n}.
\end{equation*}
\end{definition}

As a matter of fact, performing the $\mathcal C_h/\mathcal F_h$ splitting and \textcolor{reviewV2}{defining} the operators $I_h^H$ and $I_H^h$ requires choosing such strong threshold parameter $\theta$. \textcolor{reviewV2}{Even if the weights $w_{ij}^h$ do not directly depend on $\theta$, t}he performance of the AMG method will depend on \textcolor{reviewV2}{the choice of the threshold parameter}, which is empirically made a priori. 
\color{ultramarine}{To show how the choice of $\theta$ enters in the construction of the operators $I_h^H$ and $I_H^h$, we briefly recall the coarsening algorithm of \cite{trottenberg2001multigrid}. We introduce two sets that exploit Definition \ref{def:strong_conn}. The first set contains all the indexes \textcolor{reviewV2}{j} that are strongly connected to \textcolor{reviewV2}{the} index $i$, i.e.\:
\begin{equation*}
\mathcal{S}_i = \{ j \in \text{Neigh}(i): i \text{ is strongly dependent on } j \}, \quad i = 1, ..., n.
\end{equation*} 
Next, given $\mathcal{S}_i$, we introduce $\mathcal{S}_i^\top$ as:
\begin{equation*}
\mathcal{S}_i^\top = \{ j \in\mathcal{N}_h: i \in \mathcal{S}_j \}.
\end{equation*} 
The coarsening procedure follows this \textcolor{reviewV2}{algorithm: 
\begin{enumerate}
    \item Initialize the set of undecided variables $\mathcal{U}_h \leftarrow \mathcal{N}_h$.
    \item Choose a variable $i \in \mathcal{U}_h$ such that $\textcolor{reviewV2}{\eta}(i) \geq \textcolor{reviewV2}{\eta}(k)$ $\forall k \in \mathcal{U}_h$, where 
    \begin{equation*}
\eta(k) = \abs{\mathcal{S}_k^\top \cap \mathcal{U}_h} + 2\abs{\mathcal{S}_k^\top \cap \mathcal{F}_h} \quad \forall \, k \in \mathcal{U}_h,
\end{equation*}
and $\abs{\cdot}$ denotes the cardinality.
    \item Move the index $i$ from the set $\mathcal{U}_h$ to the set $\mathcal{C}_h$. 
    \item Add all the variables $j \in \mathcal{S}_i^\top$ to the set $\mathcal{F}_h$, that is add to the set $\mathcal{F}_h$ all the variables $j$ that strongly depend on $i$.
\end{enumerate}
Steps (2--4) are repeated until all the variables are either in $\mathcal{C}_h$ or $\mathcal{F}_h$. The measure $\eta$ is needed to avoid a non-uniform distribution of the variables. In this way, at each iteration, the algorithm selects as $i$ ($\mathcal{C}_h$-variable) the index such that the majority of $\mathcal{F}_h$-variables strongly dependent on.
}}

\color{black}
The last ingredient needed to define the AMG methods is a smoothing operator. In general one iteration of the smoothing can be written as:
\begin{equation*}
\mathbf{u}_h^{(k+1)} = \mathrm{S}_h \mathbf{u}_h^{(k)} + \mathbf{g}_h,  \quad k \geq 0,
\end{equation*}
where $\mathrm{S}_h \in \mathbb R^{n \times n}$ denotes the smoothing operator \textcolor{reviewV2}{to be properly chosen}. \textcolor{reviewV2}{Equivalently, it can be written in preconditioned form as:}
\color{ultramarine}
\begin{equation}
\mathbf{u}_h^{(k+1)} = \mathbf{u}_h^{(k)} + \mathrm{B}_h(\mathbf{f}_h - \mathrm{A}_h\mathbf{u}_h^{(k)}), \quad k \geq 0,
\label{eq:precon_iter}
\end{equation}
where $\mathrm{S}_h = \mathrm{I}_h - \mathrm{B}_h\mathrm{A}_h$, $\mathbf{g}_h = \mathrm{B}_h \mathbf{f}_h$, $\mathrm{B}_h = (\mathrm{I}_h - \mathrm{S}_h)\mathrm{A}_h^{-1}$ \textcolor{reviewV2}{and $\mathrm{I}_h$ is the identity operator. In practice, either $(\mathrm S_h, \mathbf g_h)$ or $\mathrm B_h$ are given and uniquely identify the smoother. }
\color{black} 
In the following, the notation:
\begin{equation*}
\mathbf{u}_h^{(l)} = \texttt{smooth}^l(\mathrm{A}_h, \mathbf{u}_h^{(0)}, \mathbf{f}_h),
\end{equation*}
means that $\mathbf{u}_h^{(l)}$ is the result of $l$ steps of \bref{eq:precon_iter}, starting from an an initial vector $\mathbf{u}_h^{(0)}$. \textcolor{reviewV2}{I}n Algorithm~\ref{a:two_grid_iter} \textcolor{reviewV2}{we report} one iteration of the two-level algorithm, where $\nu_1$ and $\nu_2$ are the \textcolor{reviewV2}{numeber of} smoothing steps that we apply before and after the error correction, respectively. The complete two-level AMG algorithm is outlined in Algorithm~\ref{a:mgcycle}. As usual, in Algorithm~\ref{a:mgcycle}, $tol$ is a user-defined tolerance that is employed \textcolor{reviewV2}{as} a stopping criterion. Analogously, $N_{max}$ is the maximum number of iterations allowed. We notice that Algorithm~\ref{a:two_grid_iter} involves selecting the parameter $\theta$ a priori.

\begin{algorithm}[t]
    \SetKwFunction{smooth}{smooth}
    \SetKwFunction{solve}{solve}
    $\mathbf{u}_h^{(*)}$ $\leftarrow$ $\smooth^{\nu_1}(\mathrm{A}_h, \mathbf{u}_h^{(k)}, \mathbf{f}_h)$\;
    $\mathbf{r}_h$ $\leftarrow$ $\mathbf{f}_h - \mathrm{A}_h \mathbf{u}_h^{(*)}$\;
    $\mathbf{r}_H$ $\leftarrow$ $\mathrm{I}_h^H \mathbf{r}_h$\;
    $\mathbf{e}_H$ $\leftarrow$ $\solve(\mathrm{A}_H, \mathbf{r}_H)$\;
    $\mathbf{u}_h^{(*)}$ $\leftarrow$ $\mathbf{u}_h^{(*)} + \mathrm{I}_H^h \mathbf{e}_H$\;
    $\mathbf{u}_h^{(k+1)}$ $\leftarrow$ $\smooth^{\nu_2}(\mathrm{A}_h, \mathbf{u}_h^{(*)}, \mathbf{f}_h)$
    \caption{One Iteration of the two-level AMG method \newline $\mathbf{u}_h^{(k+1)} = \texttt{two\_level\_iteration}(\mathbf{u}_h^{(k)}, \mathrm{A}_h, \mathbf{f}_h, \nu_1, \nu_2, \mathrm{I}_h^H, \mathrm{I}_H^h)$}
    \label{a:two_grid_iter}
\end{algorithm}

\begin{algorithm}[t]
    perform the $\mathcal{C}_h/\mathcal{F}_h$-splitting {using $\theta$}\;
    build the operators $\mathrm{I}_h^H,\, \mathrm{I}_H^h$ {using the $\mathcal{C}_h/\mathcal{F}_h$-splitting} \;
    \While{$k < N_{max}$ \textnormal{\textbf{and}} $\norm{\mathrm{A}_h\mathbf{u}_h^{(k)}-\mathbf{f}_h}/\norm{\mathbf{f}_h} < tol$} {
         $\mathbf{u}_h^{(k+1)} \leftarrow \texttt{two\_level\_iteration}(\mathbf{u}_h^{(k)}, \mathrm{A}_h, \mathbf{f}_h, \nu_1, \nu_2, \mathrm{I}_h^H, \mathrm{I}_H^h)$
    }
    \caption{Two-Level AMG algorithm \newline $\mathbf{u}_h^{(k+1)} = \texttt{AMG}(\mathbf{u}_h^{(0)}, \mathrm{A}_h, \mathbf{f}_h, {\theta,} \nu_1, \nu_2, N_{max}, tol)$}
    \label{a:mgcycle}

\end{algorithm}

As a matter of fact, the two-level AMG Algorithm~\ref{a:mgcycle} can be immediately extended to many levels by simply calling recursively Algorithm~\ref{a:mgcycle} until a sufficiently coarse level is reached (where a direct solver is employed). For the sake of the analysis carried out in the present paper, we will focus on the two-level method.

\section{Model Problem\textcolor{reviewV2}{s}}\label{sec:model_prob}
Throughout this work, we use standard notation for Sobolev spaces \cite{lions72}. Let $\Omega$ be an open, bounded domain in $\mathbb R^2$ and let $\partial \Omega = \overline{\Gamma}_D$. The \textcolor{reviewV2}{first} model problem \textcolor{reviewV2}{we consider} reads:
\begin{equation}
\begin{cases}
-\text{div}(\mu(x,y)\nabla u) = f, \quad &\text{in} \: \Omega,\\
u = g_D, &\text{on} \: \Gamma_D,
\end{cases}
\label{eq:2prob1}
\end{equation}
where $f \in L^2(\Omega)$ is a given forcing term, and $g_D \in H^{1/2}(\Gamma_D)$ is the given Dirichlet boundary data. The function $\mu \in L^{\infty}(\Omega)$ is a positive diffusion coefficient. In this work it will be a piece-wise non-negative constant function. To handle non homogeneous Dirichlet boundary condition we define $\tilde{u}$ by the means of the lifting $\tilde{u}=u-\tilde{g}$, where $\tilde{g}$ is an extension of $g_D$ in $H^1(\Omega)$.
The weak formulation of problem~\bref{eq:2prob1} reads:
\begin{equation}
\text{find } \tilde{u} \in H^1_{\Gamma_D}(\Omega) \ : \ (\tilde{u},v) = F(v) \quad \forall v \in H^1_{\Gamma_D}(\Omega),
\label{eq:weak_form_2222}
\end{equation}
where $H^1_{\Gamma_D}(\Omega) := \{ v \in H^1(\Omega): v|_{\Gamma_D} = 0 \}$ and
\begin{equation}
a(\tilde{u},v) = \int_{\Omega} \mu \nabla \tilde{u} \cdot \nabla v \: \text{d}\Omega, \quad
F(v) = \int_{\Omega} fv \: \text{d}\Omega - \int_{\Omega} \mu \nabla \tilde{g} \cdot \nabla v \: \text{d}\Omega.
\label{eq:aFdef}
\end{equation}
The well-poseness of problem~\eqref{eq:weak_form_2222} is given by the Lax-Milgram's theorem~\cite{brezis2010functional}. 

Now we pass to the FE formulation. We consider a quasi uniform mesh $\mathcal{T}_h$ of $\Omega$. We denote with the parameter $h > 0$ the mesh size of $\mathcal{T}_h$ given by $h = \max_{T \in \mathcal{T}_h} h_T$, where $h_T$ is the diameter of the element $T \in \mathcal{T}_h$. In our case, we use:
\begin{equation*}
V_h = \{ v_h \in X_h^1: v_h = 0 \ \text{on} \; \Gamma_D \}, 
\end{equation*}
where \textcolor{reviewV2}{$X_h^{1} = \{v_h \in C^0(\bar{\Omega}): v_h|_T\circ F_T\in\mathbb Q_1(\hat\Omega) \, \forall \, T \in \mathcal T_h\}$, $F_T:\hat\Omega\rightarrow T$ is an invertible function that maps the reference square $\hat\Omega=(-1,1)^2$ to the mesh element $T$, and $\mathbb Q_N$ is the space of polynomials with real coefficients and degree less than or equal to $N$ in each coordinate direction.}
The finite dimensional formulation of \bref{eq:weak_form_2222} reads:
\begin{equation}
\text{find } \tilde{u}_h \in V_{h} \: \text{s.t.:} \; a(\tilde{u}_h,v_h) = F(v_h) \quad \forall v_h \in V_{h}.
\label{eq:2feprob}
\end{equation}
By setting $n = \text{dim}(V_h)$, we denote with $\{\phi_1, ..., \phi_n\}$ the FE basis for $V_h$. Then, from Eq.~\eqref{eq:2feprob}, we obtain the linear system of equations
$A_h\mathbf u_h = \mathbf f$,
where:
\begin{equation}
a(\phi_j,\phi_i) = (\mathrm{A}_h)_{ij}, \ F(\phi_i) = (\mathbf{f})_i, \ (\mathbf{u}_h)_i=\tilde{u}_i.
\label{eq:2linear_sys}
\end{equation}

\color{reviewV2}As a second model problem we consider the Stokes equations. \color{ultramarine} Namely, we are looking for a velocity $\mathbf{u} : \mathbb{R}^2 \rightarrow \mathbb{R}^2$ and pressure $p: \mathbb{R}^2 \rightarrow \mathbb{R}$ that satisfy the Stokes equation, which reads

\begin{equation}
\begin{cases}
- \nu \triangle \mathbf{u} + \nabla p = \mathbf{0} &\text{in} \: \Omega,\\
- \text{div} \, \mathbf{u} = 0 &\text{in} \: \Omega,\\
\mathbf{u} = \mathbf{0} &\text{on} \: \Gamma_0,\\
\mathbf{u} = \mathbf u_0 &\text{on} \: \Gamma_{in}\\
\nu \partial_{\mathbf{n}}\mathbf{u} - p\mathbf{n}=0 &\text{on} \: \Gamma_{out},\\
\end{cases}
\label{eq:stokes}
\end{equation}
where $\mathbf{n}$ denotes the outer normal vector and $\mathbf u_0$ is the parabolic inflow velocity, with maximum $U > 0$.
\color{reviewV2} Here, we decompose the boundary $\Gamma$ as $\Gamma = \Gamma_{in} \cup \Gamma_{out} \cup \Gamma_0$, where $\Gamma_{in}, \Gamma_{out}, \Gamma_0$ are disjoint open sets with positive measure. To guarantee the well-posedness of the problem, we prescribe that $p \in L^2_0(\Omega)$ i.e.\ is a $L^2(\Omega)$ functions with zero average. We introduce the functional spaces:
$$V = \{v \in [H^1(\Omega)]^2: v|_{\partial\Omega} = 0\}, \quad Q = L^2_0(\Omega),$$
and endow them with the norms
$\norm{v}_V = \norm{\nu^{1/2} \nabla v}_{L^2(\Omega)}, \norm{q}_Q = \norm{q}_{L^2(\Omega)}.$
The weak formulation of problem \bref{eq:stokes} reads: find $(u, p) \in V \times Q$, such that
$$a(u, v) + b(p, v) - b(q, u) = 0 \quad \forall (v, q) \in V \times Q$$
where
$$a : V \times V \rightarrow \mathbb{R}, \quad a(u, v) = \int_\Omega \nu \nabla u : \nabla v,$$
$$b : Q \times V \rightarrow \mathbb{R}, \quad b(p, v) = -\int_\Omega p \, \text{div} v.$$
It is well-known that the bilinear form $b(\cdot, \cdot)$ satisfies a continuous inf-sup condition;
see, e.g., \cite{boffi2013mixed}. 
\color{ultramarine} We introduce a uniform quadrilateral mesh $\mathcal{T}_h$ of $\Omega$. Discretizing using the standard polynomial spaces \textcolor{reviewV2}{$V_h=\left[X_h^2\right], Q_h = X_h^1$} on $\mathcal{T}_h$ we obtain the following algebraic formulation
\begin{eqnarray*} \left(\begin{array}{cc} \mathrm{A}_h & \mathrm{B}_h^\top \\ \mathrm{B}_h & 0 \end{array}\right) \left(\begin{array}{c} \mathbf{u}_h \\ \mathbf{p}_h \end{array}\right) = \left(\begin{array}{c} \mathbf{0} \\ \mathbf{0}  \end{array}\right),
\label{eq:stokes_algebraic}
\end{eqnarray*}
\color{reviewV2}
where, setting $N_h = \textnormal{dim}(V_h)$ and $M_h = \textnormal{dim}(Q_h)$, $\mathrm A_h \in \mathbb R ^ {N_h \times N_h}$ and $\mathrm B_h \in \mathbb R ^ {M_h \times N_h}$ are the matrix representation of the bilinear forms $a(\cdot, \cdot)$ and $b(\cdot, \cdot)$, respectively.
\color{black}

\section{Artificial Neural Networks}\label{sec:ann}
An artificial neural network, is a regression (or classification) model which \textcolor{reviewV2}{given by} a function $\mathscr{F} \, :\, \mathbb{R}^N \rightarrow \mathbb{R}^M$ \textcolor{reviewV2}{defined as}
\begin{equation}
   \mathscr{F}(\mathbf{x};{\boldsymbol{\gamma}})=\mathbf{y},
\label{eq:def_general_ffnn}
\end{equation}
where $\mathbf{x}$ is the input, $\mathbf{y}$ is the predicted value of the regression and ${\boldsymbol{\gamma}}$ is the vector containing all the parameters of the model. The function $\mathscr F$ is the composition of $K$ functions $\mathscr F^{(k)}$ called \textit{layers}, the number of layers $K$ is called \textit{depth} of the model. In the case of feed-forward neural network the layer is defined as
\begin{equation}
\begin{split}
&
\begin{dcases}
\mathbf{a}^{(k)} = \mathrm{W}^{(k)} \textcolor{reviewV2}{\mathbf{z}}^{(k-1)} + \mathbf{b}^{(k)} \\
\textcolor{reviewV2}{\mathbf{z}}^{(k)} = \mathscr H^{(k)}(\mathbf{a}^{(k)})
\end{dcases}
\quad \text{for } k = 1, ..., K \\
& \mathbf{x} = \textcolor{reviewV2}{\mathbf{z}}^{(0)}, \; \mathbf{y} = \textcolor{reviewV2}{\mathbf{z}}^{(K)}, \; N_0=N, \; N_K=M, \\
&\textcolor{reviewV2}{\boldsymbol \gamma = \{(\mathrm W^{(1)}, \mathbf b^{(1)}), ..., (\mathrm W^{(K)}, \mathbf b^{(K)})\}},
\end{split}
\label{eq:layer}
\end{equation}
where $\mathrm{W}^{(k)} \in \mathbb{R}^{N_k \times N_{k-1}}$ (\textit{weights}) and $\mathbf{b}^{(k)} \in \mathbb{R}^{N_k}$ (\textit{biases}) are the parameters $\boldsymbol{\gamma}$, and $\mathscr H^{(k)}(\cdot)$ is a scalar non-linear almost everywhere differentiable function that is applied component-wise to $\mathbf{a}^{(k)}$ and called \textit{activation function}. The Rectified Linear Unit ReLU$(x)=\max\{0,x\}$ is our choice of activation function $\mathscr H^{(k)}(\cdot)$ since, in recent years, it has became very popular due to the fact that it greatly improves the convergence of the stochastic gradient descent algorithm compared to the sigmoid/tanh functions \cite{krizhevsky2012imagenet}. Moreover, it features lighter computations with a random initialization network as only about half of hidden units have a non-zero output and faster evaluation with respect to the sigmoid/tanh functions.
\textcolor{ultramarine}{Indeed, in our experiments, employing the ReLU activation function seems to lead to better results with respect to the $\tanh$ activation function.}

Next, We define the loss function $\mathcal L$.
We assume that a dataset composed by $P$ couples $(\mathbf{x}^{(i)},\mathbf{y}^{(i)})$ is available; these are realizations of random variables $\mathscr{X}$, $\mathscr{Y}$.
Once defined the ANN architecture, its training boils down to minimize the average training error, namely
\begin{equation}
J({\boldsymbol{\gamma}}) = \frac{1}{P} \sum_{i=1}^{P} \mathcal L(\mathscr F(\mathbf{x}^{(i)};{\boldsymbol{\gamma}}),\mathbf{y}^{(i)}).
\label{eq:cost_fun_J}
\end{equation}
A typical choice of the loss function $\mathcal L$ that we also use in this paper are the \textcolor{ultramarine}{M}ean \textcolor{ultramarine}{S}quare \textcolor{ultramarine}{E}rror (MSE) \textcolor{ultramarine}{and Mean Absolute Error (MAE)}.

For determining the parameters $\boldsymbol{\gamma}$, we use the Adaptive Moment Estimation (Adam) method \cite{kingma2014adam}. It is a variant of the stochastic gradient descent method that combines the Root Mean Squared propagation (RMSProp) algorithm \cite{tieleman2012lecture} and momentum method \cite{sutskever2013importance} few other significant modifications, namely the momentum is recorded in the history of the gradient and there is a correction term of the bias for the estimation of the first and second order moments of the gradient.

Finally, to prevent overfitting and minimize the generalization error we employ four regularization techniques\textcolor{ultramarine}{. Namely, we will always employ} \textcolor{reviewV2}{an} early stopping criterion, namely we stop the training at the point of smallest error with respect to the validation dataset and random parameter initialization \textcolor{ultramarine}{\cite{he2015delving}. Moreover, we will test} dropout\textcolor{reviewV2}{, which consists in randomly omitting the weights and biases of some neurons $(\mathbf z^{(k)})_l$ during the training process} \cite{srivastava2014dropout}, and batch normalization\textcolor{reviewV2}{, a transformation applied at the end of a layer that normalizes its output by the empirical mean and variance of the minibatch} \cite{ioffe2015batch}. \textcolor{ultramarine}{As we will show in the forthcoming sections, the latter seems not to lead to substantial improvements in our model.}

As we want to use the matrix of the linear system $\mathrm A_h$ as input of the network we employ \textcolor{ultramarine}{CNN}. Their characteristic is that the layer takes the form of a cross-convolution between the input and a matrix $\mathrm{K} \in \mathbb{R}^{ D \times D }$, called \textit{kernel}. Three other hyper-parameters control how the convolution is performed: number of filters, stride and zero-padding size. Moreover, in the last stage of the layer a pooling function is applied. The pooling function is a form of down-sampling that replaces the output of the net at a certain location with a summary statistic of the nearby outputs. The aim of the pooling operation is to control the number of parameters and limit the overfitting. We refer the reader to \cite{goodfellow2016deep} for more details.

\section{Numerical Assessment of the \textcolor{reviewV2}{Dependence of the Performance of the AMG Method on the }Strong Threshold Parameter}\label{sec:numeric_res}
\textcolor{ultramarine}{In this section we assess the relation between the choice of the strong threshold parameter $\theta$ and the corresponding performance of the AMG method.}

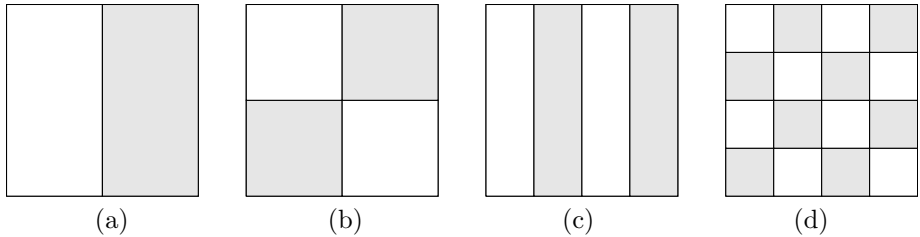
\begin{figure}[t]
\centering
\resizebox{1.0\textwidth}{!}{
%
%
\begin{tikzpicture}

\begin{axis}[%
width=4.521in,
height=0.952in,
at={(0.758in,1.788in)},
axis line style=white,
scale only axis,
xmin=-1,
xmax=8.5,
xtick={\empty},
ymin=-1,
ymax=1,
ytick={\empty},
axis background/.style={fill=white},
axis x line*=bottom,
axis y line*=left,
legend style={legend cell align=left, align=left, draw=white!15!black}
]
\addplot [color=black]
  table[row sep=crcr]{%
-1	-1\\
1	-1\\
1	1\\
-1	1\\
-1	-1\\
};

\addplot[area legend, draw=black, fill=white]
table[row sep=crcr] {%
x	y\\
-1	-1\\
-1	1\\
0	1\\
0	-1\\
}--cycle;

\addplot[area legend, draw=black, fill=black!10]
table[row sep=crcr] {%
x	y\\
0	-1\\
0	1\\
1	1\\
1	-1\\
}--cycle;

\addplot [color=black]
  table[row sep=crcr]{%
1.5	-1\\
3.5	-1\\
3.5	1\\
1.5	1\\
1.5	-1\\
};

\addplot[area legend, draw=black, fill=black!10]
table[row sep=crcr] {%
x	y\\
1.5	-1\\
1.5	0\\
2.5	0\\
2.5	-1\\
}--cycle;

\addplot[area legend, draw=black, fill=white]
table[row sep=crcr] {%
x	y\\
1.5	0\\
1.5	1\\
2.5	1\\
2.5	0\\
}--cycle;

\addplot[area legend, draw=black, fill=white]
table[row sep=crcr] {%
x	y\\
2.5	-1\\
2.5	0\\
3.5	0\\
3.5	-1\\
}--cycle;

\addplot[area legend, draw=black, fill=black!10]
table[row sep=crcr] {%
x	y\\
2.5	0\\
2.5	1\\
3.5	1\\
3.5	0\\
}--cycle;

\addplot [color=black]
  table[row sep=crcr]{%
4	-1\\
6	-1\\
6	1\\
4	1\\
4	-1\\
};

\addplot[area legend, draw=black, fill=white]
table[row sep=crcr] {%
x	y\\
4	-1\\
4	1\\
4.5	1\\
4.5	-1\\
}--cycle;

\addplot[area legend, draw=black, fill=black!10]
table[row sep=crcr] {%
x	y\\
4.5	-1\\
4.5	1\\
5	1\\
5	-1\\
}--cycle;

\addplot[area legend, draw=black, fill=white]
table[row sep=crcr] {%
x	y\\
5	-1\\
5	1\\
5.5	1\\
5.5	-1\\
}--cycle;

\addplot[area legend, draw=black, fill=black!10]
table[row sep=crcr] {%
x	y\\
5.5	-1\\
5.5	1\\
6	1\\
6	-1\\
}--cycle;

\addplot [color=black]
  table[row sep=crcr]{%
6.5	-1\\
8.5	-1\\
8.5	1\\
6.5	1\\
6.5	-1\\
};

\addplot[area legend, draw=black, fill=black!10]
table[row sep=crcr] {%
x	y\\
6.5	-1\\
6.5	-0.5\\
7	-0.5\\
7	-1\\
}--cycle;

\addplot[area legend, draw=black, fill=white]
table[row sep=crcr] {%
x	y\\
6.5	-0.5\\
6.5	0\\
7	0\\
7	-0.5\\
}--cycle;

\addplot[area legend, draw=black, fill=black!10]
table[row sep=crcr] {%
x	y\\
6.5	0\\
6.5	0.5\\
7	0.5\\
7	0\\
}--cycle;

\addplot[area legend, draw=black, fill=white]
table[row sep=crcr] {%
x	y\\
6.5	0.5\\
6.5	1\\
7	1\\
7	0.5\\
}--cycle;

\addplot[area legend, draw=black, fill=white]
table[row sep=crcr] {%
x	y\\
7	-1\\
7	-0.5\\
7.5	-0.5\\
7.5	-1\\
}--cycle;

\addplot[area legend, draw=black, fill=black!10]
table[row sep=crcr] {%
x	y\\
7	-0.5\\
7	0\\
7.5	0\\
7.5	-0.5\\
}--cycle;

\addplot[area legend, draw=black, fill=white]
table[row sep=crcr] {%
x	y\\
7	0\\
7	0.5\\
7.5	0.5\\
7.5	0\\
}--cycle;

\addplot[area legend, draw=black, fill=black!10]
table[row sep=crcr] {%
x	y\\
7	0.5\\
7	1\\
7.5	1\\
7.5	0.5\\
}--cycle;

\addplot[area legend, draw=black, fill=black!10]
table[row sep=crcr] {%
x	y\\
7.5	-1\\
7.5	-0.5\\
8	-0.5\\
8	-1\\
}--cycle;

\addplot[area legend, draw=black, fill=white]
table[row sep=crcr] {%
x	y\\
7.5	-0.5\\
7.5	0\\
8	0\\
8	-0.5\\
}--cycle;

\addplot[area legend, draw=black, fill=black!10]
table[row sep=crcr] {%
x	y\\
7.5	0\\
7.5	0.5\\
8	0.5\\
8	0\\
}--cycle;

\addplot[area legend, draw=black, fill=white]
table[row sep=crcr] {%
x	y\\
7.5	0.5\\
7.5	1\\
8	1\\
8	0.5\\
}--cycle;

\addplot[area legend, draw=black, fill=white]
table[row sep=crcr] {%
x	y\\
8	-1\\
8	-0.5\\
8.5	-0.5\\
8.5	-1\\
}--cycle;

\addplot[area legend, draw=black, fill=black!10]
table[row sep=crcr] {%
x	y\\
8	-0.5\\
8	0\\
8.5	0\\
8.5	-0.5\\
}--cycle;

\addplot[area legend, draw=black, fill=white]
table[row sep=crcr] {%
x	y\\
8	0\\
8	0.5\\
8.5	0.5\\
8.5	0\\
}--cycle;

\addplot[area legend, draw=black, fill=black!10]
table[row sep=crcr] {%
x	y\\
8	0.5\\
8	1\\
8.5	1\\
8.5	0.5\\
}--cycle;

\end{axis}
\end{tikzpicture}%
}
\begin{tabular}{cccc}

(a) \qquad \qquad & \qquad \qquad (b) \qquad \qquad & \qquad \qquad (c) \qquad \qquad & \qquad \qquad  (d)  \\
\end{tabular}
\caption{Four possible patterns of the diffusion coefficient $\mu$ of problem~\bref{eq:2prob1} \textcolor{ultramarine}{on $\Omega = \Omega_{gray} \cup \Omega_{white} = (-1, 1)^2$}: it is defined such that $\mu=1$ on the white tiles $\Omega_{white}$ and $\mu=10^\varepsilon$\textcolor{ultramarine}{, $\varepsilon > 0$,} on the gray ones $\Omega_{gray}$.}
\label{fig:mu_pattern}
\end{figure}

For our model problem \bref{eq:2prob1} we select the diffusion coefficient $\mu$ to be a piece-wise \textcolor{reviewV2}{positive} constant function. We assume that $\mu$ features different patterns, where the domain $\Omega$ splits in\textcolor{reviewV2}{to} strides or \textcolor{reviewV2}{has} a checkerboard \textcolor{reviewV2}{pattern}; see Figure~\ref{fig:mu_pattern}. The value of $\mu(x,y)$ depends on which ``tile'' \textcolor{reviewV2}{$(x, y)$} it belongs, namely
\begin{equation}
    \mu(x,y) = 
\begin{dcases}
    1 \quad &\text{ if } (x,y) \in \Omega_{gray},\\
    10^\varepsilon & \text{ if } (x,y) \in \Omega_{white},
\end{dcases}
\label{eq:mu_patt}
\end{equation}
where $\varepsilon$ is a parameter and $\Omega_{gray}$ and $\Omega_{white}$ are shown in Figure~\ref{fig:mu_pattern}. The experiments were carried out so that the exact solution $u$ of problem~\bref{eq:2prob1} is $u(x,y)=\cos(\pi x) \cos(\pi y)$ for patterns (a) and (b), while $u(x,y)=\cos(2 \pi x) \cos(2 \pi y)$ for patterns (c) and (d)\textcolor{reviewV2}{.} Dirichlet boundary conditions are set on the whole boundary $\partial \Omega$. Moreover, we employ regular cartesian meshes, so that the discontinuity of $\mu$ is aligned with mesh elements.

The implementation of the AMG method on which we rely on is the BoomerAMG of the library HYPRE \cite{falgout02hypre}. In particular, we use the AMG method as a preconditioner to accelerate the \textcolor{ultramarine}{C}onjugate \textcolor{ultramarine}{G}radient (CG) iterative method \cite{falgout96mgcg}. \textcolor{ultramarine}{The simulations were run using deal.II \cite{dealII91} with PETSc \cite{abhyankar2018petsc} on Ubuntu 18.04 LTS with CPU Intel i7-8550U. For sake of simplicity, the computations were carried out in serial. \textcolor{reviewV2}{However, even if the choice of $\theta$ might influence the parallelization, the same approach} could be extended also to the parallel case. }

To measure the performance of AMG we employ two performance indexes $p$: the elapsed CPU time and the approximate convergence factor $\rho$, defined as follows. Let $\rho^{(k)}$ be defined as
\begin{equation}
\rho^{(k)} = \left( \frac{\norm{\mathbf{r}^{(k)}}}{\norm{\mathbf{r}^{(0)}}} \right)^{\frac{1}{k}},
\label{eq:def_approx_conv_fact}
\end{equation}
where $\mathbf{r}^{(k)}$ is the residual at the $k$-th iteration and $\|\cdot\|$ is the standard euclidean norm. Then, we define $\rho$ as 
$$\rho=\rho^{(N_{it})},$$ 
where $N_{it}$ is the number of iterations reached to reduce the (relative) residual below the given tolerance of the linear solver (here it is equal to $N_{it} = \textnormal{min}_k \{ k \in \mathbb{N} \ \text{such that } \norm{\mathbf{r}^{(k)}} < 10^{-8} \} $).

\subsection{Relation between $\theta$ and the number of levels}\label{sec:theta_vs_lvls}
\begin{figure}[t]
\begin{center}
	\includegraphics[width=1.0\textwidth]{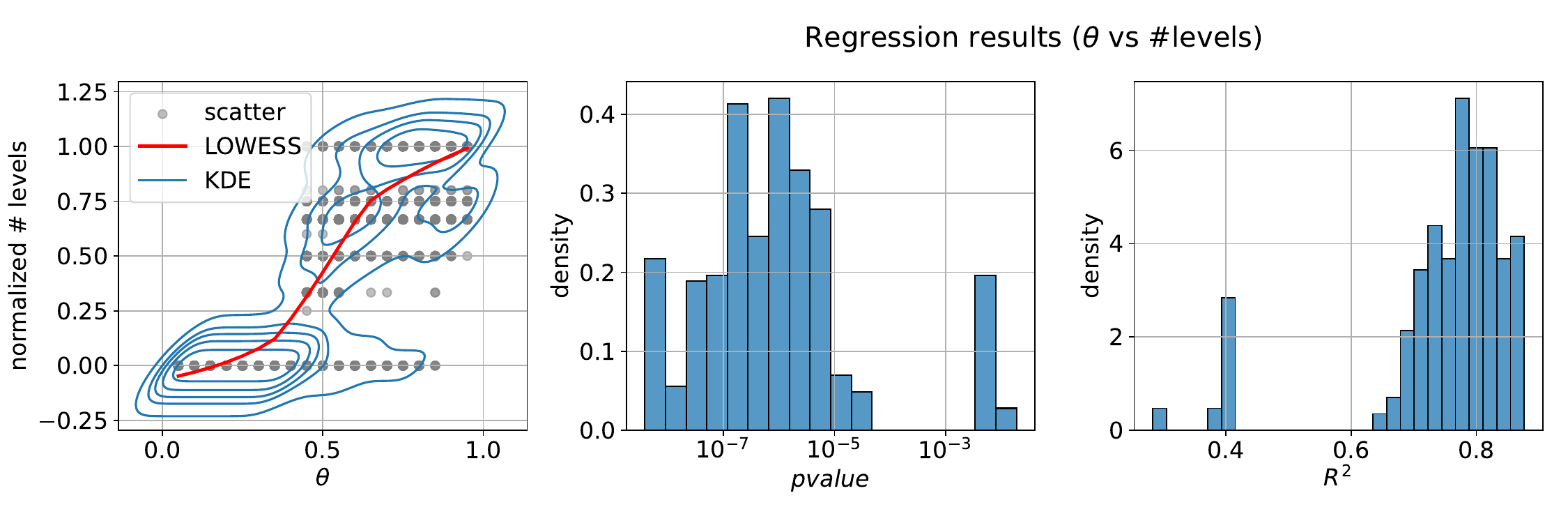}
\end{center}
\caption{\textcolor{ultramarine}{\textit{Left.} A scatter plot of the strong threshold parameter $\theta$ versus the min-max normalized number of levels of the AMG method. We superimposed a KDE and a LOWESS. \textit{Center and Right.} For each test case in the dataset we perform a least square analysis between $\theta$ and the \textcolor{reviewV2}{corresponding} number of levels of the AMG method. \textit{Center.} Histogram of the p-value of the least square analysis. \textit{Right.} Histogram of the coefficient of determination (R${}^2$). We refer to Table~\ref{tab:ls_def} for the definition of p-value and R${}^2$.}}
\label{fig:theta_vs_levels}
\end{figure}

\begin{table}[t]
\centering
 {\small\color{reviewV2}\begin{tabularx}{1\textwidth}{ llX }
  \hline
  Quantity & Definition & Formula \\ 
  \hline
$nn$ & Number of data points &\\ \addlinespace
SSE & Sum of squares of errors & ${\sum_{i=1}^{nn} (y^{(i)} - \hat x_1 x^{(i)} - \hat x_0)^2}$ 
\\ \addlinespace
TSS & Total sum of squares & $\sum_{i=1}^{nn} (y^{(i)} - \bar y)^2$, $\bar y = \frac{1}{nn}\sum_{i=1}^{nn}y^{(i)}$
\\ \addlinespace
SSR & Sum of squared residuals & $\sum_{i=1}^{nn} (\bar y - \hat x_1 x^{(i)} - \hat x_0)^2$  
\\ \addlinespace
R${}^2$        & Coefficient of determination & $1 - $ RSS / TSS \\ \addlinespace 
F-statistic    & F-statistic of the regression & $(nn - 2)$ SSR / SSE\\ \addlinespace
AIC            & Akaike's information criterion & $4 - 2\log(\hat L)$, where $\hat L$ is the log-likelihood of the model.\\ \addlinespace
SE$(\hat x_1)$& Standard error of $\hat x_1$ & $\sqrt{\textnormal{SSE}/\sum_{i=0}^{nn}(x^{(i)} - \bar x)^2/(nn-2)}$\\ \addlinespace
t-value $\hat x_1$ &t-value of $\hat x_1$& $\hat x_1 / \textnormal{SE}(\hat x_1)$\\ \addlinespace
p-value $\hat x_1$ & p-value of $\hat x_1$& $2\, \textnormal{cdf}_{t,1}(-\abs{tv_{\hat x_1}})$, where $\textnormal{cdf}_{t,1}$ is the cumulative density function of the Student's t distribution with one degree of freedom.\\
\hline
 \end{tabularx}
 }
\caption{\textcolor{reviewV2}{Definition of the quantities employed in the analysis of Section~\ref{sec:theta_vs_lvls} and Section~\ref{sec:machine_time}. The analysis aim at assessing the relation between the predictor scalar variable $x$ and the predicted scalar variable $y$ given $nn$ data points $(x^{(i)}, y^{(i)})$. We consider the linear model $y = \hat x_1 x + \hat x_0$. We refer to \cite{seabold2010statsmodels} for more details.}}
\label{tab:ls_def}
\end{table}
\color{ultramarine}{
In this section, we show how different choices of $\theta$ influence the number of levels built by the BoomerAMG algorithm and the corresponding size of the coarsest matrix $\mathrm A_h$.
Let us call \textit{test case} a fixed choice of the pattern of the diffusion coefficient $\mu$, the coefficient $\varepsilon$, and the size of the mesh $h$.
For each test case, we vary $\theta$ and record the number of levels used by the AMG algorithm to solve the problem.  \textcolor{reviewV2}{The scatter in Figure~\ref{fig:theta_vs_levels} (left) shows that} if the strong threshold is small, namely $0 < \theta < 0.3$, the number of levels is constant and \textcolor{reviewV2}{it is} equal to the minimum number of levels used to solve that test case. \textcolor{reviewV2}{In Figure~\ref{fig:theta_vs_levels} (left) we superimposed a Kernel Density Estimate (KDE), which is an estimate of the density of the distribution from where the samples are drawn that employs a standard normal kernel, and a Locally Weighted Scatterplot Smoothing (LOWESS), which is a locally weighted linear least squares giving more weight to points near the point whose response is predicted. We refer to \cite{Waskom2021} for their precise definition. The LOWESS shows that} if $\theta > 0.3$, the number of levels increases\textcolor{reviewV2}{, almost linearly,} as the strong threshold increases. \textcolor{reviewV2}{Indeed, the KDE displays higher density spots in the upper right zone.}
Indeed, a larger value of $\theta$ means that more connections are kept and \textcolor{reviewV2}{the number of unknowns between two connecting levels is only partly reduced}. In particular, in 95\% of the test cases, the number of levels is a non-decreasing function with respect to $\theta$.
For each test case, we perform the least square analysis between the value of $\theta$ and the corresponding number of levels. Figure~\ref{fig:theta_vs_levels} (center and right) shows that in most of the test cases there is a significant correlation (p-value $< 10^{-5}$) between these two variables.

\textcolor{reviewV2}{W}e carried out the same set of experiment\textcolor{reviewV2}{s} varying the value of $\theta$ and recording the corresponding size of the coarse matrix $\mathrm A_H$ built by BoomerAMG. The results of the least square \textcolor{reviewV2}{analysis seems to indicate there is no correlation} between these two variables.
}

\color{black}

\subsection{Relation between $\theta$ and $\rho$}
\textcolor{reviewV2}{In this section w}e investigate the relation between the strong threshold parameter $\theta$ and the corresponding approximated convergence factor $\rho$. The results reported in Table~\ref{tab:rho1} have been obtained with a diffusion coefficient that has a ``strides" pattern (Figure~\ref{fig:mu_pattern}(c)), while Table~\ref{tab:rho2} displays analogous results on the checkerboard pattern (Figure~\ref{fig:mu_pattern}(d)).

We have computed the value of $\rho$ and the corresponding iteration counts as a function of the \textcolor{reviewV2}{value of} $\varepsilon$ of the diffusion coefficient (Eq.~\eqref{eq:mu_patt}) and the mesh size $h$. The value of $\theta$ is kept fixed for each test.

By comparing one test with the others, we can determine if the different value of the strong threshold parameter $\theta$ has affected the convergence factor $\rho$ of the linear solver. \textcolor{ultramarine}{Twenty-five values of $\theta$ in $[0.02,0.9]$ have been chosen. In Tables~\ref{tab:rho1} and~\ref{tab:rho2} we report the results for three values of $\theta$ that are representative of the obtained results when $\theta$ is ``small", ``medium" and ``large"}\textcolor{reviewV2}{, namely $\theta=0.24, 0.48, 0.72$}. 
The values of $\varepsilon$ go from $0.0$ (yielding the standard Laplacian problem with uniform diffusion), to $9.5$, which produces a quite large discontinuity in the diffusion coefficient $\mu$. 

From the results of Tables~\ref{tab:rho1} and~\ref{tab:rho2}, it is clear that, if the choice of strong threshold $\theta$ is appropriate, there is almost always uniform convergence, independently \textcolor{reviewV2}{of} the mesh size $h$. This confirms that the AMG method works as \textcolor{reviewV2}{expected} also with a diffusion coefficient $\mu$ that presents large discontinuities, provided that $\theta$ is appropriately chosen.

The results reported in Tables~\ref{tab:rho1} and~\ref{tab:rho2} also show that for large values of the strong threshold parameter ($\theta=0.72$), the approximate convergence factor $\rho$ \textcolor{reviewV2}{increases}, i.e.\ \textcolor{reviewV2}{the convergence properties of the AMG method seems to deteriorate}. A possible explanation is the following: \textcolor{ultramarine}{as we mentioned in Section~\ref{sec:theta_vs_lvls}, from Figure~\ref{fig:theta_vs_levels} we can conclude that a larger value of $\theta$ implies that a larger number of levels will be needed by the AMG algorithm. This might lead to a deterioration of the convergence rates.}

For the test cases that present less pronounced discontinuities, the value of $\theta=0.25$ (which is almost the standard literature value) provides uniform convergence. On the other hand, we notice that in the strongly heterogeneous cases (i.e.\ when $\varepsilon$ is large) deviating from the literature value of $\theta$ can result in a significant improvement is the approximate convergence factor. 

The results shown in Figure~\ref{tab:rho1} seem to indicate that choosing $\theta$ differently from the standard value suggested in literature does not result in any significant improvement. On the other hand, the plots of the four finest mesh refinements of Figure~\ref{tab:rho2} reveal that a significant boost in the performance could be obtained. An optimal choice of the strong threshold could bring up to $33\%$ \textcolor{reviewV2}{speed-up} \textcolor{ultramarine}{w.r.t.\ the default choice of $\theta=0.25$}.

\begin{table}[htbp]
\hspace{2cm}\hspace{-9.8cm}
\begin{minipage}[c]{0.7\textwidth}
\resizebox{2.5\textwidth}{!}{
\definecolor{ccc0_0}{rgb}{0.17058823529411765,0.4946558433997788,0.9667184042691874}
\definecolor{ccc0_1}{rgb}{0.303921568627451,0.30315267411304353,0.9881654720812594}
\definecolor{ccc0_2}{rgb}{0.3666666666666667,0.20791169081775931,0.9945218953682733}
\definecolor{ccc0_3}{rgb}{0.39803921568627454,0.1594757912099808,0.9967953249171991}
\definecolor{ccc0_4}{rgb}{0.3588235294117647,0.2199463578396686,0.9938591368952737}
\definecolor{ccc0_5}{rgb}{0.3509803921568627,0.23194764145389815,0.9931586661366362}
\definecolor{ccc0_6}{rgb}{0.3431372549019608,0.24391372010837714,0.9924205096719357}
\definecolor{ccc0_7}{rgb}{0.32745098039215687,0.2677330033224679,0.9908312530915603}
\definecolor{ccc1_0}{rgb}{0.18627450980392157,0.47309355683601007,0.9697969360350095}
\definecolor{ccc1_1}{rgb}{0.31176470588235294,0.2913897468893246,0.989091608371146}
\definecolor{ccc1_2}{rgb}{0.3666666666666667,0.20791169081775931,0.9945218953682733}
\definecolor{ccc1_3}{rgb}{0.38235294117647056,0.18374951781657034,0.9957341762950345}
\definecolor{ccc1_4}{rgb}{0.3588235294117647,0.2199463578396686,0.9938591368952737}
\definecolor{ccc1_5}{rgb}{0.3509803921568627,0.23194764145389815,0.9931586661366362}
\definecolor{ccc1_6}{rgb}{0.3509803921568627,0.23194764145389815,0.9931586661366362}
\definecolor{ccc1_7}{rgb}{0.32745098039215687,0.2677330033224679,0.9908312530915603}
\definecolor{ccc2_0}{rgb}{0.20980392156862748,0.44021574083098725,0.9741386021045101}
\definecolor{ccc2_1}{rgb}{0.32745098039215687,0.2677330033224679,0.9908312530915603}
\definecolor{ccc2_2}{rgb}{0.3666666666666667,0.20791169081775931,0.9945218953682733}
\definecolor{ccc2_3}{rgb}{0.37450980392156863,0.19584546700716696,0.9951469164070644}
\definecolor{ccc2_4}{rgb}{0.3666666666666667,0.20791169081775931,0.9945218953682733}
\definecolor{ccc2_5}{rgb}{0.3666666666666667,0.20791169081775931,0.9945218953682733}
\definecolor{ccc2_6}{rgb}{0.3588235294117647,0.2199463578396686,0.9938591368952737}
\definecolor{ccc2_7}{rgb}{0.3196078431372549,0.2795825925967438,0.989980213280707}
\definecolor{ccc3_0}{rgb}{0.22549019607843135,0.41796034488678346,0.9768483177596007}
\definecolor{ccc3_1}{rgb}{0.32745098039215687,0.2677330033224679,0.9908312530915603}
\definecolor{ccc3_2}{rgb}{0.3666666666666667,0.20791169081775931,0.9945218953682733}
\definecolor{ccc3_3}{rgb}{0.3666666666666667,0.20791169081775931,0.9945218953682733}
\definecolor{ccc3_4}{rgb}{0.3588235294117647,0.2199463578396686,0.9938591368952737}
\definecolor{ccc3_5}{rgb}{0.3666666666666667,0.20791169081775931,0.9945218953682733}
\definecolor{ccc3_6}{rgb}{0.3588235294117647,0.2199463578396686,0.9938591368952737}
\definecolor{ccc3_7}{rgb}{0.31176470588235294,0.2913897468893246,0.989091608371146}
\definecolor{ccc4_0}{rgb}{0.22549019607843135,0.41796034488678346,0.9768483177596007}
\definecolor{ccc4_1}{rgb}{0.33529411764705885,0.2558427775944356,0.9916446955107427}
\definecolor{ccc4_2}{rgb}{0.3588235294117647,0.2199463578396686,0.9938591368952737}
\definecolor{ccc4_3}{rgb}{0.32745098039215687,0.2677330033224679,0.9908312530915603}
\definecolor{ccc4_4}{rgb}{0.3509803921568627,0.23194764145389815,0.9931586661366362}
\definecolor{ccc4_5}{rgb}{0.3588235294117647,0.2199463578396686,0.9938591368952737}
\definecolor{ccc4_6}{rgb}{0.3509803921568627,0.23194764145389815,0.9931586661366362}
\definecolor{ccc4_7}{rgb}{0.303921568627451,0.30315267411304353,0.9881654720812594}
\definecolor{ccc5_0}{rgb}{0.22549019607843135,0.41796034488678346,0.9768483177596007}
\definecolor{ccc5_1}{rgb}{0.33529411764705885,0.2558427775944356,0.9916446955107427}
\definecolor{ccc5_2}{rgb}{0.3509803921568627,0.23194764145389815,0.9931586661366362}
\definecolor{ccc5_3}{rgb}{0.3196078431372549,0.2795825925967438,0.989980213280707}
\definecolor{ccc5_4}{rgb}{0.3509803921568627,0.23194764145389815,0.9931586661366362}
\definecolor{ccc5_5}{rgb}{0.3588235294117647,0.2199463578396686,0.9938591368952737}
\definecolor{ccc5_6}{rgb}{0.3509803921568627,0.23194764145389815,0.9931586661366362}
\definecolor{ccc5_7}{rgb}{0.303921568627451,0.30315267411304353,0.9881654720812594}
\definecolor{ccc6_0}{rgb}{0.22549019607843135,0.41796034488678346,0.9768483177596007}
\definecolor{ccc6_1}{rgb}{0.33529411764705885,0.2558427775944356,0.9916446955107427}
\definecolor{ccc6_2}{rgb}{0.3509803921568627,0.23194764145389815,0.9931586661366362}
\definecolor{ccc6_3}{rgb}{0.3196078431372549,0.2795825925967438,0.989980213280707}
\definecolor{ccc6_4}{rgb}{0.3509803921568627,0.23194764145389815,0.9931586661366362}
\definecolor{ccc6_5}{rgb}{0.3509803921568627,0.23194764145389815,0.9931586661366362}
\definecolor{ccc6_6}{rgb}{0.3509803921568627,0.23194764145389815,0.9931586661366362}
\definecolor{ccc6_7}{rgb}{0.303921568627451,0.30315267411304353,0.9881654720812594}
\definecolor{ccc7_0}{rgb}{0.22549019607843135,0.41796034488678346,0.9768483177596007}
\definecolor{ccc7_1}{rgb}{0.33529411764705885,0.2558427775944356,0.9916446955107427}
\definecolor{ccc7_2}{rgb}{0.3509803921568627,0.23194764145389815,0.9931586661366362}
\definecolor{ccc7_3}{rgb}{0.31176470588235294,0.2913897468893246,0.989091608371146}
\definecolor{ccc7_4}{rgb}{0.3509803921568627,0.23194764145389815,0.9931586661366362}
\definecolor{ccc7_5}{rgb}{0.3509803921568627,0.23194764145389815,0.9931586661366362}
\definecolor{ccc7_6}{rgb}{0.3509803921568627,0.23194764145389815,0.9931586661366362}
\definecolor{ccc7_7}{rgb}{0.303921568627451,0.30315267411304353,0.9881654720812594}
\definecolor{ccc8_0}{rgb}{0.22549019607843135,0.41796034488678346,0.9768483177596007}
\definecolor{ccc8_1}{rgb}{0.33529411764705885,0.2558427775944356,0.9916446955107427}
\definecolor{ccc8_2}{rgb}{0.3509803921568627,0.23194764145389815,0.9931586661366362}
\definecolor{ccc8_3}{rgb}{0.31176470588235294,0.2913897468893246,0.989091608371146}
\definecolor{ccc8_4}{rgb}{0.3509803921568627,0.23194764145389815,0.9931586661366362}
\definecolor{ccc8_5}{rgb}{0.3509803921568627,0.23194764145389815,0.9931586661366362}
\definecolor{ccc8_6}{rgb}{0.3509803921568627,0.23194764145389815,0.9931586661366362}
\definecolor{ccc8_7}{rgb}{0.303921568627451,0.30315267411304353,0.9881654720812594}
\definecolor{ccc9_0}{rgb}{0.22549019607843135,0.41796034488678346,0.9768483177596007}
\definecolor{ccc9_1}{rgb}{0.33529411764705885,0.2558427775944356,0.9916446955107427}
\definecolor{ccc9_2}{rgb}{0.3509803921568627,0.23194764145389815,0.9931586661366362}
\definecolor{ccc9_3}{rgb}{0.31176470588235294,0.2913897468893246,0.989091608371146}
\definecolor{ccc9_4}{rgb}{0.3509803921568627,0.23194764145389815,0.9931586661366362}
\definecolor{ccc9_5}{rgb}{0.3509803921568627,0.23194764145389815,0.9931586661366362}
\definecolor{ccc9_6}{rgb}{0.3509803921568627,0.23194764145389815,0.9931586661366362}
\definecolor{ccc9_7}{rgb}{0.303921568627451,0.30315267411304353,0.9881654720812594}
\definecolor{ccc10_0}{rgb}{0.22549019607843135,0.41796034488678346,0.9768483177596007}
\definecolor{ccc10_1}{rgb}{0.33529411764705885,0.2558427775944356,0.9916446955107427}
\definecolor{ccc10_2}{rgb}{0.3509803921568627,0.23194764145389815,0.9931586661366362}
\definecolor{ccc10_3}{rgb}{0.31176470588235294,0.2913897468893246,0.989091608371146}
\definecolor{ccc10_4}{rgb}{0.3509803921568627,0.23194764145389815,0.9931586661366362}
\definecolor{ccc10_5}{rgb}{0.3509803921568627,0.23194764145389815,0.9931586661366362}
\definecolor{ccc10_6}{rgb}{0.3509803921568627,0.23194764145389815,0.9931586661366362}
\definecolor{ccc10_7}{rgb}{0.303921568627451,0.30315267411304353,0.9881654720812594}
\definecolor{ccc11_0}{rgb}{0.22549019607843135,0.41796034488678346,0.9768483177596007}
\definecolor{ccc11_1}{rgb}{0.33529411764705885,0.2558427775944356,0.9916446955107427}
\definecolor{ccc11_2}{rgb}{0.3509803921568627,0.23194764145389815,0.9931586661366362}
\definecolor{ccc11_3}{rgb}{0.31176470588235294,0.2913897468893246,0.989091608371146}
\definecolor{ccc11_4}{rgb}{0.3509803921568627,0.23194764145389815,0.9931586661366362}
\definecolor{ccc11_5}{rgb}{0.3509803921568627,0.23194764145389815,0.9931586661366362}
\definecolor{ccc11_6}{rgb}{0.3509803921568627,0.23194764145389815,0.9931586661366362}
\definecolor{ccc11_7}{rgb}{0.303921568627451,0.30315267411304353,0.9881654720812594}
\begin{tabular}{|c||c|c|c|c|c|c|c|c|}
\hline
$\varepsilon \backslash h$ & 1.25e-01& 6.25e-02& 3.12e-02& 1.56e-02& 7.81e-03& 3.91e-03& 1.95e-03& 9.77e-04\\
\hline\hline
0.0 & \cellcolor{ccc0_0!20} 0.094(9) & \cellcolor{ccc0_1!20} 0.071(8) & \cellcolor{ccc0_2!20} 0.060(8) & \cellcolor{ccc0_3!20} 0.054(8) & \cellcolor{ccc0_4!20} 0.061(9) & \cellcolor{ccc0_5!20} 0.063(9) & \cellcolor{ccc0_6!20} 0.064(9) & \cellcolor{ccc0_7!20} 0.066(10)\\
\hline
0.4 & \cellcolor{ccc1_0!20} 0.091(9) & \cellcolor{ccc1_1!20} 0.069(8) & \cellcolor{ccc1_2!20} 0.059(8) & \cellcolor{ccc1_3!20} 0.057(8) & \cellcolor{ccc1_4!20} 0.061(9) & \cellcolor{ccc1_5!20} 0.062(9) & \cellcolor{ccc1_6!20} 0.063(9) & \cellcolor{ccc1_7!20} 0.066(10)\\
\hline
0.8 & \cellcolor{ccc2_0!20} 0.087(9) & \cellcolor{ccc2_1!20} 0.066(8) & \cellcolor{ccc2_2!20} 0.059(8) & \cellcolor{ccc2_3!20} 0.058(8) & \cellcolor{ccc2_4!20} 0.059(9) & \cellcolor{ccc2_5!20} 0.060(9) & \cellcolor{ccc2_6!20} 0.061(9) & \cellcolor{ccc2_7!20} 0.068(10)\\
\hline
1.2 & \cellcolor{ccc3_0!20} 0.085(9) & \cellcolor{ccc3_1!20} 0.066(8) & \cellcolor{ccc3_2!20} 0.060(8) & \cellcolor{ccc3_3!20} 0.059(8) & \cellcolor{ccc3_4!20} 0.061(9) & \cellcolor{ccc3_5!20} 0.060(9) & \cellcolor{ccc3_6!20} 0.061(9) & \cellcolor{ccc3_7!20} 0.069(10)\\
\hline
1.6 & \cellcolor{ccc4_0!20} 0.085(9) & \cellcolor{ccc4_1!20} 0.065(8) & \cellcolor{ccc4_2!20} 0.061(8) & \cellcolor{ccc4_3!20} 0.067(9) & \cellcolor{ccc4_4!20} 0.062(9) & \cellcolor{ccc4_5!20} 0.061(9) & \cellcolor{ccc4_6!20} 0.062(9) & \cellcolor{ccc4_7!20} 0.070(10)\\
\hline
2.0 & \cellcolor{ccc5_0!20} 0.084(9) & \cellcolor{ccc5_1!20} 0.065(8) & \cellcolor{ccc5_2!20} 0.062(8) & \cellcolor{ccc5_3!20} 0.068(9) & \cellcolor{ccc5_4!20} 0.062(9) & \cellcolor{ccc5_5!20} 0.061(9) & \cellcolor{ccc5_6!20} 0.062(9) & \cellcolor{ccc5_7!20} 0.070(10)\\
\hline
2.4 & \cellcolor{ccc6_0!20} 0.084(9) & \cellcolor{ccc6_1!20} 0.065(8) & \cellcolor{ccc6_2!20} 0.062(8) & \cellcolor{ccc6_3!20} 0.068(9) & \cellcolor{ccc6_4!20} 0.062(9) & \cellcolor{ccc6_5!20} 0.062(9) & \cellcolor{ccc6_6!20} 0.062(9) & \cellcolor{ccc6_7!20} 0.070(10)\\
\hline
2.8 & \cellcolor{ccc7_0!20} 0.084(9) & \cellcolor{ccc7_1!20} 0.065(8) & \cellcolor{ccc7_2!20} 0.062(8) & \cellcolor{ccc7_3!20} 0.069(9) & \cellcolor{ccc7_4!20} 0.062(9) & \cellcolor{ccc7_5!20} 0.062(9) & \cellcolor{ccc7_6!20} 0.063(9) & \cellcolor{ccc7_7!20} 0.070(10)\\
\hline
3.5 & \cellcolor{ccc8_0!20} 0.084(9) & \cellcolor{ccc8_1!20} 0.065(8) & \cellcolor{ccc8_2!20} 0.062(8) & \cellcolor{ccc8_3!20} 0.069(9) & \cellcolor{ccc8_4!20} 0.062(9) & \cellcolor{ccc8_5!20} 0.062(9) & \cellcolor{ccc8_6!20} 0.063(9) & \cellcolor{ccc8_7!20} 0.070(10)\\
\hline
5.0 & \cellcolor{ccc9_0!20} 0.084(9) & \cellcolor{ccc9_1!20} 0.065(8) & \cellcolor{ccc9_2!20} 0.062(8) & \cellcolor{ccc9_3!20} 0.069(9) & \cellcolor{ccc9_4!20} 0.062(9) & \cellcolor{ccc9_5!20} 0.062(9) & \cellcolor{ccc9_6!20} 0.063(9) & \cellcolor{ccc9_7!20} 0.070(10)\\
\hline
7.0 & \cellcolor{ccc10_0!20} 0.084(9) & \cellcolor{ccc10_1!20} 0.065(8) & \cellcolor{ccc10_2!20} 0.062(8) & \cellcolor{ccc10_3!20} 0.069(9) & \cellcolor{ccc10_4!20} 0.062(9) & \cellcolor{ccc10_5!20} 0.062(9) & \cellcolor{ccc10_6!20} 0.063(9) & \cellcolor{ccc10_7!20} 0.070(10)\\
\hline
9.5 & \cellcolor{ccc11_0!20} 0.084(9) & \cellcolor{ccc11_1!20} 0.065(8) & \cellcolor{ccc11_2!20} 0.062(8) & \cellcolor{ccc11_3!20} 0.069(9) & \cellcolor{ccc11_4!20} 0.062(9) & \cellcolor{ccc11_5!20} 0.062(9) & \cellcolor{ccc11_6!20} 0.063(9) & \cellcolor{ccc11_7!20} 0.070(10)\\
\hline
\end{tabular}
}
\end{minipage}
\noindent
\hspace{-2.3cm}
\begin{minipage}[c]{0.3\textwidth}\hspace{0.1cm}$\theta=0.24$\\~\\
\resizebox{0.4\textwidth}{!}{
%
%
\begin{tikzpicture}

\begin{axis}[%
width=3.566in,
height=3.566in,
at={(1.236in,0.481in)},
scale only axis,
xmin=-1,
xmax=1,
xtick={\empty},
ymin=-1,
ymax=1,
ytick={\empty},
axis background/.style={fill=white},
axis x line*=bottom,
axis y line*=left,
legend style={legend cell align=left, align=left, draw=white!15!black}
]

\addplot[area legend, draw=black, fill=white]
table[row sep=crcr] {%
x	y\\
-1	-1\\
-1	1\\
-0.5	1\\
-0.5	-1\\
}--cycle;

\addplot[area legend, draw=black, fill=black!10]
table[row sep=crcr] {%
x	y\\
-0.5	-1\\
-0.5	1\\
0	1\\
0	-1\\
}--cycle;

\addplot[area legend, draw=black, fill=white]
table[row sep=crcr] {%
x	y\\
0	-1\\
0	1\\
0.5	1\\
0.5	-1\\
}--cycle;

\addplot[area legend, draw=black, fill=black!10]
table[row sep=crcr] {%
x	y\\
0.5	-1\\
0.5	1\\
1	1\\
1	-1\\
}--cycle;

\end{axis}
\end{tikzpicture}%
}
\end{minipage}\\

\vspace{1.5cm}\hspace{2cm}\hspace{-9.8cm}
\begin{minipage}[c]{0.7\textwidth}
\resizebox{2.5\textwidth}{!}{
\definecolor{ccc0_0}{rgb}{0.17058823529411765,0.4946558433997788,0.9667184042691874}
\definecolor{ccc0_1}{rgb}{0.303921568627451,0.30315267411304353,0.9881654720812594}
\definecolor{ccc0_2}{rgb}{0.3666666666666667,0.20791169081775931,0.9945218953682733}
\definecolor{ccc0_3}{rgb}{0.39803921568627454,0.1594757912099808,0.9967953249171991}
\definecolor{ccc0_4}{rgb}{0.3588235294117647,0.2199463578396686,0.9938591368952737}
\definecolor{ccc0_5}{rgb}{0.3509803921568627,0.23194764145389815,0.9931586661366362}
\definecolor{ccc0_6}{rgb}{0.3431372549019608,0.24391372010837714,0.9924205096719357}
\definecolor{ccc0_7}{rgb}{0.32745098039215687,0.2677330033224679,0.9908312530915603}
\definecolor{ccc1_0}{rgb}{0.18627450980392157,0.47309355683601007,0.9697969360350095}
\definecolor{ccc1_1}{rgb}{0.3196078431372549,0.2795825925967438,0.989980213280707}
\definecolor{ccc1_2}{rgb}{0.38235294117647056,0.18374951781657034,0.9957341762950345}
\definecolor{ccc1_3}{rgb}{0.38235294117647056,0.18374951781657034,0.9957341762950345}
\definecolor{ccc1_4}{rgb}{0.3666666666666667,0.20791169081775931,0.9945218953682733}
\definecolor{ccc1_5}{rgb}{0.3666666666666667,0.20791169081775931,0.9945218953682733}
\definecolor{ccc1_6}{rgb}{0.3588235294117647,0.2199463578396686,0.9938591368952737}
\definecolor{ccc1_7}{rgb}{0.32745098039215687,0.2677330033224679,0.9908312530915603}
\definecolor{ccc2_0}{rgb}{0.20980392156862748,0.44021574083098725,0.9741386021045101}
\definecolor{ccc2_1}{rgb}{0.32745098039215687,0.2677330033224679,0.9908312530915603}
\definecolor{ccc2_2}{rgb}{0.3666666666666667,0.20791169081775931,0.9945218953682733}
\definecolor{ccc2_3}{rgb}{0.37450980392156863,0.19584546700716696,0.9951469164070644}
\definecolor{ccc2_4}{rgb}{0.3666666666666667,0.20791169081775931,0.9945218953682733}
\definecolor{ccc2_5}{rgb}{0.3666666666666667,0.20791169081775931,0.9945218953682733}
\definecolor{ccc2_6}{rgb}{0.3588235294117647,0.2199463578396686,0.9938591368952737}
\definecolor{ccc2_7}{rgb}{0.3196078431372549,0.2795825925967438,0.989980213280707}
\definecolor{ccc3_0}{rgb}{0.2647058823529412,0.3612416661871529,0.9829730996839018}
\definecolor{ccc3_1}{rgb}{0.3196078431372549,0.2795825925967438,0.989980213280707}
\definecolor{ccc3_2}{rgb}{0.3666666666666667,0.20791169081775931,0.9945218953682733}
\definecolor{ccc3_3}{rgb}{0.28823529411764703,0.3265387128400833,0.9862007473534026}
\definecolor{ccc3_4}{rgb}{0.3509803921568627,0.23194764145389815,0.9931586661366362}
\definecolor{ccc3_5}{rgb}{0.296078431372549,0.31486958889350786,0.987201839553569}
\definecolor{ccc3_6}{rgb}{0.2019607843137255,0.4512440570453228,0.9727282722446048}
\definecolor{ccc3_7}{rgb}{0.18627450980392157,0.47309355683601007,0.9697969360350095}
\definecolor{ccc4_0}{rgb}{0.27254901960784317,0.34972651120626114,0.9840863373026044}
\definecolor{ccc4_1}{rgb}{0.3196078431372549,0.2795825925967438,0.989980213280707}
\definecolor{ccc4_2}{rgb}{0.3666666666666667,0.20791169081775931,0.9945218953682733}
\definecolor{ccc4_3}{rgb}{0.2803921568627451,0.33815827481581706,0.9851622334675065}
\definecolor{ccc4_4}{rgb}{0.3509803921568627,0.23194764145389815,0.9931586661366362}
\definecolor{ccc4_5}{rgb}{0.3196078431372549,0.2795825925967438,0.989980213280707}
\definecolor{ccc4_6}{rgb}{0.23333333333333334,0.40673664307580015,0.9781476007338057}
\definecolor{ccc4_7}{rgb}{0.17843137254901964,0.48391142410030147,0.9682760409157589}
\definecolor{ccc5_0}{rgb}{0.27254901960784317,0.34972651120626114,0.9840863373026044}
\definecolor{ccc5_1}{rgb}{0.3196078431372549,0.2795825925967438,0.989980213280707}
\definecolor{ccc5_2}{rgb}{0.3666666666666667,0.20791169081775931,0.9945218953682733}
\definecolor{ccc5_3}{rgb}{0.2803921568627451,0.33815827481581706,0.9851622334675065}
\definecolor{ccc5_4}{rgb}{0.3509803921568627,0.23194764145389815,0.9931586661366362}
\definecolor{ccc5_5}{rgb}{0.32745098039215687,0.2677330033224679,0.9908312530915603}
\definecolor{ccc5_6}{rgb}{0.22549019607843135,0.41796034488678346,0.9768483177596007}
\definecolor{ccc5_7}{rgb}{0.20980392156862748,0.44021574083098725,0.9741386021045101}
\definecolor{ccc6_0}{rgb}{0.2803921568627451,0.33815827481581706,0.9851622334675065}
\definecolor{ccc6_1}{rgb}{0.3196078431372549,0.2795825925967438,0.989980213280707}
\definecolor{ccc6_2}{rgb}{0.3666666666666667,0.20791169081775931,0.9945218953682733}
\definecolor{ccc6_3}{rgb}{0.27254901960784317,0.34972651120626114,0.9840863373026044}
\definecolor{ccc6_4}{rgb}{0.3509803921568627,0.23194764145389815,0.9931586661366362}
\definecolor{ccc6_5}{rgb}{0.32745098039215687,0.2677330033224679,0.9908312530915603}
\definecolor{ccc6_6}{rgb}{0.22549019607843135,0.41796034488678346,0.9768483177596007}
\definecolor{ccc6_7}{rgb}{0.21764705882352942,0.42912060877260894,0.9755119679804366}
\definecolor{ccc7_0}{rgb}{0.2803921568627451,0.33815827481581706,0.9851622334675065}
\definecolor{ccc7_1}{rgb}{0.3196078431372549,0.2795825925967438,0.989980213280707}
\definecolor{ccc7_2}{rgb}{0.3666666666666667,0.20791169081775931,0.9945218953682733}
\definecolor{ccc7_3}{rgb}{0.27254901960784317,0.34972651120626114,0.9840863373026044}
\definecolor{ccc7_4}{rgb}{0.3509803921568627,0.23194764145389815,0.9931586661366362}
\definecolor{ccc7_5}{rgb}{0.32745098039215687,0.2677330033224679,0.9908312530915603}
\definecolor{ccc7_6}{rgb}{0.22549019607843135,0.41796034488678346,0.9768483177596007}
\definecolor{ccc7_7}{rgb}{0.22549019607843135,0.41796034488678346,0.9768483177596007}
\definecolor{ccc8_0}{rgb}{0.2803921568627451,0.33815827481581706,0.9851622334675065}
\definecolor{ccc8_1}{rgb}{0.31176470588235294,0.2913897468893246,0.989091608371146}
\definecolor{ccc8_2}{rgb}{0.3666666666666667,0.20791169081775931,0.9945218953682733}
\definecolor{ccc8_3}{rgb}{0.27254901960784317,0.34972651120626114,0.9840863373026044}
\definecolor{ccc8_4}{rgb}{0.3509803921568627,0.23194764145389815,0.9931586661366362}
\definecolor{ccc8_5}{rgb}{0.32745098039215687,0.2677330033224679,0.9908312530915603}
\definecolor{ccc8_6}{rgb}{0.23333333333333334,0.40673664307580015,0.9781476007338057}
\definecolor{ccc8_7}{rgb}{0.2568627450980392,0.37270199199091397,0.9818225628535369}
\definecolor{ccc9_0}{rgb}{0.2803921568627451,0.33815827481581706,0.9851622334675065}
\definecolor{ccc9_1}{rgb}{0.31176470588235294,0.2913897468893246,0.989091608371146}
\definecolor{ccc9_2}{rgb}{0.3666666666666667,0.20791169081775931,0.9945218953682733}
\definecolor{ccc9_3}{rgb}{0.27254901960784317,0.34972651120626114,0.9840863373026044}
\definecolor{ccc9_4}{rgb}{0.3509803921568627,0.23194764145389815,0.9931586661366362}
\definecolor{ccc9_5}{rgb}{0.32745098039215687,0.2677330033224679,0.9908312530915603}
\definecolor{ccc9_6}{rgb}{0.23333333333333334,0.40673664307580015,0.9781476007338057}
\definecolor{ccc9_7}{rgb}{0.2568627450980392,0.37270199199091397,0.9818225628535369}
\definecolor{ccc10_0}{rgb}{0.2803921568627451,0.33815827481581706,0.9851622334675065}
\definecolor{ccc10_1}{rgb}{0.31176470588235294,0.2913897468893246,0.989091608371146}
\definecolor{ccc10_2}{rgb}{0.3666666666666667,0.20791169081775931,0.9945218953682733}
\definecolor{ccc10_3}{rgb}{0.27254901960784317,0.34972651120626114,0.9840863373026044}
\definecolor{ccc10_4}{rgb}{0.3509803921568627,0.23194764145389815,0.9931586661366362}
\definecolor{ccc10_5}{rgb}{0.32745098039215687,0.2677330033224679,0.9908312530915603}
\definecolor{ccc10_6}{rgb}{0.23333333333333334,0.40673664307580015,0.9781476007338057}
\definecolor{ccc10_7}{rgb}{0.2568627450980392,0.37270199199091397,0.9818225628535369}
\definecolor{ccc11_0}{rgb}{0.2803921568627451,0.33815827481581706,0.9851622334675065}
\definecolor{ccc11_1}{rgb}{0.31176470588235294,0.2913897468893246,0.989091608371146}
\definecolor{ccc11_2}{rgb}{0.3666666666666667,0.20791169081775931,0.9945218953682733}
\definecolor{ccc11_3}{rgb}{0.27254901960784317,0.34972651120626114,0.9840863373026044}
\definecolor{ccc11_4}{rgb}{0.3509803921568627,0.23194764145389815,0.9931586661366362}
\definecolor{ccc11_5}{rgb}{0.32745098039215687,0.2677330033224679,0.9908312530915603}
\definecolor{ccc11_6}{rgb}{0.23333333333333334,0.40673664307580015,0.9781476007338057}
\definecolor{ccc11_7}{rgb}{0.2568627450980392,0.37270199199091397,0.9818225628535369}
\begin{tabular}{|c||c|c|c|c|c|c|c|c|}
\hline
$\varepsilon \backslash h$ & 1.25e-01& 6.25e-02& 3.12e-02& 1.56e-02& 7.81e-03& 3.91e-03& 1.95e-03& 9.77e-04\\
\hline\hline
0.0 & \cellcolor{ccc0_0!20} 0.094(9) & \cellcolor{ccc0_1!20} 0.071(8) & \cellcolor{ccc0_2!20} 0.060(8) & \cellcolor{ccc0_3!20} 0.054(8) & \cellcolor{ccc0_4!20} 0.061(9) & \cellcolor{ccc0_5!20} 0.063(9) & \cellcolor{ccc0_6!20} 0.064(9) & \cellcolor{ccc0_7!20} 0.066(10)\\
\hline
0.4 & \cellcolor{ccc1_0!20} 0.091(9) & \cellcolor{ccc1_1!20} 0.068(8) & \cellcolor{ccc1_2!20} 0.057(8) & \cellcolor{ccc1_3!20} 0.057(8) & \cellcolor{ccc1_4!20} 0.059(9) & \cellcolor{ccc1_5!20} 0.060(9) & \cellcolor{ccc1_6!20} 0.061(9) & \cellcolor{ccc1_7!20} 0.066(10)\\
\hline
0.8 & \cellcolor{ccc2_0!20} 0.087(9) & \cellcolor{ccc2_1!20} 0.066(8) & \cellcolor{ccc2_2!20} 0.059(8) & \cellcolor{ccc2_3!20} 0.058(8) & \cellcolor{ccc2_4!20} 0.059(9) & \cellcolor{ccc2_5!20} 0.060(9) & \cellcolor{ccc2_6!20} 0.061(9) & \cellcolor{ccc2_7!20} 0.068(10)\\
\hline
1.2 & \cellcolor{ccc3_0!20} 0.077(8) & \cellcolor{ccc3_1!20} 0.068(8) & \cellcolor{ccc3_2!20} 0.059(8) & \cellcolor{ccc3_3!20} 0.073(9) & \cellcolor{ccc3_4!20} 0.063(9) & \cellcolor{ccc3_5!20} 0.072(10) & \cellcolor{ccc3_6!20} 0.089(11) & \cellcolor{ccc3_7!20} 0.091(11)\\
\hline
1.6 & \cellcolor{ccc4_0!20} 0.076(8) & \cellcolor{ccc4_1!20} 0.068(8) & \cellcolor{ccc4_2!20} 0.059(8) & \cellcolor{ccc4_3!20} 0.075(9) & \cellcolor{ccc4_4!20} 0.063(9) & \cellcolor{ccc4_5!20} 0.068(9) & \cellcolor{ccc4_6!20} 0.083(11) & \cellcolor{ccc4_7!20} 0.092(11)\\
\hline
2.0 & \cellcolor{ccc5_0!20} 0.076(8) & \cellcolor{ccc5_1!20} 0.068(8) & \cellcolor{ccc5_2!20} 0.059(8) & \cellcolor{ccc5_3!20} 0.075(9) & \cellcolor{ccc5_4!20} 0.063(9) & \cellcolor{ccc5_5!20} 0.067(9) & \cellcolor{ccc5_6!20} 0.085(10) & \cellcolor{ccc5_7!20} 0.087(11)\\
\hline
2.4 & \cellcolor{ccc6_0!20} 0.075(8) & \cellcolor{ccc6_1!20} 0.068(8) & \cellcolor{ccc6_2!20} 0.059(8) & \cellcolor{ccc6_3!20} 0.076(9) & \cellcolor{ccc6_4!20} 0.063(9) & \cellcolor{ccc6_5!20} 0.067(9) & \cellcolor{ccc6_6!20} 0.084(10) & \cellcolor{ccc6_7!20} 0.086(11)\\
\hline
2.8 & \cellcolor{ccc7_0!20} 0.075(8) & \cellcolor{ccc7_1!20} 0.068(8) & \cellcolor{ccc7_2!20} 0.059(8) & \cellcolor{ccc7_3!20} 0.076(9) & \cellcolor{ccc7_4!20} 0.063(9) & \cellcolor{ccc7_5!20} 0.067(9) & \cellcolor{ccc7_6!20} 0.084(10) & \cellcolor{ccc7_7!20} 0.085(11)\\
\hline
3.5 & \cellcolor{ccc8_0!20} 0.075(8) & \cellcolor{ccc8_1!20} 0.069(8) & \cellcolor{ccc8_2!20} 0.059(8) & \cellcolor{ccc8_3!20} 0.076(9) & \cellcolor{ccc8_4!20} 0.062(9) & \cellcolor{ccc8_5!20} 0.067(9) & \cellcolor{ccc8_6!20} 0.083(10) & \cellcolor{ccc8_7!20} 0.079(10)\\
\hline
5.0 & \cellcolor{ccc9_0!20} 0.075(8) & \cellcolor{ccc9_1!20} 0.069(8) & \cellcolor{ccc9_2!20} 0.059(8) & \cellcolor{ccc9_3!20} 0.076(9) & \cellcolor{ccc9_4!20} 0.062(9) & \cellcolor{ccc9_5!20} 0.067(9) & \cellcolor{ccc9_6!20} 0.083(10) & \cellcolor{ccc9_7!20} 0.079(10)\\
\hline
7.0 & \cellcolor{ccc10_0!20} 0.075(8) & \cellcolor{ccc10_1!20} 0.069(8) & \cellcolor{ccc10_2!20} 0.059(8) & \cellcolor{ccc10_3!20} 0.076(9) & \cellcolor{ccc10_4!20} 0.062(9) & \cellcolor{ccc10_5!20} 0.067(9) & \cellcolor{ccc10_6!20} 0.083(10) & \cellcolor{ccc10_7!20} 0.079(10)\\
\hline
9.5 & \cellcolor{ccc11_0!20} 0.075(8) & \cellcolor{ccc11_1!20} 0.069(8) & \cellcolor{ccc11_2!20} 0.059(8) & \cellcolor{ccc11_3!20} 0.076(9) & \cellcolor{ccc11_4!20} 0.062(9) & \cellcolor{ccc11_5!20} 0.067(9) & \cellcolor{ccc11_6!20} 0.083(10) & \cellcolor{ccc11_7!20} 0.079(10)\\
\hline
\end{tabular}
}
\end{minipage}
\noindent
\hspace{-2.3cm}
\begin{minipage}[c]{0.3\textwidth}\hspace{0.1cm}$\theta=0.48$\\~\\
\resizebox{0.4\textwidth}{!}{
%
%
\begin{tikzpicture}

\begin{axis}[%
width=3.566in,
height=3.566in,
at={(1.236in,0.481in)},
scale only axis,
xmin=-1,
xmax=1,
xtick={\empty},
ymin=-1,
ymax=1,
ytick={\empty},
axis background/.style={fill=white},
axis x line*=bottom,
axis y line*=left,
legend style={legend cell align=left, align=left, draw=white!15!black}
]

\addplot[area legend, draw=black, fill=white]
table[row sep=crcr] {%
x	y\\
-1	-1\\
-1	1\\
-0.5	1\\
-0.5	-1\\
}--cycle;

\addplot[area legend, draw=black, fill=black!10]
table[row sep=crcr] {%
x	y\\
-0.5	-1\\
-0.5	1\\
0	1\\
0	-1\\
}--cycle;

\addplot[area legend, draw=black, fill=white]
table[row sep=crcr] {%
x	y\\
0	-1\\
0	1\\
0.5	1\\
0.5	-1\\
}--cycle;

\addplot[area legend, draw=black, fill=black!10]
table[row sep=crcr] {%
x	y\\
0.5	-1\\
0.5	1\\
1	1\\
1	-1\\
}--cycle;

\end{axis}
\end{tikzpicture}%
}
\end{minipage}\\

\vspace{1.5cm}\hspace{2cm}\hspace{-9.8cm}
\begin{minipage}[c]{0.7\textwidth}
\resizebox{2.5\textwidth}{!}{
\definecolor{ccc0_0}{rgb}{0.17058823529411765,0.4946558433997788,0.9667184042691874}
\definecolor{ccc0_1}{rgb}{0.303921568627451,0.30315267411304353,0.9881654720812594}
\definecolor{ccc0_2}{rgb}{0.3666666666666667,0.20791169081775931,0.9945218953682733}
\definecolor{ccc0_3}{rgb}{0.39803921568627454,0.1594757912099808,0.9967953249171991}
\definecolor{ccc0_4}{rgb}{0.3588235294117647,0.2199463578396686,0.9938591368952737}
\definecolor{ccc0_5}{rgb}{0.3509803921568627,0.23194764145389815,0.9931586661366362}
\definecolor{ccc0_6}{rgb}{0.3431372549019608,0.24391372010837714,0.9924205096719357}
\definecolor{ccc0_7}{rgb}{0.32745098039215687,0.2677330033224679,0.9908312530915603}
\definecolor{ccc1_0}{rgb}{0.39803921568627454,0.1594757912099808,0.9967953249171991}
\definecolor{ccc1_1}{rgb}{0.17450980392156867,0.8721195109836108,0.8629289996673897}
\definecolor{ccc1_2}{rgb}{0.06862745098039214,0.6269238058941066,0.9431544344712774}
\definecolor{ccc1_3}{rgb}{0.23725490196078436,0.9160336012803335,0.8369891083319778}
\definecolor{ccc1_4}{rgb}{0.15098039215686276,0.8534437988883159,0.8721195109836108}
\definecolor{ccc1_5}{rgb}{0.40980392156862744,0.989980213280707,0.7553827347189938}
\definecolor{ccc1_6}{rgb}{0.37058823529411766,0.9794097676013659,0.7752039761111298}
\definecolor{ccc1_7}{rgb}{0.38627450980392153,0.9840863373026044,0.767362681448697}
\definecolor{ccc2_0}{rgb}{0.4372549019607843,0.0984002782796427,0.99878599219429}
\definecolor{ccc2_1}{rgb}{0.15098039215686276,0.8534437988883159,0.8721195109836108}
\definecolor{ccc2_2}{rgb}{0.09999999999999998,0.5877852522924731,0.9510565162951535}
\definecolor{ccc2_3}{rgb}{0.19019607843137254,0.883909710213612,0.8566380778638628}
\definecolor{ccc2_4}{rgb}{0.16666666666666663,0.8660254037844386,0.8660254037844387}
\definecolor{ccc2_5}{rgb}{0.41764705882352937,0.9916446955107427,0.7513318895568734}
\definecolor{ccc2_6}{rgb}{0.27647058823529413,0.9389883606150565,0.8197404829072211}
\definecolor{ccc2_7}{rgb}{0.32352941176470584,0.961825643172819,0.7980172272802395}
\definecolor{ccc3_0}{rgb}{0.4137254901960784,0.13510524748139296,0.9977051801738729}
\definecolor{ccc3_1}{rgb}{0.1941176470588235,0.46220388354031317,0.9712810319161138}
\definecolor{ccc3_2}{rgb}{0.09215686274509804,0.5977074592660936,0.9491349440359013}
\definecolor{ccc3_3}{rgb}{0.00588235294117645,0.700543037593291,0.9256376597815563}
\definecolor{ccc3_4}{rgb}{0.0725490196078431,0.7829276104921027,0.9005867023006374}
\definecolor{ccc3_5}{rgb}{0.11176470588235299,0.8197404829072211,0.8867736859200619}
\definecolor{ccc3_6}{rgb}{0.2529411764705882,0.9256376597815562,0.8301840308155507}
\definecolor{ccc3_7}{rgb}{0.4019607843137255,0.9881654720812594,0.7594049166547071}
\definecolor{ccc4_0}{rgb}{0.39803921568627454,0.1594757912099808,0.9967953249171991}
\definecolor{ccc4_1}{rgb}{0.24117647058823533,0.3954512068705425,0.9794097676013659}
\definecolor{ccc4_2}{rgb}{0.07647058823529412,0.6172782212897929,0.9451838281608196}
\definecolor{ccc4_3}{rgb}{0.0725490196078431,0.7829276104921027,0.9005867023006374}
\definecolor{ccc4_4}{rgb}{0.0490196078431373,0.7594049166547072,0.9084652718195236}
\definecolor{ccc4_5}{rgb}{0.10392156862745094,0.8126223709664456,0.8896040127307095}
\definecolor{ccc4_6}{rgb}{0.2058823529411764,0.8951632913550623,0.8502171357296142}
\definecolor{ccc4_7}{rgb}{0.3549019607843137,0.9741386021045101,0.7829276104921028}
\definecolor{ccc5_0}{rgb}{0.39803921568627454,0.1594757912099808,0.9967953249171991}
\definecolor{ccc5_1}{rgb}{0.24901960784313726,0.38410574917192586,0.9806347704689777}
\definecolor{ccc5_2}{rgb}{0.08431372549019611,0.6075389463388168,0.9471773565640402}
\definecolor{ccc5_3}{rgb}{0.009803921568627416,0.7179119230644418,0.9209055179449537}
\definecolor{ccc5_4}{rgb}{0.0490196078431373,0.7594049166547072,0.9084652718195236}
\definecolor{ccc5_5}{rgb}{0.19019607843137254,0.883909710213612,0.8566380778638628}
\definecolor{ccc5_6}{rgb}{0.2529411764705882,0.9256376597815562,0.8301840308155507}
\definecolor{ccc5_7}{rgb}{0.5352941176470587,0.9984636039674339,0.6872366859692628}
\definecolor{ccc6_0}{rgb}{0.39803921568627454,0.1594757912099808,0.9967953249171991}
\definecolor{ccc6_1}{rgb}{0.2568627450980392,0.37270199199091397,0.9818225628535369}
\definecolor{ccc6_2}{rgb}{0.07647058823529412,0.6172782212897929,0.9451838281608196}
\definecolor{ccc6_3}{rgb}{0.025490196078431393,0.7348449670469757,0.9160336012803335}
\definecolor{ccc6_4}{rgb}{0.0725490196078431,0.7829276104921027,0.9005867023006374}
\definecolor{ccc6_5}{rgb}{0.2058823529411764,0.8951632913550623,0.8502171357296142}
\definecolor{ccc6_6}{rgb}{0.26078431372549016,0.9302293085467402,0.8267341748257635}
\definecolor{ccc6_7}{rgb}{0.5666666666666667,0.9945218953682734,0.6691306063588582}
\definecolor{ccc7_0}{rgb}{0.39803921568627454,0.1594757912099808,0.9967953249171991}
\definecolor{ccc7_1}{rgb}{0.2568627450980392,0.37270199199091397,0.9818225628535369}
\definecolor{ccc7_2}{rgb}{0.07647058823529412,0.6172782212897929,0.9451838281608196}
\definecolor{ccc7_3}{rgb}{0.06470588235294117,0.7752039761111298,0.9032471993461288}
\definecolor{ccc7_4}{rgb}{0.08823529411764708,0.7980172272802395,0.8951632913550623}
\definecolor{ccc7_5}{rgb}{0.23725490196078436,0.9160336012803335,0.8369891083319778}
\definecolor{ccc7_6}{rgb}{0.2686274509803921,0.9346797673211106,0.8232529481575873}
\definecolor{ccc7_7}{rgb}{0.6058823529411765,0.9862007473534026,0.6459280624867872}
\definecolor{ccc8_0}{rgb}{0.39803921568627454,0.1594757912099808,0.9967953249171991}
\definecolor{ccc8_1}{rgb}{0.2568627450980392,0.37270199199091397,0.9818225628535369}
\definecolor{ccc8_2}{rgb}{0.06862745098039214,0.6269238058941066,0.9431544344712774}
\definecolor{ccc8_3}{rgb}{0.08823529411764708,0.7980172272802395,0.8951632913550623}
\definecolor{ccc8_4}{rgb}{0.11176470588235299,0.8197404829072211,0.8867736859200619}
\definecolor{ccc8_5}{rgb}{0.27647058823529413,0.9389883606150565,0.8197404829072211}
\definecolor{ccc8_6}{rgb}{0.2686274509803921,0.9346797673211106,0.8232529481575873}
\definecolor{ccc8_7}{rgb}{0.6215686274509804,0.9818225628535369,0.6364742361471414}
\definecolor{ccc9_0}{rgb}{0.39803921568627454,0.1594757912099808,0.9967953249171991}
\definecolor{ccc9_1}{rgb}{0.2568627450980392,0.37270199199091397,0.9818225628535369}
\definecolor{ccc9_2}{rgb}{0.06862745098039214,0.6269238058941066,0.9431544344712774}
\definecolor{ccc9_3}{rgb}{0.10392156862745094,0.8126223709664456,0.8896040127307095}
\definecolor{ccc9_4}{rgb}{0.11176470588235299,0.8197404829072211,0.8867736859200619}
\definecolor{ccc9_5}{rgb}{0.28431372549019607,0.9431544344712774,0.8161969123562217}
\definecolor{ccc9_6}{rgb}{0.27647058823529413,0.9389883606150565,0.8197404829072211}
\definecolor{ccc9_7}{rgb}{0.41764705882352937,0.9916446955107427,0.7513318895568734}
\definecolor{ccc10_0}{rgb}{0.39803921568627454,0.1594757912099808,0.9967953249171991}
\definecolor{ccc10_1}{rgb}{0.2568627450980392,0.37270199199091397,0.9818225628535369}
\definecolor{ccc10_2}{rgb}{0.06862745098039214,0.6269238058941066,0.9431544344712774}
\definecolor{ccc10_3}{rgb}{0.10392156862745094,0.8126223709664456,0.8896040127307095}
\definecolor{ccc10_4}{rgb}{0.11176470588235299,0.8197404829072211,0.8867736859200619}
\definecolor{ccc10_5}{rgb}{0.292156862745098,0.9471773565640402,0.8126223709664456}
\definecolor{ccc10_6}{rgb}{0.27647058823529413,0.9389883606150565,0.8197404829072211}
\definecolor{ccc10_7}{rgb}{0.41764705882352937,0.9916446955107427,0.7513318895568734}
\definecolor{ccc11_0}{rgb}{0.39803921568627454,0.1594757912099808,0.9967953249171991}
\definecolor{ccc11_1}{rgb}{0.2568627450980392,0.37270199199091397,0.9818225628535369}
\definecolor{ccc11_2}{rgb}{0.06862745098039214,0.6269238058941066,0.9431544344712774}
\definecolor{ccc11_3}{rgb}{0.10392156862745094,0.8126223709664456,0.8896040127307095}
\definecolor{ccc11_4}{rgb}{0.11176470588235299,0.8197404829072211,0.8867736859200619}
\definecolor{ccc11_5}{rgb}{0.292156862745098,0.9471773565640402,0.8126223709664456}
\definecolor{ccc11_6}{rgb}{0.27647058823529413,0.9389883606150565,0.8197404829072211}
\definecolor{ccc11_7}{rgb}{0.41764705882352937,0.9916446955107427,0.7513318895568734}
\begin{tabular}{|c||c|c|c|c|c|c|c|c|}
\hline
$\varepsilon \backslash h$ & 1.25e-01& 6.25e-02& 3.12e-02& 1.56e-02& 7.81e-03& 3.91e-03& 1.95e-03& 9.77e-04\\
\hline\hline
0.0 & \cellcolor{ccc0_0!20} 0.094(9) & \cellcolor{ccc0_1!20} 0.071(8) & \cellcolor{ccc0_2!20} 0.060(8) & \cellcolor{ccc0_3!20} 0.054(8) & \cellcolor{ccc0_4!20} 0.061(9) & \cellcolor{ccc0_5!20} 0.063(9) & \cellcolor{ccc0_6!20} 0.064(9) & \cellcolor{ccc0_7!20} 0.066(10)\\
\hline
0.4 & \cellcolor{ccc1_0!20} 0.054(7) & \cellcolor{ccc1_1!20} 0.155(12) & \cellcolor{ccc1_2!20} 0.112(10) & \cellcolor{ccc1_3!20} 0.165(13) & \cellcolor{ccc1_4!20} 0.150(13) & \cellcolor{ccc1_5!20} 0.196(15) & \cellcolor{ccc1_6!20} 0.189(15) & \cellcolor{ccc1_7!20} 0.192(16)\\
\hline
0.8 & \cellcolor{ccc2_0!20} 0.047(7) & \cellcolor{ccc2_1!20} 0.151(12) & \cellcolor{ccc2_2!20} 0.107(10) & \cellcolor{ccc2_3!20} 0.157(13) & \cellcolor{ccc2_4!20} 0.153(13) & \cellcolor{ccc2_5!20} 0.197(15) & \cellcolor{ccc2_6!20} 0.172(14) & \cellcolor{ccc2_7!20} 0.181(15)\\
\hline
1.2 & \cellcolor{ccc3_0!20} 0.052(7) & \cellcolor{ccc3_1!20} 0.090(9) & \cellcolor{ccc3_2!20} 0.108(10) & \cellcolor{ccc3_3!20} 0.123(11) & \cellcolor{ccc3_4!20} 0.136(12) & \cellcolor{ccc3_5!20} 0.143(13) & \cellcolor{ccc3_6!20} 0.169(14) & \cellcolor{ccc3_7!20} 0.195(16)\\
\hline
1.6 & \cellcolor{ccc4_0!20} 0.054(7) & \cellcolor{ccc4_1!20} 0.081(9) & \cellcolor{ccc4_2!20} 0.110(10) & \cellcolor{ccc4_3!20} 0.136(12) & \cellcolor{ccc4_4!20} 0.132(12) & \cellcolor{ccc4_5!20} 0.142(13) & \cellcolor{ccc4_6!20} 0.160(14) & \cellcolor{ccc4_7!20} 0.186(15)\\
\hline
2.0 & \cellcolor{ccc5_0!20} 0.054(7) & \cellcolor{ccc5_1!20} 0.080(9) & \cellcolor{ccc5_2!20} 0.109(10) & \cellcolor{ccc5_3!20} 0.126(11) & \cellcolor{ccc5_4!20} 0.133(12) & \cellcolor{ccc5_5!20} 0.158(13) & \cellcolor{ccc5_6!20} 0.168(14) & \cellcolor{ccc5_7!20} 0.218(17)\\
\hline
2.4 & \cellcolor{ccc6_0!20} 0.054(7) & \cellcolor{ccc6_1!20} 0.079(9) & \cellcolor{ccc6_2!20} 0.110(10) & \cellcolor{ccc6_3!20} 0.128(12) & \cellcolor{ccc6_4!20} 0.137(12) & \cellcolor{ccc6_5!20} 0.160(14) & \cellcolor{ccc6_6!20} 0.170(14) & \cellcolor{ccc6_7!20} 0.223(17)\\
\hline
2.8 & \cellcolor{ccc7_0!20} 0.054(7) & \cellcolor{ccc7_1!20} 0.079(9) & \cellcolor{ccc7_2!20} 0.111(10) & \cellcolor{ccc7_3!20} 0.135(12) & \cellcolor{ccc7_4!20} 0.139(12) & \cellcolor{ccc7_5!20} 0.166(14) & \cellcolor{ccc7_6!20} 0.171(14) & \cellcolor{ccc7_7!20} 0.230(18)\\
\hline
3.5 & \cellcolor{ccc8_0!20} 0.054(7) & \cellcolor{ccc8_1!20} 0.079(9) & \cellcolor{ccc8_2!20} 0.112(10) & \cellcolor{ccc8_3!20} 0.140(12) & \cellcolor{ccc8_4!20} 0.143(12) & \cellcolor{ccc8_5!20} 0.172(14) & \cellcolor{ccc8_6!20} 0.171(14) & \cellcolor{ccc8_7!20} 0.233(18)\\
\hline
5.0 & \cellcolor{ccc9_0!20} 0.054(7) & \cellcolor{ccc9_1!20} 0.079(9) & \cellcolor{ccc9_2!20} 0.112(10) & \cellcolor{ccc9_3!20} 0.142(12) & \cellcolor{ccc9_4!20} 0.144(12) & \cellcolor{ccc9_5!20} 0.174(14) & \cellcolor{ccc9_6!20} 0.172(14) & \cellcolor{ccc9_7!20} 0.197(16)\\
\hline
7.0 & \cellcolor{ccc10_0!20} 0.054(7) & \cellcolor{ccc10_1!20} 0.079(9) & \cellcolor{ccc10_2!20} 0.112(10) & \cellcolor{ccc10_3!20} 0.142(12) & \cellcolor{ccc10_4!20} 0.144(12) & \cellcolor{ccc10_5!20} 0.175(14) & \cellcolor{ccc10_6!20} 0.172(14) & \cellcolor{ccc10_7!20} 0.197(16)\\
\hline
9.5 & \cellcolor{ccc11_0!20} 0.054(7) & \cellcolor{ccc11_1!20} 0.079(9) & \cellcolor{ccc11_2!20} 0.112(10) & \cellcolor{ccc11_3!20} 0.142(12) & \cellcolor{ccc11_4!20} 0.144(12) & \cellcolor{ccc11_5!20} 0.175(14) & \cellcolor{ccc11_6!20} 0.172(14) & \cellcolor{ccc11_7!20} 0.197(16)\\
\hline
\end{tabular}
}
\end{minipage}
\noindent
\hspace{-2.3cm}
\begin{minipage}[c]{0.3\textwidth}\hspace{0.1cm}$\theta=0.72$\\~\\
\resizebox{0.4\textwidth}{!}{
%
%
\begin{tikzpicture}

\begin{axis}[%
width=3.566in,
height=3.566in,
at={(1.236in,0.481in)},
scale only axis,
xmin=-1,
xmax=1,
xtick={\empty},
ymin=-1,
ymax=1,
ytick={\empty},
axis background/.style={fill=white},
axis x line*=bottom,
axis y line*=left,
legend style={legend cell align=left, align=left, draw=white!15!black}
]

\addplot[area legend, draw=black, fill=white]
table[row sep=crcr] {%
x	y\\
-1	-1\\
-1	1\\
-0.5	1\\
-0.5	-1\\
}--cycle;

\addplot[area legend, draw=black, fill=black!10]
table[row sep=crcr] {%
x	y\\
-0.5	-1\\
-0.5	1\\
0	1\\
0	-1\\
}--cycle;

\addplot[area legend, draw=black, fill=white]
table[row sep=crcr] {%
x	y\\
0	-1\\
0	1\\
0.5	1\\
0.5	-1\\
}--cycle;

\addplot[area legend, draw=black, fill=black!10]
table[row sep=crcr] {%
x	y\\
0.5	-1\\
0.5	1\\
1	1\\
1	-1\\
}--cycle;

\end{axis}
\end{tikzpicture}%
}
\end{minipage}

\caption{Computed values of the approximate convergence factor \textcolor{ultramarine}{$\rho$ and corresponding preconditioned CG iteration counts (between parenthesis)} w.r.t.\ parameters $\varepsilon$ on rows and mesh size $h$ on columns. In each table \textcolor{ultramarine}{the} pattern of diffusion coefficient $\mu$ and the strong threshold $\theta$ (shown on the left) is fixed. The background color depends on $\rho$ with colormap \protect\includegraphics[height=0.2cm]{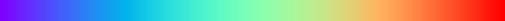}. It is scaled to range between \textcolor{reviewV2}{the} minimum and maximum value (among all the tables) of $\rho$. Strides pattern for $\mu$ as in Figure~\ref{fig:mu_pattern}(c).}
\label{tab:rho1}
\end{table}

\begin{table}[htbp]
\hspace{2cm}\hspace{-9.8cm}
\begin{minipage}[c]{0.7\textwidth}
\resizebox{2.5\textwidth}{!}{
\definecolor{ccc0_0}{rgb}{0.17058823529411765,0.4946558433997788,0.9667184042691874}
\definecolor{ccc0_1}{rgb}{0.303921568627451,0.30315267411304353,0.9881654720812594}
\definecolor{ccc0_2}{rgb}{0.3666666666666667,0.20791169081775931,0.9945218953682733}
\definecolor{ccc0_3}{rgb}{0.39803921568627454,0.1594757912099808,0.9967953249171991}
\definecolor{ccc0_4}{rgb}{0.3588235294117647,0.2199463578396686,0.9938591368952737}
\definecolor{ccc0_5}{rgb}{0.3509803921568627,0.23194764145389815,0.9931586661366362}
\definecolor{ccc0_6}{rgb}{0.3431372549019608,0.24391372010837714,0.9924205096719357}
\definecolor{ccc0_7}{rgb}{0.32745098039215687,0.2677330033224679,0.9908312530915603}
\definecolor{ccc1_0}{rgb}{0.2019607843137255,0.4512440570453228,0.9727282722446048}
\definecolor{ccc1_1}{rgb}{0.303921568627451,0.30315267411304353,0.9881654720812594}
\definecolor{ccc1_2}{rgb}{0.3431372549019608,0.24391372010837714,0.9924205096719357}
\definecolor{ccc1_3}{rgb}{0.32745098039215687,0.2677330033224679,0.9908312530915603}
\definecolor{ccc1_4}{rgb}{0.33529411764705885,0.2558427775944356,0.9916446955107427}
\definecolor{ccc1_5}{rgb}{0.3431372549019608,0.24391372010837714,0.9924205096719357}
\definecolor{ccc1_6}{rgb}{0.28823529411764703,0.3265387128400833,0.9862007473534026}
\definecolor{ccc1_7}{rgb}{0.2803921568627451,0.33815827481581706,0.9851622334675065}
\definecolor{ccc2_0}{rgb}{0.15490196078431373,0.5159178263577511,0.9634931442059831}
\definecolor{ccc2_1}{rgb}{0.11568627450980395,0.5676747161445901,0.9547913248866443}
\definecolor{ccc2_2}{rgb}{0.2019607843137255,0.4512440570453228,0.9727282722446048}
\definecolor{ccc2_3}{rgb}{0.13137254901960782,0.5472195469221112,0.9583812154701222}
\definecolor{ccc2_4}{rgb}{0.0607843137254902,0.6364742361471413,0.9410892525013715}
\definecolor{ccc2_5}{rgb}{0.009803921568627416,0.7179119230644418,0.9209055179449537}
\definecolor{ccc2_6}{rgb}{0.01764705882352935,0.7264335740162241,0.918486985745923}
\definecolor{ccc2_7}{rgb}{0.056862745098039236,0.767362681448697,0.9058734224033671}
\definecolor{ccc3_0}{rgb}{0.10392156862745094,0.8126223709664456,0.8896040127307095}
\definecolor{ccc3_1}{rgb}{0.2058823529411764,0.8951632913550623,0.8502171357296142}
\definecolor{ccc3_2}{rgb}{0.2529411764705882,0.9256376597815562,0.8301840308155507}
\definecolor{ccc3_3}{rgb}{0.2686274509803921,0.9346797673211106,0.8232529481575873}
\definecolor{ccc3_4}{rgb}{0.3313725490196078,0.965124085200289,0.7942898895752861}
\definecolor{ccc3_5}{rgb}{0.4019607843137255,0.9881654720812594,0.7594049166547071}
\definecolor{ccc3_6}{rgb}{0.39411764705882346,0.9862007473534026,0.763398282741103}
\definecolor{ccc3_7}{rgb}{0.4647058823529412,0.9984636039674339,0.7264335740162242}
\definecolor{ccc4_0}{rgb}{0.23725490196078436,0.9160336012803335,0.8369891083319778}
\definecolor{ccc4_1}{rgb}{0.40980392156862744,0.989980213280707,0.7553827347189938}
\definecolor{ccc4_2}{rgb}{0.5509803921568628,0.9967953249171991,0.678235117349234}
\definecolor{ccc4_3}{rgb}{0.5901960784313725,0.989980213280707,0.6552838500134537}
\definecolor{ccc4_4}{rgb}{0.6607843137254901,0.9682760409157589,0.6124202030492502}
\definecolor{ccc4_5}{rgb}{0.7470588235294118,0.9256376597815563,0.5574894393428855}
\definecolor{ccc4_6}{rgb}{0.8176470588235294,0.878081248083698,0.510631193180907}
\definecolor{ccc4_7}{rgb}{0.8725490196078431,0.8336023852211194,0.47309355683601007}
\definecolor{ccc5_0}{rgb}{0.30000000000000004,0.9510565162951535,0.8090169943749475}
\definecolor{ccc5_1}{rgb}{0.5509803921568628,0.9967953249171991,0.678235117349234}
\definecolor{ccc5_2}{rgb}{0.7,0.9510565162951536,0.5877852522924731}
\definecolor{ccc5_3}{rgb}{0.7784313725490195,0.9058734224033673,0.5368665976441802}
\definecolor{ccc5_4}{rgb}{0.9352941176470588,0.77520397611113,0.4291206087726091}
\definecolor{ccc5_5}{rgb}{1.0,0.6916984393193701,0.3727019919909141}
\definecolor{ccc5_6}{rgb}{1.0,0.5159178263577514,0.2677330033224681}
\definecolor{ccc5_7}{rgb}{1.0,0.37270199199091436,0.18980109344182594}
\definecolor{ccc6_0}{rgb}{0.3156862745098039,0.9583812154701222,0.8017142839800667}
\definecolor{ccc6_1}{rgb}{0.6294117647058823,0.9794097676013659,0.631711006253251}
\definecolor{ccc6_2}{rgb}{0.7392156862745098,0.9302293085467404,0.5625927516198231}
\definecolor{ccc6_3}{rgb}{0.9431372549019608,0.7673626814486971,0.4235485126679244}
\definecolor{ccc6_4}{rgb}{1.0,0.6552838500134535,0.3497265112062611}
\definecolor{ccc6_5}{rgb}{1.0,0.5159178263577514,0.2677330033224681}
\definecolor{ccc6_6}{rgb}{1.0,0.38410574917192575,0.1958454670071669}
\definecolor{ccc6_7}{rgb}{1.0,0.21994635783966857,0.11065268189150082}
\definecolor{ccc7_0}{rgb}{0.3313725490196078,0.965124085200289,0.7942898895752861}
\definecolor{ccc7_1}{rgb}{0.6372549019607843,0.9768483177596008,0.6269238058941066}
\definecolor{ccc7_2}{rgb}{0.8490196078431373,0.8534437988883159,0.4892929169339235}
\definecolor{ccc7_3}{rgb}{0.9666666666666666,0.7431448254773945,0.4067366430758004}
\definecolor{ccc7_4}{rgb}{1.0,0.6269238058941065,0.3323547994796596}
\definecolor{ccc7_5}{rgb}{1.0,0.4622038835403133,0.2379351950426188}
\definecolor{ccc7_6}{rgb}{1.0,0.3265387128400838,0.1655538761841302}
\definecolor{ccc7_7}{rgb}{1.0,0.1106526818915011,0.05541147491597008}
\definecolor{ccc8_0}{rgb}{0.40980392156862744,0.989980213280707,0.7553827347189938}
\definecolor{ccc8_1}{rgb}{0.6450980392156862,0.9741386021045102,0.622112816721474}
\definecolor{ccc8_2}{rgb}{0.8568627450980393,0.8469582108244671,0.4839114241003015}
\definecolor{ccc8_3}{rgb}{0.9823529411764707,0.7264335740162241,0.3954512068705426}
\definecolor{ccc8_4}{rgb}{1.0,0.5877852522924732,0.30901699437494745}
\definecolor{ccc8_5}{rgb}{1.0,0.4622038835403133,0.2379351950426188}
\definecolor{ccc8_6}{rgb}{1.0,0.279582592596744,0.14120615182309149}
\definecolor{ccc8_7}{rgb}{1.0,0.012319659535238529,0.006159946638138691}
\definecolor{ccc9_0}{rgb}{0.40980392156862744,0.989980213280707,0.7553827347189938}
\definecolor{ccc9_1}{rgb}{0.6529411764705881,0.9712810319161139,0.617278221289793}
\definecolor{ccc9_2}{rgb}{0.8568627450980393,0.8469582108244671,0.4839114241003015}
\definecolor{ccc9_3}{rgb}{0.9901960784313726,0.717911923064442,0.3897858732926794}
\definecolor{ccc9_4}{rgb}{1.0,0.5777738314082512,0.3031526741130436}
\definecolor{ccc9_5}{rgb}{1.0,0.4622038835403133,0.2379351950426188}
\definecolor{ccc9_6}{rgb}{1.0,0.279582592596744,0.14120615182309149}
\definecolor{ccc9_7}{rgb}{1.0,0.0,0.0}
\definecolor{ccc10_0}{rgb}{0.40980392156862744,0.989980213280707,0.7553827347189938}
\definecolor{ccc10_1}{rgb}{0.6529411764705881,0.9712810319161139,0.617278221289793}
\definecolor{ccc10_2}{rgb}{0.8568627450980393,0.8469582108244671,0.4839114241003015}
\definecolor{ccc10_3}{rgb}{0.9901960784313726,0.717911923064442,0.3897858732926794}
\definecolor{ccc10_4}{rgb}{1.0,0.5777738314082512,0.3031526741130436}
\definecolor{ccc10_5}{rgb}{1.0,0.4622038835403133,0.2379351950426188}
\definecolor{ccc10_6}{rgb}{1.0,0.279582592596744,0.14120615182309149}
\definecolor{ccc10_7}{rgb}{1.0,0.0,0.0}
\definecolor{ccc11_0}{rgb}{0.40980392156862744,0.989980213280707,0.7553827347189938}
\definecolor{ccc11_1}{rgb}{0.6529411764705881,0.9712810319161139,0.617278221289793}
\definecolor{ccc11_2}{rgb}{0.8568627450980393,0.8469582108244671,0.4839114241003015}
\definecolor{ccc11_3}{rgb}{0.9901960784313726,0.717911923064442,0.3897858732926794}
\definecolor{ccc11_4}{rgb}{1.0,0.5777738314082512,0.3031526741130436}
\definecolor{ccc11_5}{rgb}{1.0,0.4622038835403133,0.2379351950426188}
\definecolor{ccc11_6}{rgb}{1.0,0.279582592596744,0.14120615182309149}
\definecolor{ccc11_7}{rgb}{1.0,0.0,0.0}
\begin{tabular}{|c||c|c|c|c|c|c|c|c|}
\hline
$\varepsilon \backslash h$ & 1.25e-01& 6.25e-02& 3.12e-02& 1.56e-02& 7.81e-03& 3.91e-03& 1.95e-03& 9.77e-04\\
\hline\hline
0.0 & \cellcolor{ccc0_0!20} 0.094(9) & \cellcolor{ccc0_1!20} 0.071(8) & \cellcolor{ccc0_2!20} 0.060(8) & \cellcolor{ccc0_3!20} 0.054(8) & \cellcolor{ccc0_4!20} 0.061(9) & \cellcolor{ccc0_5!20} 0.063(9) & \cellcolor{ccc0_6!20} 0.064(9) & \cellcolor{ccc0_7!20} 0.066(10)\\
\hline
0.4 & \cellcolor{ccc1_0!20} 0.088(9) & \cellcolor{ccc1_1!20} 0.070(8) & \cellcolor{ccc1_2!20} 0.064(8) & \cellcolor{ccc1_3!20} 0.067(9) & \cellcolor{ccc1_4!20} 0.065(9) & \cellcolor{ccc1_5!20} 0.064(9) & \cellcolor{ccc1_6!20} 0.074(10) & \cellcolor{ccc1_7!20} 0.075(10)\\
\hline
0.8 & \cellcolor{ccc2_0!20} 0.097(9) & \cellcolor{ccc2_1!20} 0.103(10) & \cellcolor{ccc2_2!20} 0.088(9) & \cellcolor{ccc2_3!20} 0.101(10) & \cellcolor{ccc2_4!20} 0.113(11) & \cellcolor{ccc2_5!20} 0.125(12) & \cellcolor{ccc2_6!20} 0.127(12) & \cellcolor{ccc2_7!20} 0.134(13)\\
\hline
1.2 & \cellcolor{ccc3_0!20} 0.142(11) & \cellcolor{ccc3_1!20} 0.160(12) & \cellcolor{ccc3_2!20} 0.169(13) & \cellcolor{ccc3_3!20} 0.171(13) & \cellcolor{ccc3_4!20} 0.182(14) & \cellcolor{ccc3_5!20} 0.194(15) & \cellcolor{ccc3_6!20} 0.193(15) & \cellcolor{ccc3_7!20} 0.205(16)\\
\hline
1.6 & \cellcolor{ccc4_0!20} 0.166(12) & \cellcolor{ccc4_1!20} 0.196(13) & \cellcolor{ccc4_2!20} 0.220(15) & \cellcolor{ccc4_3!20} 0.228(16) & \cellcolor{ccc4_4!20} 0.240(17) & \cellcolor{ccc4_5!20} 0.255(18) & \cellcolor{ccc4_6!20} 0.268(19) & \cellcolor{ccc4_7!20} 0.277(20)\\
\hline
2.0 & \cellcolor{ccc5_0!20} 0.176(12) & \cellcolor{ccc5_1!20} 0.221(14) & \cellcolor{ccc5_2!20} 0.247(16) & \cellcolor{ccc5_3!20} 0.261(17) & \cellcolor{ccc5_4!20} 0.288(19) & \cellcolor{ccc5_5!20} 0.302(20) & \cellcolor{ccc5_6!20} 0.326(22) & \cellcolor{ccc5_7!20} 0.344(24)\\
\hline
2.4 & \cellcolor{ccc6_0!20} 0.180(12) & \cellcolor{ccc6_1!20} 0.234(15) & \cellcolor{ccc6_2!20} 0.254(16) & \cellcolor{ccc6_3!20} 0.289(19) & \cellcolor{ccc6_4!20} 0.307(20) & \cellcolor{ccc6_5!20} 0.326(22) & \cellcolor{ccc6_6!20} 0.343(24) & \cellcolor{ccc6_7!20} 0.362(25)\\
\hline
2.8 & \cellcolor{ccc7_0!20} 0.182(12) & \cellcolor{ccc7_1!20} 0.236(15) & \cellcolor{ccc7_2!20} 0.273(17) & \cellcolor{ccc7_3!20} 0.294(19) & \cellcolor{ccc7_4!20} 0.312(20) & \cellcolor{ccc7_5!20} 0.333(22) & \cellcolor{ccc7_6!20} 0.350(24) & \cellcolor{ccc7_7!20} 0.375(26)\\
\hline
3.5 & \cellcolor{ccc8_0!20} 0.196(13) & \cellcolor{ccc8_1!20} 0.237(15) & \cellcolor{ccc8_2!20} 0.275(17) & \cellcolor{ccc8_3!20} 0.297(19) & \cellcolor{ccc8_4!20} 0.317(21) & \cellcolor{ccc8_5!20} 0.333(23) & \cellcolor{ccc8_6!20} 0.355(24) & \cellcolor{ccc8_7!20} 0.386(27)\\
\hline
5.0 & \cellcolor{ccc9_0!20} 0.196(13) & \cellcolor{ccc9_1!20} 0.238(15) & \cellcolor{ccc9_2!20} 0.275(17) & \cellcolor{ccc9_3!20} 0.298(19) & \cellcolor{ccc9_4!20} 0.318(21) & \cellcolor{ccc9_5!20} 0.333(23) & \cellcolor{ccc9_6!20} 0.356(24) & \cellcolor{ccc9_7!20} 0.388(27)\\
\hline
7.0 & \cellcolor{ccc10_0!20} 0.196(13) & \cellcolor{ccc10_1!20} 0.238(15) & \cellcolor{ccc10_2!20} 0.275(17) & \cellcolor{ccc10_3!20} 0.298(19) & \cellcolor{ccc10_4!20} 0.318(21) & \cellcolor{ccc10_5!20} 0.333(23) & \cellcolor{ccc10_6!20} 0.356(24) & \cellcolor{ccc10_7!20} 0.388(27)\\
\hline
9.5 & \cellcolor{ccc11_0!20} 0.196(13) & \cellcolor{ccc11_1!20} 0.238(15) & \cellcolor{ccc11_2!20} 0.275(17) & \cellcolor{ccc11_3!20} 0.298(19) & \cellcolor{ccc11_4!20} 0.318(21) & \cellcolor{ccc11_5!20} 0.333(23) & \cellcolor{ccc11_6!20} 0.356(24) & \cellcolor{ccc11_7!20} 0.388(27)\\
\hline
\end{tabular}
}
\end{minipage}
\noindent
\hspace{-2.3cm}
\begin{minipage}[c]{0.3\textwidth}\hspace{0.1cm}$\theta=0.24$\\~\\
\resizebox{0.4\textwidth}{!}{
%
%
\begin{tikzpicture}

\begin{axis}[%
width=3.566in,
height=3.566in,
at={(1.236in,0.481in)},
scale only axis,
xmin=-1,
xmax=1,
xtick={\empty},
ymin=-1,
ymax=1,
ytick={\empty},
axis background/.style={fill=white},
axis x line*=bottom,
axis y line*=left,
legend style={legend cell align=left, align=left, draw=white!15!black}
]

\addplot[area legend, draw=black, fill=black!10]
table[row sep=crcr] {%
x	y\\
-1	-1\\
-1	-0.5\\
-0.5	-0.5\\
-0.5	-1\\
}--cycle;

\addplot[area legend, draw=black, fill=white]
table[row sep=crcr] {%
x	y\\
-1	-0.5\\
-1	0\\
-0.5	0\\
-0.5	-0.5\\
}--cycle;

\addplot[area legend, draw=black, fill=black!10]
table[row sep=crcr] {%
x	y\\
-1	0\\
-1	0.5\\
-0.5	0.5\\
-0.5	0\\
}--cycle;

\addplot[area legend, draw=black, fill=white]
table[row sep=crcr] {%
x	y\\
-1	0.5\\
-1	1\\
-0.5	1\\
-0.5	0.5\\
}--cycle;

\addplot[area legend, draw=black, fill=white]
table[row sep=crcr] {%
x	y\\
-0.5	-1\\
-0.5	-0.5\\
0	-0.5\\
0	-1\\
}--cycle;

\addplot[area legend, draw=black, fill=black!10]
table[row sep=crcr] {%
x	y\\
-0.5	-0.5\\
-0.5	0\\
0	0\\
0	-0.5\\
}--cycle;

\addplot[area legend, draw=black, fill=white]
table[row sep=crcr] {%
x	y\\
-0.5	0\\
-0.5	0.5\\
0	0.5\\
0	0\\
}--cycle;

\addplot[area legend, draw=black, fill=black!10]
table[row sep=crcr] {%
x	y\\
-0.5	0.5\\
-0.5	1\\
0	1\\
0	0.5\\
}--cycle;

\addplot[area legend, draw=black, fill=black!10]
table[row sep=crcr] {%
x	y\\
0	-1\\
0	-0.5\\
0.5	-0.5\\
0.5	-1\\
}--cycle;

\addplot[area legend, draw=black, fill=white]
table[row sep=crcr] {%
x	y\\
0	-0.5\\
0	0\\
0.5	0\\
0.5	-0.5\\
}--cycle;

\addplot[area legend, draw=black, fill=black!10]
table[row sep=crcr] {%
x	y\\
0	0\\
0	0.5\\
0.5	0.5\\
0.5	0\\
}--cycle;

\addplot[area legend, draw=black, fill=white]
table[row sep=crcr] {%
x	y\\
0	0.5\\
0	1\\
0.5	1\\
0.5	0.5\\
}--cycle;

\addplot[area legend, draw=black, fill=white]
table[row sep=crcr] {%
x	y\\
0.5	-1\\
0.5	-0.5\\
1	-0.5\\
1	-1\\
}--cycle;

\addplot[area legend, draw=black, fill=black!10]
table[row sep=crcr] {%
x	y\\
0.5	-0.5\\
0.5	0\\
1	0\\
1	-0.5\\
}--cycle;

\addplot[area legend, draw=black, fill=white]
table[row sep=crcr] {%
x	y\\
0.5	0\\
0.5	0.5\\
1	0.5\\
1	0\\
}--cycle;

\addplot[area legend, draw=black, fill=black!10]
table[row sep=crcr] {%
x	y\\
0.5	0.5\\
0.5	1\\
1	1\\
1	0.5\\
}--cycle;

\end{axis}
\end{tikzpicture}%
}
\end{minipage}\\

\vspace{1.5cm}\hspace{2cm}\hspace{-9.8cm}
\begin{minipage}[c]{0.7\textwidth}
\resizebox{2.5\textwidth}{!}{
\definecolor{ccc0_0}{rgb}{0.17058823529411765,0.4946558433997788,0.9667184042691874}
\definecolor{ccc0_1}{rgb}{0.303921568627451,0.30315267411304353,0.9881654720812594}
\definecolor{ccc0_2}{rgb}{0.3666666666666667,0.20791169081775931,0.9945218953682733}
\definecolor{ccc0_3}{rgb}{0.39803921568627454,0.1594757912099808,0.9967953249171991}
\definecolor{ccc0_4}{rgb}{0.3588235294117647,0.2199463578396686,0.9938591368952737}
\definecolor{ccc0_5}{rgb}{0.3509803921568627,0.23194764145389815,0.9931586661366362}
\definecolor{ccc0_6}{rgb}{0.3431372549019608,0.24391372010837714,0.9924205096719357}
\definecolor{ccc0_7}{rgb}{0.32745098039215687,0.2677330033224679,0.9908312530915603}
\definecolor{ccc1_0}{rgb}{0.2019607843137255,0.4512440570453228,0.9727282722446048}
\definecolor{ccc1_1}{rgb}{0.2647058823529412,0.3612416661871529,0.9829730996839018}
\definecolor{ccc1_2}{rgb}{0.3431372549019608,0.24391372010837714,0.9924205096719357}
\definecolor{ccc1_3}{rgb}{0.32745098039215687,0.2677330033224679,0.9908312530915603}
\definecolor{ccc1_4}{rgb}{0.33529411764705885,0.2558427775944356,0.9916446955107427}
\definecolor{ccc1_5}{rgb}{0.3509803921568627,0.23194764145389815,0.9931586661366362}
\definecolor{ccc1_6}{rgb}{0.3509803921568627,0.23194764145389815,0.9931586661366362}
\definecolor{ccc1_7}{rgb}{0.303921568627451,0.30315267411304353,0.9881654720812594}
\definecolor{ccc2_0}{rgb}{0.15490196078431373,0.5159178263577511,0.9634931442059831}
\definecolor{ccc2_1}{rgb}{0.11568627450980395,0.5676747161445901,0.9547913248866443}
\definecolor{ccc2_2}{rgb}{0.2019607843137255,0.4512440570453228,0.9727282722446048}
\definecolor{ccc2_3}{rgb}{0.13137254901960782,0.5472195469221112,0.9583812154701222}
\definecolor{ccc2_4}{rgb}{0.0607843137254902,0.6364742361471413,0.9410892525013715}
\definecolor{ccc2_5}{rgb}{0.009803921568627416,0.7179119230644418,0.9209055179449537}
\definecolor{ccc2_6}{rgb}{0.01764705882352935,0.7264335740162241,0.918486985745923}
\definecolor{ccc2_7}{rgb}{0.056862745098039236,0.767362681448697,0.9058734224033671}
\definecolor{ccc3_0}{rgb}{0.0607843137254902,0.6364742361471413,0.9410892525013715}
\definecolor{ccc3_1}{rgb}{0.21372549019607845,0.9005867023006374,0.8469582108244671}
\definecolor{ccc3_2}{rgb}{0.09999999999999998,0.5877852522924731,0.9510565162951535}
\definecolor{ccc3_3}{rgb}{0.03725490196078429,0.664540178707858,0.9346797673211106}
\definecolor{ccc3_4}{rgb}{0.0725490196078431,0.7829276104921027,0.9005867023006374}
\definecolor{ccc3_5}{rgb}{0.19019607843137254,0.883909710213612,0.8566380778638628}
\definecolor{ccc3_6}{rgb}{0.3392156862745098,0.9682760409157589,0.7905324123001634}
\definecolor{ccc3_7}{rgb}{0.3156862745098039,0.9583812154701222,0.8017142839800667}
\definecolor{ccc4_0}{rgb}{0.025490196078431393,0.7348449670469757,0.9160336012803335}
\definecolor{ccc4_1}{rgb}{0.38627450980392153,0.9840863373026044,0.767362681448697}
\definecolor{ccc4_2}{rgb}{0.03725490196078429,0.664540178707858,0.9346797673211106}
\definecolor{ccc4_3}{rgb}{0.0490196078431373,0.7594049166547072,0.9084652718195236}
\definecolor{ccc4_4}{rgb}{0.08823529411764708,0.7980172272802395,0.8951632913550623}
\definecolor{ccc4_5}{rgb}{0.23725490196078436,0.9160336012803335,0.8369891083319778}
\definecolor{ccc4_6}{rgb}{0.3784313725490196,0.9818225628535369,0.7712979623471807}
\definecolor{ccc4_7}{rgb}{0.4647058823529412,0.9984636039674339,0.7264335740162242}
\definecolor{ccc5_0}{rgb}{0.12745098039215685,0.8336023852211195,0.8810121942857845}
\definecolor{ccc5_1}{rgb}{0.4882352941176471,0.9998292504580527,0.7136101544117524}
\definecolor{ccc5_2}{rgb}{0.0019607843137254832,0.7092813076058535,0.9232891061054894}
\definecolor{ccc5_3}{rgb}{0.009803921568627416,0.7179119230644418,0.9209055179449537}
\definecolor{ccc5_4}{rgb}{0.1588235294117647,0.8597998514483723,0.8690889463055284}
\definecolor{ccc5_5}{rgb}{0.2529411764705882,0.9256376597815562,0.8301840308155507}
\definecolor{ccc5_6}{rgb}{0.4019607843137255,0.9881654720812594,0.7594049166547071}
\definecolor{ccc5_7}{rgb}{0.5745098039215686,0.9931586661366362,0.6645401787078581}
\definecolor{ccc6_0}{rgb}{0.15098039215686276,0.8534437988883159,0.8721195109836108}
\definecolor{ccc6_1}{rgb}{0.5274509803921568,0.9990704811844932,0.69169843931937}
\definecolor{ccc6_2}{rgb}{0.01764705882352935,0.7264335740162241,0.918486985745923}
\definecolor{ccc6_3}{rgb}{0.09607843137254901,0.8053809193888326,0.8924005832479478}
\definecolor{ccc6_4}{rgb}{0.17450980392156867,0.8721195109836108,0.8629289996673897}
\definecolor{ccc6_5}{rgb}{0.26078431372549016,0.9302293085467402,0.8267341748257635}
\definecolor{ccc6_6}{rgb}{0.3392156862745098,0.9682760409157589,0.7905324123001634}
\definecolor{ccc6_7}{rgb}{0.7313725490196079,0.9346797673211106,0.5676747161445901}
\definecolor{ccc7_0}{rgb}{0.1588235294117647,0.8597998514483723,0.8690889463055284}
\definecolor{ccc7_1}{rgb}{0.5352941176470587,0.9984636039674339,0.6872366859692628}
\definecolor{ccc7_2}{rgb}{0.009803921568627416,0.7179119230644418,0.9209055179449537}
\definecolor{ccc7_3}{rgb}{0.08823529411764708,0.7980172272802395,0.8951632913550623}
\definecolor{ccc7_4}{rgb}{0.11176470588235299,0.8197404829072211,0.8867736859200619}
\definecolor{ccc7_5}{rgb}{0.2686274509803921,0.9346797673211106,0.8232529481575873}
\definecolor{ccc7_6}{rgb}{0.40980392156862744,0.989980213280707,0.7553827347189938}
\definecolor{ccc7_7}{rgb}{0.7313725490196079,0.9346797673211106,0.5676747161445901}
\definecolor{ccc8_0}{rgb}{0.16666666666666663,0.8660254037844386,0.8660254037844387}
\definecolor{ccc8_1}{rgb}{0.5509803921568628,0.9967953249171991,0.678235117349234}
\definecolor{ccc8_2}{rgb}{0.009803921568627416,0.7179119230644418,0.9209055179449537}
\definecolor{ccc8_3}{rgb}{0.08039215686274503,0.7905324123001632,0.8978920322202582}
\definecolor{ccc8_4}{rgb}{0.1588235294117647,0.8597998514483723,0.8690889463055284}
\definecolor{ccc8_5}{rgb}{0.40980392156862744,0.989980213280707,0.7553827347189938}
\definecolor{ccc8_6}{rgb}{0.4490196078431372,0.9967953249171991,0.7348449670469757}
\definecolor{ccc8_7}{rgb}{0.5431372549019609,0.9977051801738729,0.6827488552151855}
\definecolor{ccc9_0}{rgb}{0.16666666666666663,0.8660254037844386,0.8660254037844387}
\definecolor{ccc9_1}{rgb}{0.5509803921568628,0.9967953249171991,0.678235117349234}
\definecolor{ccc9_2}{rgb}{0.009803921568627416,0.7179119230644418,0.9209055179449537}
\definecolor{ccc9_3}{rgb}{0.08039215686274503,0.7905324123001632,0.8978920322202582}
\definecolor{ccc9_4}{rgb}{0.16666666666666663,0.8660254037844386,0.8660254037844387}
\definecolor{ccc9_5}{rgb}{0.40980392156862744,0.989980213280707,0.7553827347189938}
\definecolor{ccc9_6}{rgb}{0.49607843137254903,0.9999810273487268,0.7092813076058535}
\definecolor{ccc9_7}{rgb}{0.5745098039215686,0.9931586661366362,0.6645401787078581}
\definecolor{ccc10_0}{rgb}{0.16666666666666663,0.8660254037844386,0.8660254037844387}
\definecolor{ccc10_1}{rgb}{0.5509803921568628,0.9967953249171991,0.678235117349234}
\definecolor{ccc10_2}{rgb}{0.009803921568627416,0.7179119230644418,0.9209055179449537}
\definecolor{ccc10_3}{rgb}{0.08039215686274503,0.7905324123001632,0.8978920322202582}
\definecolor{ccc10_4}{rgb}{0.16666666666666663,0.8660254037844386,0.8660254037844387}
\definecolor{ccc10_5}{rgb}{0.40980392156862744,0.989980213280707,0.7553827347189938}
\definecolor{ccc10_6}{rgb}{0.49607843137254903,0.9999810273487268,0.7092813076058535}
\definecolor{ccc10_7}{rgb}{0.5745098039215686,0.9931586661366362,0.6645401787078581}
\definecolor{ccc11_0}{rgb}{0.16666666666666663,0.8660254037844386,0.8660254037844387}
\definecolor{ccc11_1}{rgb}{0.5509803921568628,0.9967953249171991,0.678235117349234}
\definecolor{ccc11_2}{rgb}{0.009803921568627416,0.7179119230644418,0.9209055179449537}
\definecolor{ccc11_3}{rgb}{0.08039215686274503,0.7905324123001632,0.8978920322202582}
\definecolor{ccc11_4}{rgb}{0.16666666666666663,0.8660254037844386,0.8660254037844387}
\definecolor{ccc11_5}{rgb}{0.40980392156862744,0.989980213280707,0.7553827347189938}
\definecolor{ccc11_6}{rgb}{0.49607843137254903,0.9999810273487268,0.7092813076058535}
\definecolor{ccc11_7}{rgb}{0.5745098039215686,0.9931586661366362,0.6645401787078581}
\begin{tabular}{|c||c|c|c|c|c|c|c|c|}
\hline
$\varepsilon \backslash h$ & 1.25e-01& 6.25e-02& 3.12e-02& 1.56e-02& 7.81e-03& 3.91e-03& 1.95e-03& 9.77e-04\\
\hline\hline
0.0 & \cellcolor{ccc0_0!20} 0.094(9) & \cellcolor{ccc0_1!20} 0.071(8) & \cellcolor{ccc0_2!20} 0.060(8) & \cellcolor{ccc0_3!20} 0.054(8) & \cellcolor{ccc0_4!20} 0.061(9) & \cellcolor{ccc0_5!20} 0.063(9) & \cellcolor{ccc0_6!20} 0.064(9) & \cellcolor{ccc0_7!20} 0.066(10)\\
\hline
0.4 & \cellcolor{ccc1_0!20} 0.088(9) & \cellcolor{ccc1_1!20} 0.077(9) & \cellcolor{ccc1_2!20} 0.064(8) & \cellcolor{ccc1_3!20} 0.066(9) & \cellcolor{ccc1_4!20} 0.065(9) & \cellcolor{ccc1_5!20} 0.063(9) & \cellcolor{ccc1_6!20} 0.063(9) & \cellcolor{ccc1_7!20} 0.071(10)\\
\hline
0.8 & \cellcolor{ccc2_0!20} 0.097(9) & \cellcolor{ccc2_1!20} 0.103(10) & \cellcolor{ccc2_2!20} 0.088(9) & \cellcolor{ccc2_3!20} 0.101(10) & \cellcolor{ccc2_4!20} 0.113(11) & \cellcolor{ccc2_5!20} 0.125(12) & \cellcolor{ccc2_6!20} 0.127(12) & \cellcolor{ccc2_7!20} 0.134(13)\\
\hline
1.2 & \cellcolor{ccc3_0!20} 0.113(10) & \cellcolor{ccc3_1!20} 0.161(12) & \cellcolor{ccc3_2!20} 0.107(10) & \cellcolor{ccc3_3!20} 0.117(11) & \cellcolor{ccc3_4!20} 0.137(12) & \cellcolor{ccc3_5!20} 0.158(13) & \cellcolor{ccc3_6!20} 0.184(15) & \cellcolor{ccc3_7!20} 0.180(15)\\
\hline
1.6 & \cellcolor{ccc4_0!20} 0.129(10) & \cellcolor{ccc4_1!20} 0.192(13) & \cellcolor{ccc4_2!20} 0.118(11) & \cellcolor{ccc4_3!20} 0.132(12) & \cellcolor{ccc4_4!20} 0.140(12) & \cellcolor{ccc4_5!20} 0.166(14) & \cellcolor{ccc4_6!20} 0.191(15) & \cellcolor{ccc4_7!20} 0.205(16)\\
\hline
2.0 & \cellcolor{ccc5_0!20} 0.147(11) & \cellcolor{ccc5_1!20} 0.210(14) & \cellcolor{ccc5_2!20} 0.124(11) & \cellcolor{ccc5_3!20} 0.126(12) & \cellcolor{ccc5_4!20} 0.152(13) & \cellcolor{ccc5_5!20} 0.169(14) & \cellcolor{ccc5_6!20} 0.194(15) & \cellcolor{ccc5_7!20} 0.225(18)\\
\hline
2.4 & \cellcolor{ccc6_0!20} 0.150(11) & \cellcolor{ccc6_1!20} 0.216(14) & \cellcolor{ccc6_2!20} 0.127(11) & \cellcolor{ccc6_3!20} 0.141(12) & \cellcolor{ccc6_4!20} 0.155(13) & \cellcolor{ccc6_5!20} 0.170(14) & \cellcolor{ccc6_6!20} 0.184(15) & \cellcolor{ccc6_7!20} 0.252(19)\\
\hline
2.8 & \cellcolor{ccc7_0!20} 0.152(11) & \cellcolor{ccc7_1!20} 0.218(14) & \cellcolor{ccc7_2!20} 0.125(11) & \cellcolor{ccc7_3!20} 0.140(12) & \cellcolor{ccc7_4!20} 0.144(12) & \cellcolor{ccc7_5!20} 0.171(14) & \cellcolor{ccc7_6!20} 0.196(16) & \cellcolor{ccc7_7!20} 0.252(19)\\
\hline
3.5 & \cellcolor{ccc8_0!20} 0.153(11) & \cellcolor{ccc8_1!20} 0.220(14) & \cellcolor{ccc8_2!20} 0.126(11) & \cellcolor{ccc8_3!20} 0.138(12) & \cellcolor{ccc8_4!20} 0.152(13) & \cellcolor{ccc8_5!20} 0.196(15) & \cellcolor{ccc8_6!20} 0.203(16) & \cellcolor{ccc8_7!20} 0.219(17)\\
\hline
5.0 & \cellcolor{ccc9_0!20} 0.153(11) & \cellcolor{ccc9_1!20} 0.220(14) & \cellcolor{ccc9_2!20} 0.126(11) & \cellcolor{ccc9_3!20} 0.138(12) & \cellcolor{ccc9_4!20} 0.153(13) & \cellcolor{ccc9_5!20} 0.196(15) & \cellcolor{ccc9_6!20} 0.211(16) & \cellcolor{ccc9_7!20} 0.225(17)\\
\hline
7.0 & \cellcolor{ccc10_0!20} 0.153(11) & \cellcolor{ccc10_1!20} 0.220(14) & \cellcolor{ccc10_2!20} 0.126(11) & \cellcolor{ccc10_3!20} 0.138(12) & \cellcolor{ccc10_4!20} 0.153(13) & \cellcolor{ccc10_5!20} 0.196(15) & \cellcolor{ccc10_6!20} 0.211(16) & \cellcolor{ccc10_7!20} 0.225(17)\\
\hline
9.5 & \cellcolor{ccc11_0!20} 0.153(11) & \cellcolor{ccc11_1!20} 0.220(14) & \cellcolor{ccc11_2!20} 0.126(11) & \cellcolor{ccc11_3!20} 0.138(12) & \cellcolor{ccc11_4!20} 0.153(13) & \cellcolor{ccc11_5!20} 0.196(15) & \cellcolor{ccc11_6!20} 0.211(16) & \cellcolor{ccc11_7!20} 0.225(17)\\
\hline
\end{tabular}
}
\end{minipage}
\noindent
\hspace{-2.3cm}
\begin{minipage}[c]{0.3\textwidth}\hspace{0.1cm}$\theta=0.48$\\~\\
\resizebox{0.4\textwidth}{!}{
%
%
\begin{tikzpicture}

\begin{axis}[%
width=3.566in,
height=3.566in,
at={(1.236in,0.481in)},
scale only axis,
xmin=-1,
xmax=1,
xtick={\empty},
ymin=-1,
ymax=1,
ytick={\empty},
axis background/.style={fill=white},
axis x line*=bottom,
axis y line*=left,
legend style={legend cell align=left, align=left, draw=white!15!black}
]

\addplot[area legend, draw=black, fill=black!10]
table[row sep=crcr] {%
x	y\\
-1	-1\\
-1	-0.5\\
-0.5	-0.5\\
-0.5	-1\\
}--cycle;

\addplot[area legend, draw=black, fill=white]
table[row sep=crcr] {%
x	y\\
-1	-0.5\\
-1	0\\
-0.5	0\\
-0.5	-0.5\\
}--cycle;

\addplot[area legend, draw=black, fill=black!10]
table[row sep=crcr] {%
x	y\\
-1	0\\
-1	0.5\\
-0.5	0.5\\
-0.5	0\\
}--cycle;

\addplot[area legend, draw=black, fill=white]
table[row sep=crcr] {%
x	y\\
-1	0.5\\
-1	1\\
-0.5	1\\
-0.5	0.5\\
}--cycle;

\addplot[area legend, draw=black, fill=white]
table[row sep=crcr] {%
x	y\\
-0.5	-1\\
-0.5	-0.5\\
0	-0.5\\
0	-1\\
}--cycle;

\addplot[area legend, draw=black, fill=black!10]
table[row sep=crcr] {%
x	y\\
-0.5	-0.5\\
-0.5	0\\
0	0\\
0	-0.5\\
}--cycle;

\addplot[area legend, draw=black, fill=white]
table[row sep=crcr] {%
x	y\\
-0.5	0\\
-0.5	0.5\\
0	0.5\\
0	0\\
}--cycle;

\addplot[area legend, draw=black, fill=black!10]
table[row sep=crcr] {%
x	y\\
-0.5	0.5\\
-0.5	1\\
0	1\\
0	0.5\\
}--cycle;

\addplot[area legend, draw=black, fill=black!10]
table[row sep=crcr] {%
x	y\\
0	-1\\
0	-0.5\\
0.5	-0.5\\
0.5	-1\\
}--cycle;

\addplot[area legend, draw=black, fill=white]
table[row sep=crcr] {%
x	y\\
0	-0.5\\
0	0\\
0.5	0\\
0.5	-0.5\\
}--cycle;

\addplot[area legend, draw=black, fill=black!10]
table[row sep=crcr] {%
x	y\\
0	0\\
0	0.5\\
0.5	0.5\\
0.5	0\\
}--cycle;

\addplot[area legend, draw=black, fill=white]
table[row sep=crcr] {%
x	y\\
0	0.5\\
0	1\\
0.5	1\\
0.5	0.5\\
}--cycle;

\addplot[area legend, draw=black, fill=white]
table[row sep=crcr] {%
x	y\\
0.5	-1\\
0.5	-0.5\\
1	-0.5\\
1	-1\\
}--cycle;

\addplot[area legend, draw=black, fill=black!10]
table[row sep=crcr] {%
x	y\\
0.5	-0.5\\
0.5	0\\
1	0\\
1	-0.5\\
}--cycle;

\addplot[area legend, draw=black, fill=white]
table[row sep=crcr] {%
x	y\\
0.5	0\\
0.5	0.5\\
1	0.5\\
1	0\\
}--cycle;

\addplot[area legend, draw=black, fill=black!10]
table[row sep=crcr] {%
x	y\\
0.5	0.5\\
0.5	1\\
1	1\\
1	0.5\\
}--cycle;

\end{axis}
\end{tikzpicture}%
}
\end{minipage}\\

\vspace{1.5cm}\hspace{2cm}\hspace{-9.8cm}
\begin{minipage}[c]{0.7\textwidth}
\resizebox{2.5\textwidth}{!}{
\definecolor{ccc0_0}{rgb}{0.17058823529411765,0.4946558433997788,0.9667184042691874}
\definecolor{ccc0_1}{rgb}{0.303921568627451,0.30315267411304353,0.9881654720812594}
\definecolor{ccc0_2}{rgb}{0.3666666666666667,0.20791169081775931,0.9945218953682733}
\definecolor{ccc0_3}{rgb}{0.39803921568627454,0.1594757912099808,0.9967953249171991}
\definecolor{ccc0_4}{rgb}{0.3588235294117647,0.2199463578396686,0.9938591368952737}
\definecolor{ccc0_5}{rgb}{0.3509803921568627,0.23194764145389815,0.9931586661366362}
\definecolor{ccc0_6}{rgb}{0.3431372549019608,0.24391372010837714,0.9924205096719357}
\definecolor{ccc0_7}{rgb}{0.32745098039215687,0.2677330033224679,0.9908312530915603}
\definecolor{ccc1_0}{rgb}{0.4843137254901961,0.02463744919538197,0.9999241101148306}
\definecolor{ccc1_1}{rgb}{0.09607843137254901,0.8053809193888326,0.8924005832479478}
\definecolor{ccc1_2}{rgb}{0.11176470588235299,0.8197404829072211,0.8867736859200619}
\definecolor{ccc1_3}{rgb}{0.27647058823529413,0.9389883606150565,0.8197404829072211}
\definecolor{ccc1_4}{rgb}{0.3392156862745098,0.9682760409157589,0.7905324123001634}
\definecolor{ccc1_5}{rgb}{0.4019607843137255,0.9881654720812594,0.7594049166547071}
\definecolor{ccc1_6}{rgb}{0.48039215686274506,0.9995257197133659,0.717911923064442}
\definecolor{ccc1_7}{rgb}{0.5196078431372548,0.9995257197133659,0.6961339459629267}
\definecolor{ccc2_0}{rgb}{0.45294117647058824,0.07385252747487396,0.9993170601430229}
\definecolor{ccc2_1}{rgb}{0.02941176470588236,0.6736956436465572,0.9324722294043558}
\definecolor{ccc2_2}{rgb}{0.04117647058823526,0.7513318895568732,0.9110226492460883}
\definecolor{ccc2_3}{rgb}{0.27647058823529413,0.9389883606150565,0.8197404829072211}
\definecolor{ccc2_4}{rgb}{0.3784313725490196,0.9818225628535369,0.7712979623471807}
\definecolor{ccc2_5}{rgb}{0.503921568627451,0.9999810273487268,0.7049255469061472}
\definecolor{ccc2_6}{rgb}{0.503921568627451,0.9999810273487268,0.7049255469061472}
\definecolor{ccc2_7}{rgb}{0.6607843137254901,0.9682760409157589,0.6124202030492502}
\definecolor{ccc3_0}{rgb}{0.4137254901960784,0.13510524748139296,0.9977051801738729}
\definecolor{ccc3_1}{rgb}{0.05294117647058827,0.6459280624867872,0.9389883606150565}
\definecolor{ccc3_2}{rgb}{0.292156862745098,0.9471773565640402,0.8126223709664456}
\definecolor{ccc3_3}{rgb}{0.19803921568627447,0.8896040127307095,0.853443798888316}
\definecolor{ccc3_4}{rgb}{0.4411764705882353,0.9957341762950345,0.7390089172206591}
\definecolor{ccc3_5}{rgb}{0.6607843137254901,0.9682760409157589,0.6124202030492502}
\definecolor{ccc3_6}{rgb}{0.6764705882352942,0.961825643172819,0.6026346363792564}
\definecolor{ccc3_7}{rgb}{0.8098039215686275,0.8839097102136121,0.5159178263577511}
\definecolor{ccc4_0}{rgb}{0.38235294117647056,0.18374951781657034,0.9957341762950345}
\definecolor{ccc4_1}{rgb}{0.09215686274509804,0.5977074592660936,0.9491349440359013}
\definecolor{ccc4_2}{rgb}{0.1588235294117647,0.8597998514483723,0.8690889463055284}
\definecolor{ccc4_3}{rgb}{0.2686274509803921,0.9346797673211106,0.8232529481575873}
\definecolor{ccc4_4}{rgb}{0.503921568627451,0.9999810273487268,0.7049255469061472}
\definecolor{ccc4_5}{rgb}{0.6843137254901961,0.9583812154701222,0.5977074592660936}
\definecolor{ccc4_6}{rgb}{0.8333333333333333,0.8660254037844387,0.5000000000000001}
\definecolor{ccc4_7}{rgb}{0.8568627450980393,0.8469582108244671,0.4839114241003015}
\definecolor{ccc5_0}{rgb}{0.4137254901960784,0.13510524748139296,0.9977051801738729}
\definecolor{ccc5_1}{rgb}{0.05294117647058827,0.6459280624867872,0.9389883606150565}
\definecolor{ccc5_2}{rgb}{0.19019607843137254,0.883909710213612,0.8566380778638628}
\definecolor{ccc5_3}{rgb}{0.2529411764705882,0.9256376597815562,0.8301840308155507}
\definecolor{ccc5_4}{rgb}{0.4647058823529412,0.9984636039674339,0.7264335740162242}
\definecolor{ccc5_5}{rgb}{0.5588235294117647,0.9957341762950346,0.6736956436465572}
\definecolor{ccc5_6}{rgb}{0.7941176470588236,0.8951632913550623,0.5264321628773558}
\definecolor{ccc5_7}{rgb}{0.8098039215686275,0.8839097102136121,0.5159178263577511}
\definecolor{ccc6_0}{rgb}{0.4137254901960784,0.13510524748139296,0.9977051801738729}
\definecolor{ccc6_1}{rgb}{0.05294117647058827,0.6459280624867872,0.9389883606150565}
\definecolor{ccc6_2}{rgb}{0.12745098039215685,0.8336023852211195,0.8810121942857845}
\definecolor{ccc6_3}{rgb}{0.16666666666666663,0.8660254037844386,0.8660254037844387}
\definecolor{ccc6_4}{rgb}{0.3549019607843137,0.9741386021045101,0.7829276104921028}
\definecolor{ccc6_5}{rgb}{0.5509803921568628,0.9967953249171991,0.678235117349234}
\definecolor{ccc6_6}{rgb}{0.6843137254901961,0.9583812154701222,0.5977074592660936}
\definecolor{ccc6_7}{rgb}{0.8490196078431373,0.8534437988883159,0.4892929169339235}
\definecolor{ccc7_0}{rgb}{0.4137254901960784,0.13510524748139296,0.9977051801738729}
\definecolor{ccc7_1}{rgb}{0.04509803921568628,0.6552838500134536,0.9368518385313106}
\definecolor{ccc7_2}{rgb}{0.12745098039215685,0.8336023852211195,0.8810121942857845}
\definecolor{ccc7_3}{rgb}{0.21372549019607845,0.9005867023006374,0.8469582108244671}
\definecolor{ccc7_4}{rgb}{0.36274509803921573,0.9768483177596007,0.7790805745256704}
\definecolor{ccc7_5}{rgb}{0.5666666666666667,0.9945218953682734,0.6691306063588582}
\definecolor{ccc7_6}{rgb}{0.6686274509803922,0.965124085200289,0.6075389463388169}
\definecolor{ccc7_7}{rgb}{0.7941176470588236,0.8951632913550623,0.5264321628773558}
\definecolor{ccc8_0}{rgb}{0.40588235294117647,0.1473016980546375,0.997269173385788}
\definecolor{ccc8_1}{rgb}{0.03725490196078429,0.664540178707858,0.9346797673211106}
\definecolor{ccc8_2}{rgb}{0.09607843137254901,0.8053809193888326,0.8924005832479478}
\definecolor{ccc8_3}{rgb}{0.21372549019607845,0.9005867023006374,0.8469582108244671}
\definecolor{ccc8_4}{rgb}{0.36274509803921573,0.9768483177596007,0.7790805745256704}
\definecolor{ccc8_5}{rgb}{0.6294117647058823,0.9794097676013659,0.631711006253251}
\definecolor{ccc8_6}{rgb}{0.6372549019607843,0.9768483177596008,0.6269238058941066}
\definecolor{ccc8_7}{rgb}{0.7784313725490195,0.9058734224033673,0.5368665976441802}
\definecolor{ccc9_0}{rgb}{0.40588235294117647,0.1473016980546375,0.997269173385788}
\definecolor{ccc9_1}{rgb}{0.03725490196078429,0.664540178707858,0.9346797673211106}
\definecolor{ccc9_2}{rgb}{0.09607843137254901,0.8053809193888326,0.8924005832479478}
\definecolor{ccc9_3}{rgb}{0.21372549019607845,0.9005867023006374,0.8469582108244671}
\definecolor{ccc9_4}{rgb}{0.37058823529411766,0.9794097676013659,0.7752039761111298}
\definecolor{ccc9_5}{rgb}{0.6843137254901961,0.9583812154701222,0.5977074592660936}
\definecolor{ccc9_6}{rgb}{0.6372549019607843,0.9768483177596008,0.6269238058941066}
\definecolor{ccc9_7}{rgb}{0.7705882352941176,0.9110226492460884,0.5420533564724495}
\definecolor{ccc10_0}{rgb}{0.40588235294117647,0.1473016980546375,0.997269173385788}
\definecolor{ccc10_1}{rgb}{0.021568627450980427,0.6827488552151855,0.9302293085467404}
\definecolor{ccc10_2}{rgb}{0.09607843137254901,0.8053809193888326,0.8924005832479478}
\definecolor{ccc10_3}{rgb}{0.21372549019607845,0.9005867023006374,0.8469582108244671}
\definecolor{ccc10_4}{rgb}{0.37058823529411766,0.9794097676013659,0.7752039761111298}
\definecolor{ccc10_5}{rgb}{0.692156862745098,0.9547913248866443,0.5927576019625549}
\definecolor{ccc10_6}{rgb}{0.6372549019607843,0.9768483177596008,0.6269238058941066}
\definecolor{ccc10_7}{rgb}{0.7705882352941176,0.9110226492460884,0.5420533564724495}
\definecolor{ccc11_0}{rgb}{0.40588235294117647,0.1473016980546375,0.997269173385788}
\definecolor{ccc11_1}{rgb}{0.021568627450980427,0.6827488552151855,0.9302293085467404}
\definecolor{ccc11_2}{rgb}{0.09607843137254901,0.8053809193888326,0.8924005832479478}
\definecolor{ccc11_3}{rgb}{0.21372549019607845,0.9005867023006374,0.8469582108244671}
\definecolor{ccc11_4}{rgb}{0.37058823529411766,0.9794097676013659,0.7752039761111298}
\definecolor{ccc11_5}{rgb}{0.692156862745098,0.9547913248866443,0.5927576019625549}
\definecolor{ccc11_6}{rgb}{0.6372549019607843,0.9768483177596008,0.6269238058941066}
\definecolor{ccc11_7}{rgb}{0.7705882352941176,0.9110226492460884,0.5420533564724495}
\begin{tabular}{|c||c|c|c|c|c|c|c|c|}
\hline
$\varepsilon \backslash h$ & 1.25e-01& 6.25e-02& 3.12e-02& 1.56e-02& 7.81e-03& 3.91e-03& 1.95e-03& 9.77e-04\\
\hline\hline
0.0 & \cellcolor{ccc0_0!20} 0.094(9) & \cellcolor{ccc0_1!20} 0.071(8) & \cellcolor{ccc0_2!20} 0.060(8) & \cellcolor{ccc0_3!20} 0.054(8) & \cellcolor{ccc0_4!20} 0.061(9) & \cellcolor{ccc0_5!20} 0.063(9) & \cellcolor{ccc0_6!20} 0.064(9) & \cellcolor{ccc0_7!20} 0.066(10)\\
\hline
0.4 & \cellcolor{ccc1_0!20} 0.039(7) & \cellcolor{ccc1_1!20} 0.141(11) & \cellcolor{ccc1_2!20} 0.143(12) & \cellcolor{ccc1_3!20} 0.172(13) & \cellcolor{ccc1_4!20} 0.183(14) & \cellcolor{ccc1_5!20} 0.195(15) & \cellcolor{ccc1_6!20} 0.208(16) & \cellcolor{ccc1_7!20} 0.215(17)\\
\hline
0.8 & \cellcolor{ccc2_0!20} 0.045(7) & \cellcolor{ccc2_1!20} 0.119(10) & \cellcolor{ccc2_2!20} 0.131(11) & \cellcolor{ccc2_3!20} 0.173(13) & \cellcolor{ccc2_4!20} 0.190(14) & \cellcolor{ccc2_5!20} 0.213(16) & \cellcolor{ccc2_6!20} 0.212(16) & \cellcolor{ccc2_7!20} 0.240(18)\\
\hline
1.2 & \cellcolor{ccc3_0!20} 0.052(7) & \cellcolor{ccc3_1!20} 0.115(10) & \cellcolor{ccc3_2!20} 0.175(13) & \cellcolor{ccc3_3!20} 0.159(13) & \cellcolor{ccc3_4!20} 0.202(15) & \cellcolor{ccc3_5!20} 0.240(17) & \cellcolor{ccc3_6!20} 0.243(18) & \cellcolor{ccc3_7!20} 0.266(20)\\
\hline
1.6 & \cellcolor{ccc4_0!20} 0.057(8) & \cellcolor{ccc4_1!20} 0.108(10) & \cellcolor{ccc4_2!20} 0.152(12) & \cellcolor{ccc4_3!20} 0.171(13) & \cellcolor{ccc4_4!20} 0.212(15) & \cellcolor{ccc4_5!20} 0.244(17) & \cellcolor{ccc4_6!20} 0.270(19) & \cellcolor{ccc4_7!20} 0.274(20)\\
\hline
2.0 & \cellcolor{ccc5_0!20} 0.052(7) & \cellcolor{ccc5_1!20} 0.115(10) & \cellcolor{ccc5_2!20} 0.158(12) & \cellcolor{ccc5_3!20} 0.169(13) & \cellcolor{ccc5_4!20} 0.206(15) & \cellcolor{ccc5_5!20} 0.222(16) & \cellcolor{ccc5_6!20} 0.264(19) & \cellcolor{ccc5_7!20} 0.266(19)\\
\hline
2.4 & \cellcolor{ccc6_0!20} 0.052(7) & \cellcolor{ccc6_1!20} 0.115(10) & \cellcolor{ccc6_2!20} 0.147(12) & \cellcolor{ccc6_3!20} 0.153(12) & \cellcolor{ccc6_4!20} 0.186(14) & \cellcolor{ccc6_5!20} 0.221(16) & \cellcolor{ccc6_6!20} 0.244(18) & \cellcolor{ccc6_7!20} 0.273(20)\\
\hline
2.8 & \cellcolor{ccc7_0!20} 0.052(7) & \cellcolor{ccc7_1!20} 0.116(10) & \cellcolor{ccc7_2!20} 0.147(12) & \cellcolor{ccc7_3!20} 0.162(13) & \cellcolor{ccc7_4!20} 0.187(14) & \cellcolor{ccc7_5!20} 0.224(16) & \cellcolor{ccc7_6!20} 0.242(18) & \cellcolor{ccc7_7!20} 0.264(19)\\
\hline
3.5 & \cellcolor{ccc8_0!20} 0.053(7) & \cellcolor{ccc8_1!20} 0.117(10) & \cellcolor{ccc8_2!20} 0.141(12) & \cellcolor{ccc8_3!20} 0.162(13) & \cellcolor{ccc8_4!20} 0.188(14) & \cellcolor{ccc8_5!20} 0.235(17) & \cellcolor{ccc8_6!20} 0.236(17) & \cellcolor{ccc8_7!20} 0.260(19)\\
\hline
5.0 & \cellcolor{ccc9_0!20} 0.053(7) & \cellcolor{ccc9_1!20} 0.117(10) & \cellcolor{ccc9_2!20} 0.141(12) & \cellcolor{ccc9_3!20} 0.162(13) & \cellcolor{ccc9_4!20} 0.189(14) & \cellcolor{ccc9_5!20} 0.244(17) & \cellcolor{ccc9_6!20} 0.236(17) & \cellcolor{ccc9_7!20} 0.259(19)\\
\hline
7.0 & \cellcolor{ccc10_0!20} 0.053(7) & \cellcolor{ccc10_1!20} 0.120(10) & \cellcolor{ccc10_2!20} 0.141(12) & \cellcolor{ccc10_3!20} 0.162(13) & \cellcolor{ccc10_4!20} 0.189(14) & \cellcolor{ccc10_5!20} 0.245(17) & \cellcolor{ccc10_6!20} 0.236(17) & \cellcolor{ccc10_7!20} 0.259(19)\\
\hline
9.5 & \cellcolor{ccc11_0!20} 0.053(7) & \cellcolor{ccc11_1!20} 0.120(10) & \cellcolor{ccc11_2!20} 0.141(12) & \cellcolor{ccc11_3!20} 0.162(13) & \cellcolor{ccc11_4!20} 0.189(14) & \cellcolor{ccc11_5!20} 0.245(17) & \cellcolor{ccc11_6!20} 0.236(17) & \cellcolor{ccc11_7!20} 0.259(19)\\
\hline
\end{tabular}
}
\end{minipage}
\noindent
\hspace{-2.3cm}
\begin{minipage}[c]{0.3\textwidth}\hspace{0.1cm}$\theta=0.72$\\~\\
\resizebox{0.4\textwidth}{!}{
%
%
\begin{tikzpicture}

\begin{axis}[%
width=3.566in,
height=3.566in,
at={(1.236in,0.481in)},
scale only axis,
xmin=-1,
xmax=1,
xtick={\empty},
ymin=-1,
ymax=1,
ytick={\empty},
axis background/.style={fill=white},
axis x line*=bottom,
axis y line*=left,
legend style={legend cell align=left, align=left, draw=white!15!black}
]

\addplot[area legend, draw=black, fill=black!10]
table[row sep=crcr] {%
x	y\\
-1	-1\\
-1	-0.5\\
-0.5	-0.5\\
-0.5	-1\\
}--cycle;

\addplot[area legend, draw=black, fill=white]
table[row sep=crcr] {%
x	y\\
-1	-0.5\\
-1	0\\
-0.5	0\\
-0.5	-0.5\\
}--cycle;

\addplot[area legend, draw=black, fill=black!10]
table[row sep=crcr] {%
x	y\\
-1	0\\
-1	0.5\\
-0.5	0.5\\
-0.5	0\\
}--cycle;

\addplot[area legend, draw=black, fill=white]
table[row sep=crcr] {%
x	y\\
-1	0.5\\
-1	1\\
-0.5	1\\
-0.5	0.5\\
}--cycle;

\addplot[area legend, draw=black, fill=white]
table[row sep=crcr] {%
x	y\\
-0.5	-1\\
-0.5	-0.5\\
0	-0.5\\
0	-1\\
}--cycle;

\addplot[area legend, draw=black, fill=black!10]
table[row sep=crcr] {%
x	y\\
-0.5	-0.5\\
-0.5	0\\
0	0\\
0	-0.5\\
}--cycle;

\addplot[area legend, draw=black, fill=white]
table[row sep=crcr] {%
x	y\\
-0.5	0\\
-0.5	0.5\\
0	0.5\\
0	0\\
}--cycle;

\addplot[area legend, draw=black, fill=black!10]
table[row sep=crcr] {%
x	y\\
-0.5	0.5\\
-0.5	1\\
0	1\\
0	0.5\\
}--cycle;

\addplot[area legend, draw=black, fill=black!10]
table[row sep=crcr] {%
x	y\\
0	-1\\
0	-0.5\\
0.5	-0.5\\
0.5	-1\\
}--cycle;

\addplot[area legend, draw=black, fill=white]
table[row sep=crcr] {%
x	y\\
0	-0.5\\
0	0\\
0.5	0\\
0.5	-0.5\\
}--cycle;

\addplot[area legend, draw=black, fill=black!10]
table[row sep=crcr] {%
x	y\\
0	0\\
0	0.5\\
0.5	0.5\\
0.5	0\\
}--cycle;

\addplot[area legend, draw=black, fill=white]
table[row sep=crcr] {%
x	y\\
0	0.5\\
0	1\\
0.5	1\\
0.5	0.5\\
}--cycle;

\addplot[area legend, draw=black, fill=white]
table[row sep=crcr] {%
x	y\\
0.5	-1\\
0.5	-0.5\\
1	-0.5\\
1	-1\\
}--cycle;

\addplot[area legend, draw=black, fill=black!10]
table[row sep=crcr] {%
x	y\\
0.5	-0.5\\
0.5	0\\
1	0\\
1	-0.5\\
}--cycle;

\addplot[area legend, draw=black, fill=white]
table[row sep=crcr] {%
x	y\\
0.5	0\\
0.5	0.5\\
1	0.5\\
1	0\\
}--cycle;

\addplot[area legend, draw=black, fill=black!10]
table[row sep=crcr] {%
x	y\\
0.5	0.5\\
0.5	1\\
1	1\\
1	0.5\\
}--cycle;

\end{axis}
\end{tikzpicture}%
}
\end{minipage}

\caption{Computed values of the approximate convergence factor\textcolor{ultramarine}{$\rho$ and corresponding preconditioned CG iteration counts (between parenthesis)} w.r.t.\ parameters $\varepsilon$ on rows and mesh size $h$ on columns. In each table \textcolor{ultramarine}{the} pattern of diffusion coefficient $\mu$ and the strong threshold $\theta$ (shown on the left) is fixed. The background color depends on $\rho$ with colormap \protect\includegraphics[height=0.2cm]{rainbow.png}. It is scaled to range between \textcolor{reviewV2}{the} minimum and maximum value (among all the tables) of $\rho$. Checkerboard pattern for $\mu$ as in Figure~\ref{fig:mu_pattern}(d).}
\label{tab:rho2}
\end{table}

\subsection{Relation between $\theta$ and computational costs}\label{sec:machine_time}
\textcolor{reviewV2}{We also investigate the relation between $\theta$ and the} CPU time $t$ needed to solve the linear system. Indeed, this is the quantity that we want to minimize \textcolor{reviewV2}{in practice.} 

\begin{figure}[htbp]
\begin{center}
	\includegraphics[width=1.0\textwidth]{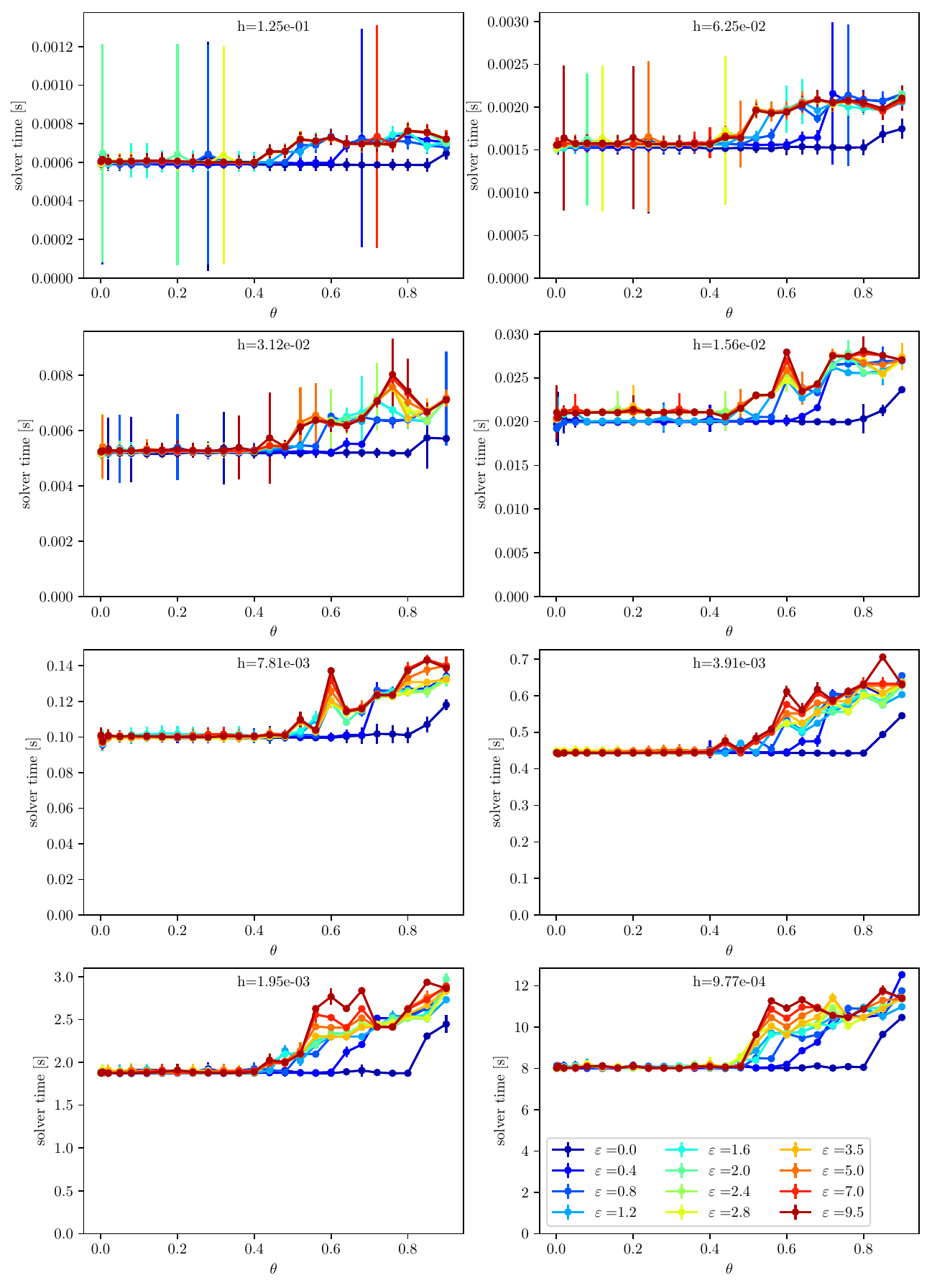}
\end{center}
\caption{Mean and standard deviation \textcolor{ultramarine}{(visualized as errorbar)} of the elapsed CPU time \textcolor{ultramarine}{$t$} to solve the linear system of equations (in seconds) based on employing the AMG preconditioned CG. In each plot we have fixed a different mesh size \textcolor{reviewV2}{$h=1.25e$-$1, ..., 9.773$-$4$}. Each line represents \textcolor{ultramarine}{the solver CPU time $t$ for a fixed} \textcolor{reviewV2}{choice of $\varepsilon$ entering in the definition} of the diffusion coefficient $\mu$. The pattern of $\mu$ is reported in Figure~\ref{fig:mu_pattern}(c).}
\label{fig:theta_vs_time_m1r2}
\end{figure}

\begin{figure}[htbp]
\begin{center}
	\includegraphics[width=1.0\textwidth]{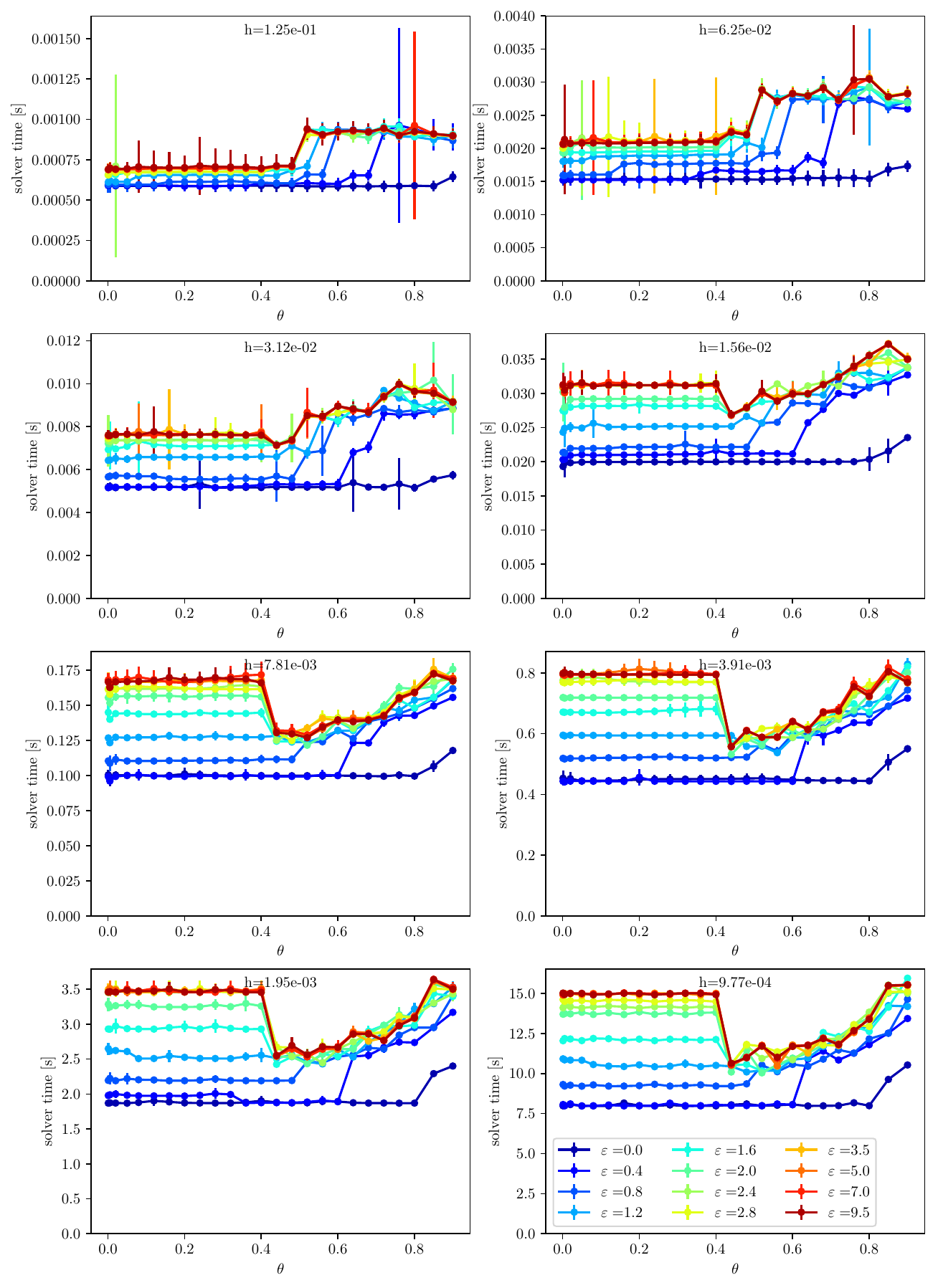}
\end{center}

\caption{Mean and standard deviation \textcolor{ultramarine}{(visualized as errorbar)} of the elapsed CPU time \textcolor{ultramarine}{$t$} to solve the linear system of equations (in seconds) based on employing the AMG preconditioned CG. In each plot we have fixed a different mesh size \textcolor{reviewV2}{$h=1.25e$-$1, ..., 9.773$-$4$}. Each line represents \textcolor{ultramarine}{the solver CPU time $t$ for a fixed} \textcolor{reviewV2}{choice of $\varepsilon$ entering in the definition} of the diffusion coefficient $\mu$. The pattern of $\mu$ is reported in Figure~\ref{fig:mu_pattern}(d).}
\label{fig:theta_vs_time_m2r2}
\end{figure}

\textcolor{reviewV2}{In order to have an accurate estimate of the CPU time $t$ we gather multiple samples by repeating each simulation. The number of iteration we choose was the minimum number such that the standard deviation of $t$ did not change significantly when increasing the number of samples. Namely, we employed $200$, $100$, $50$, $20$, $10$, $7$, $5$, and $4$ iterations for each mesh refinement from the coarsest to the finest, respectively.} We report the plots of $t$ vs $\rho$ in Figures~\ref{fig:theta_vs_time_m1r2} and~\ref{fig:theta_vs_time_m2r2}. Notice that despite the large number of repetitions the standard deviation (shown as an errorbar) in some cases is still large. \textcolor{ultramarine}{These results are also useful to analyze the \textcolor{reviewV2}{relation} between $\theta$ and $\rho$ since we will show that the least square analysis \textcolor{reviewV2}{seems to indicate} that they have a linear \textcolor{reviewV2}{relation}.}
We also observe that for small values of the strong threshold parameter ($\theta \leq 0.3$) there is an interval where the CPU time is almost constant: this appears to be true for all the test cases addressed. Since a smaller strong threshold parameter means that more connections are discarded in the coarsening phase, one would expect that as $\theta$ gets smaller, then the approximate convergence factor $\rho$ may deteriorate, which in turn would lead to larger CPU times. A possible motivation of this behaviour is to consider that among the settings of BoomerAMG there is parameter that prevents the coarsening from being too small (in the present test its value has been set as default, i.e. equal to $1$). Thus, the coarse system associated to $\mathrm{A}_H$ is still effectively damping the smooth components. Indeed, from Figure~\ref{fig:theta_vs_levels} we can see how the number of levels and the size of the coarse system is constant for small $\theta$.

\subsubsection{Choice of the performance index $p$}

We are now interested in finding a scalar $p$ that evaluates how good the AMG configuration is. Two possible choices for such performance index are the approximate convergence factor $\rho$, which measures how rapidly the linear solver converges, and the elapsed CPU time $t$. 

\begin{figure}[t]
    \makebox[\textwidth][c]{\includegraphics[trim={0 1cm 0 3cm},clip,width=10cm]{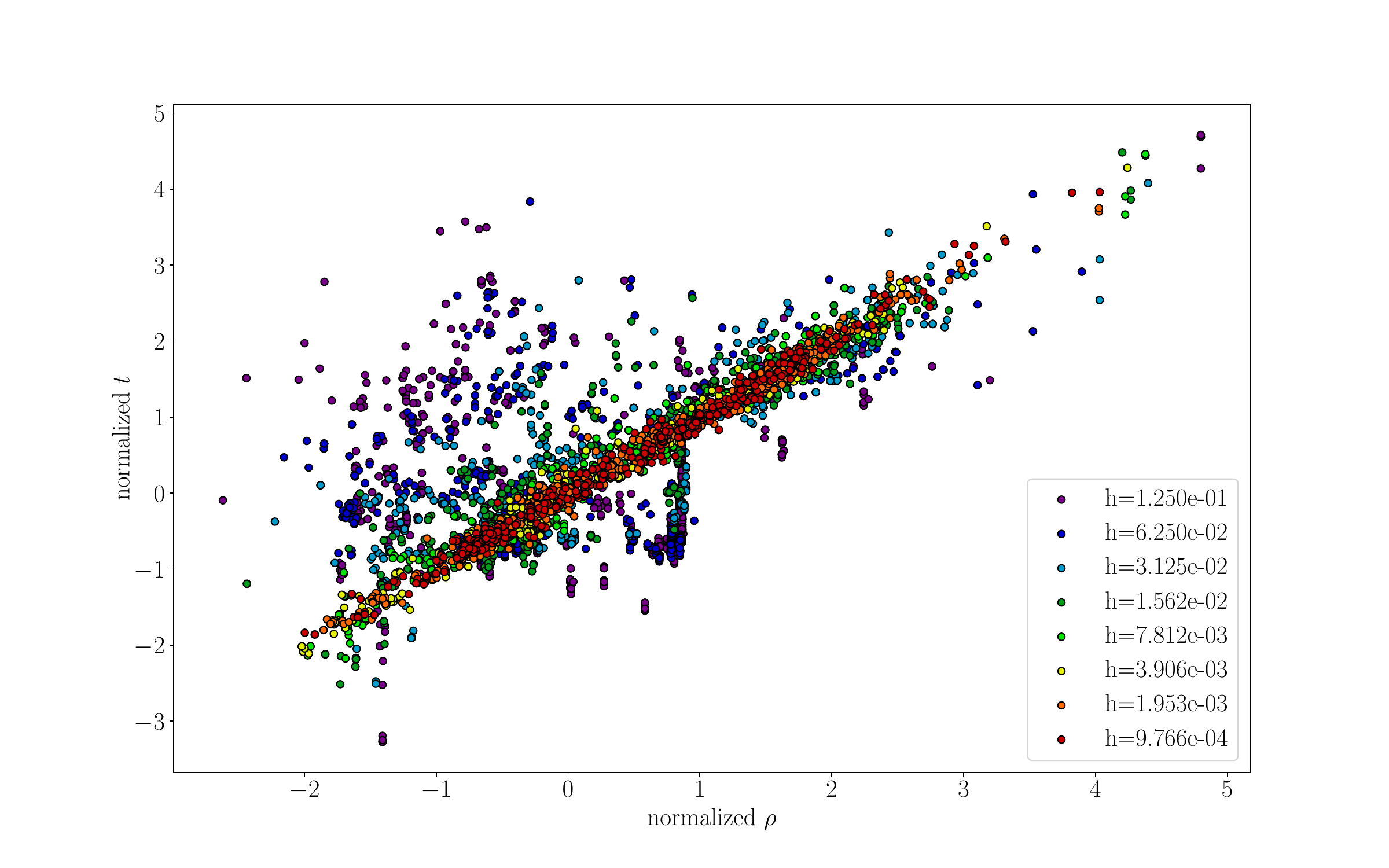}}
    \caption{Scatter plot of the (average) elapsed CPU time ($t$) \textcolor{reviewV2}{versus the} approximate convergence factor $\rho$. \textcolor{reviewV2}{Different colors identify different mesh sizes $h$.} \textcolor{reviewV2}{The data are} normalized (in both components) with respect to the corresponding data in the same test case.}
    \label{fig:normal_scatter}
\end{figure}

We now proceed to analyze the relation between the elapsed CPU time $t$ and the approximate convergence factor $\rho$. In Figure~\ref{fig:normal_scatter}, we show a scattered plot of the elapsed CPU time ($t$) as a function of $\rho$, for different values of the mesh size $h$. The results are normalized with respect to the data that belong to the same test case. A linear \textcolor{reviewV2}{relation} between $t$ and $\rho$ \textcolor{reviewV2}{can be clearly identified}. This is also confirmed by the results shown in Table~\ref{tab:ols_rho_t}, where we report the \textcolor{ultramarine}{least square} analysis of the data of Figure~\ref{fig:normal_scatter}. We highlight that these results support the hypothesis of a \textcolor{reviewV2}{relation} between $\rho$ and $t$. We can explain the poor correlation for coarse mesh sizes $h$ due to the higher relative uncertainty of the measure. Indeed, as $h$ gets smaller the coefficient of determination $\text{R}^2$ improves.

\begin{table}[t]
\small
\makebox[\textwidth][c]{
     \begin{tabular}{ l r r r r r r r r }
      \hline
      $h$ & 1.25e-1 & 6.25e-2 & 3.12e-2 & 1.56e-2 & 7.81e-3 & 3.91e-3 & 1.95e-3 & 9.77e-4\\
      \hline
      \textcolor{reviewV2}{$nn$}             & 1200 & 1200 & 1200 & 1200  & 1200 & 1200 & 1200 & 1200\\
      R${}^2$          & 0.155 & 0.564 & 0.762 & 0.912 & 0.985 & 0.992 & 0.991 & 0.993\\
      F-statistic      & 220.5 & 1550 & 3836 & 1.24e4  & 7.88e4 & 1.58e5 & 1.34e5 & 1.68e5\\
      AIC              & $-$2.15e4 & $-$1.79e4 & $-$1.53e4 & $-$1.30e+4 & $-$1.09e4 & $-$7864 & $-$4048 & $-$800.1\\
      \textcolor{reviewV2}{$\hat \rho_1$}    & 3.49e-4 & 3.05e-3 & 1.25e-2 & 4.99e-2 & 0.263 & 1.245 & 5.259 & 20.95\\
      \textcolor{reviewV2}{SE($\hat \rho_1$)}        & 2.35e-5 & 7.76e-5 & 2.02e-4 & 4.48e-4   & 9.35e-4 & 3.13e-3 & 1.42e-2 & 5.11e-2\\
      t-value \textcolor{reviewV2}{$\hat \rho_1$}   & 14.848 & 39.369 & 61.937 & 111.299  & 280.677 & 397.857 & 366.596 & 409.565\\
      p-value \textcolor{reviewV2}{$\hat \rho_1$}   & $<0.001$ & $<0.001$ & $<0.001$ & $<0.001$ & $<0.001$ & $<0.001$ & $<0.001$ & $<0.001$\\
      \hline
     \end{tabular}
     }
\caption{\textcolor{reviewV2}{Linear least square analysis of the model $t = \hat \rho_1 \rho + \hat \rho_0$. Data are grouped by the mesh size $h$.} $t$ is the elapsed CPU time and $\rho$ is the approximate converge factor. \textcolor{ultramarine}{For the definition of the quantities appearing in the first column we refer to Table~\ref{tab:ls_def}}.}
\label{tab:ols_rho_t}
\end{table}

\textcolor{reviewV2}{Therefore, in the following we use as performance index the convergence factor $\rho$. I}ndeed, $\rho$ is not machine nor implementation dependent, thus leading to reproducible results.

The dataset that we build contains numerical simulations made with every combination of parameters among $8$ mesh sizes $h$, $25$ values of $\theta$, $12$ values of $\varepsilon$ and $4$ patterns of $\mu$ (Figure~\ref{fig:mu_pattern}) for a total of $9600$ samples.

\section{ANN-enhanced AMG Method} \label{sec:ann-enhancedAMG}

\textcolor{ultramarine}{In this section, we design and use ANNs to predict the value of the strong threshold parameter $\theta^*$ that maximizes the performance of the AMG method\textcolor{reviewV2}{, measured in terms of corresponding convergence factor $\rho$}.} Our goal is to design a model, namely an ANN, that enables predictions of the optimal strong threshold $\theta$ for a given model problem. We remark that, in the framework discussed in Section~\ref{sec:model_prob}, fixing a test case (model problem) is equivalent to fixing the matrix {$\mathrm{A}_h$ defined in Eq.~\bref{eq:2linear_sys}.} We define the optimal value of strong threshold parameter $\theta^*$ for a certain test case as the minimizers of \textcolor{ultramarine}{the convergence factor $\rho$ = $\rho(\mathrm A_h, \theta)$}.

Then, we build our model (the ANN) to predict the convergence factor $\rho$ of the AMG in a fixed test case and with a fixed strong threshold parameter $\theta$\textcolor{reviewV2}{. More precisely, $\mathscr F$} is the ANN such that
\begin{equation}
    \mathscr F(\texttt{normalize}(\texttt{pooling}(\mathrm{A}_h, m\textcolor{reviewV2}{, \texttt{op}})), - \log_2(h), \theta; \boldsymbol{\gamma}) = \rho.
\end{equation}
\textcolor{reviewV2}{Here, the} \textcolor{ultramarine}{\texttt{pooling} will be introduced in Section~\ref{sec:pooling}, \textcolor{reviewV2}{whereas} $m$, \textcolor{reviewV2}{\texttt{normalize} and \texttt{op} are \textcolor{reviewV2}{so-called} hyperparameters of the model} (\textcolor{reviewV2}{that will be introduced and discussed} in detail in Section~\ref{sec:normalization})}\textcolor{reviewV2}{. Finally,} $\boldsymbol{\gamma}$ are the parameters that define the ANN \textcolor{reviewV2}{(see Eq.~\bref{eq:layer})}. There are two reasons to adopt this approach: first, it is possible to quantify the improvement on the performance that we expect; second, each numerical simulation can be added to the dataset making this process less computationally expensive and more flexible.

If otherwise not stated, we will use a $60\%$-$20\%$-$20\%$ split of the dataset into training-validation-test.

\subsection{ANN-based prediction of the optimal strong threshold parameter $\theta$}
In this section, we discuss how to predict the optimal strong threshold parameter $\theta$ to be used in the ANN-enhanced AMG Algorithm~\ref{a:ann_amg} without using any prior assumption on the diffusion coefficient $\mu$. In other words, we do not rely on the fact that $\mu$ shows a finite number of patterns.

The variables that we use as inputs of the ANN $\mathscr F$ are $\theta$, $-\log_2(h)$ and a set of variables $\hat{\mathrm V} = \texttt{normalize}(\texttt{pooling}(\mathrm{A}_h, m \textcolor{reviewV2}{, \texttt{op}}))$ that is extracted from the matrix $\mathrm{A}_h$ of the linear system \bref{eq:2linear_sys} by means of the pooling Algorithm~\ref{a:view} and a normalization algorithm (extraction step). This extraction process should be computationally cheap, indeed this approach is worthwhile only until the process of predicting the optimal value of $\theta$ has a negligible computational cost with respect to the elapsed CPU time to solve the linear system. \textcolor{ultramarine}{We use $-\log_2(h)$ instead of $h$ as input of the ANN since, by construction, $h$ is not linearly distributed. It is a well known that ANNs can be more easily trained if the inputs have the same order of magnitude.}

We recall at this stage that the optimal parameter $\theta^*$ to be used in the linear solver with AMG preconditioner (step \textcolor{reviewV2}{4} of Algorithm~\ref{a:ann_amg}) is such that
$$
\theta^* = \argmin_{\, \theta \in (0,1]} \, \, \mathscr F(\hat{\mathrm{V}}, -\log_2(h), \theta;\boldsymbol{\gamma}).
$$
\textcolor{ultramarine}{In practice, $\theta^*$ is found by first evaluating
$$
\theta^{(0)} = \argmin_{\, k \in \{20,30,...,900\}} \, \, \mathscr F(\hat{\mathrm{V}}, -\log_2(h), 0.001k;\boldsymbol{\gamma}),
$$
and then applying a suitable number of steps of the gradient descent algorithm
$$
\theta^{(k)} = \theta^{(k-1)} - \alpha \nabla_\theta \mathscr F(\hat{\mathrm{V}}, -\log_2(h), \theta^{(k-1)};\boldsymbol{\gamma}), \quad k \geq 0,
$$
where $\alpha = 10^{-5}$ is the learning rate. The gradient $\nabla_\theta \mathscr F$ can be computed by the automatic differentiation algorithm of Tensorflow.
However, we empirically found that this second step appears to be unnecessary since it gives small to negligible improvements.
}

\subsubsection{Pooling (step 1 of Algorithm~\ref{a:ann_amg})}\label{sec:pooling}
We introduce what we call the \textit{view} $\mathrm V$ of the matrix A${}_h$. 

\textcolor{reviewV2}{First, let us define the following hyperparameters. Let $m \in \mathbb N$ be a positive integer that describes the size of V$\in \mathbb R^{m\times m}$. It must be large enough so that the features of $\mathrm{A}_h$ are not lost. At the same time, \textcolor{reviewV2}{$m$ should} not be too large to avoid expensive computations in the forward propagation step. Let \texttt{op}$: \mathbb R \times \mathbb R \rightarrow \mathbb R$ be a function that combines two values. In the field of computer vision this function is usually the sum \texttt{op=sum} where \texttt{sum}$(v_1, v_2) = v_1 + v_2$, or the maximum of two numbers.}

\textcolor{reviewV2}{We define the view $\mathrm V \in \mathbb R^{m\times m}$ of the matrix A${}_h$ and the non-zeros count $\mathrm C \in \mathbb N^{m\times m}$ (i.e.\ the matrix where each entry $(\mathrm C)_{ij}$ is the number of non-zero elements of $\mathrm{A}_h$ used to compute $(\mathrm V)_{ij}$) as (V, C) = \texttt{pooling}$(\mathrm A_h, m, \texttt{op})$, where \texttt{pooling} is defined in Algoritm~\ref{a:view}}. \textcolor{ultramarine}{The insight in this algorithm \textcolor{reviewV2}{stems} from the operator used in the pooling layer of CNNs. On one hand we downscale the input, significantly reducing the computational cost, and on the other hand we also gain translation invariance. Moreover, we also prune details that may not be useful for the task.}

\begin{algorithm}[t]
    access $\mathrm{A}_h$ in COO form and extract its size: \texttt{val}, \texttt{row}, \texttt{col}, $n$ $\leftarrow$ $\mathrm{A}_h$\;
    initialize $\mathrm{V}$ to an $m \times m$ dense matrix with all zero entries\;
    initialize $\mathrm{C}$ to an $m \times m$ dense matrix with all zero entries\;
    $q$ $\leftarrow$ $n / m$\;
    $p$ $\leftarrow$ $n$ \textnormal{mod} $m$\;
    $t$ $\leftarrow$ $(q+1)p$\;
    \For{$k = 0$ \KwTo $\textnormal{\texttt{val.size()}}-1$} {
        $i\,$ $\leftarrow$ \texttt{row}[$\,k\,$]$/(q+1)$ \textnormal{\textbf{if}} (\texttt{row}[$\,k\,$] $< t$) \textnormal{\textbf{else}} $($\texttt{row}[$\,k\,$] $- \,t)/q + p$\;
        $j$ $\leftarrow$ \texttt{col}[$\,k\,$]$/(q+1)$ \textnormal{\textbf{if}} (\texttt{col}[$\,k\,$] $< t$) \textnormal{\textbf{else}} $($\texttt{col}[$\,k\,$] $- \,t)/q + p$\;
        $\mathrm{V}_{ij}$ $\leftarrow$ \textcolor{reviewV2}{\texttt{op}($\mathrm{V}_{ij}, \texttt{val}[\,k\,]$)}\;
        $\mathrm{C}_{ij}$ $\leftarrow$ $\mathrm{C}_{ij} + 1$\;
    }
    \KwRet{$\mathrm V, \mathrm C$}\;

    \caption{Pooling algorithm \newline $(\mathrm V, \mathrm C)=\texttt{pooling}(\mathrm{A}_h,m\textcolor{reviewV2}{,\texttt{op}})$}
    \label{a:view}
\end{algorithm}

\textcolor{ultramarine}{By exploiting the storage data structure of sparse matrices (for instance compressed row storage or coordinate lists) to access the elements of $\mathrm{A}_h$, we realize the pooling with complexity $O(nnz)$, where $nnz$ is the number of non-zero elements in the matrix $\mathrm{A}_h$. \textcolor{reviewV2}{In Algorithm~\ref{a:view},  for the sake of simplicity,} it is assumed that the matrix $\mathrm{A}_h$ is stored in coordinate lists format.} 

We have also measured the elapsed CPU time by \textcolor{reviewV2}{the \texttt{pooling} algorithm. In each simulation it seems that Algorithm~\ref{a:view} requires a negligible CPU time compared to the global one. We notice that this algorithm could easily be extended to work in parallel.} \textcolor{ultramarine}{\textcolor{reviewV2}{We also point out that Algorithm~\ref{a:view}} does not rely on the connectivity of the mesh nor on the definition of the coefficients, thus it should also work in more complex cases, as for example in the case of unstructured meshes. }

\subsubsection{Normalization (step 2 in Algorithm~\ref{a:ann_amg})}\label{sec:normalization}

\begin{table}[t]
\begin{minipage}[b]{0.47\linewidth}
\centering
{\color{reviewV2}\begin{tabular}{ lll  }
  \hline
  Name & \texttt{normalize}  \\ 
  \hline
  \texttt{std+id}    & $\hat{\mathrm V} = $\texttt{std}  $(\mathrm V)                         $ \\
  \texttt{std+avg}   & $\hat{\mathrm V} = $\texttt{std}  $(\texttt{avg}(\mathrm V, \mathrm C))$ \\
  \texttt{scale+id}  & $\hat{\mathrm V} = $\texttt{scale}$(\mathrm V)                         $ \\
  \texttt{scale+avg} & $\hat{\mathrm V} = $\texttt{scale}$(\texttt{avg}(\mathrm V, \mathrm C))$ \\
  \texttt{log+id}    & $\hat{\mathrm V} = $\texttt{log}  $(\mathrm V)                         $ \\
  \texttt{log+avg}   & $\hat{\mathrm V} = $\texttt{log}  $(\texttt{avg}(\mathrm V, \mathrm C))$ \\
  \hline
 \end{tabular}
 }
\caption{\textcolor{reviewV2}{Summary of the normalization techniques for the view $\mathrm V$. See Eq.~\bref{eq:norm_modes} and Eq.~\bref{eq:mean46} for further details.}}
\label{tab:norm_modes}
\end{minipage}\hfill
\begin{minipage}[b]{0.47\linewidth}
	\centering
    \includegraphics[width=0.8\textwidth]{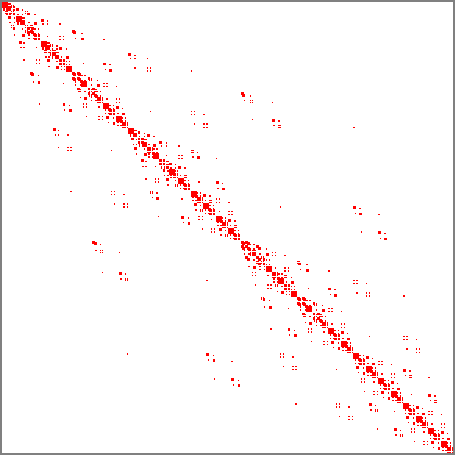}\\
    \captionof{figure}{\textcolor{ultramarine}{Sparsity pattern of the matrix $\mathrm{A}_h$ for $h=1.25$e-$1$ of the case 1 (Section~\ref{sec:case_1}).}}
    \label{fig:sparsity_pattern}
    
\end{minipage}
\end{table}

\begin{figure}[t]
\begin{center}
    \makebox[\textwidth][c]{\includegraphics[trim={0 35mm 0 35mm},clip,width=1.3\textwidth]{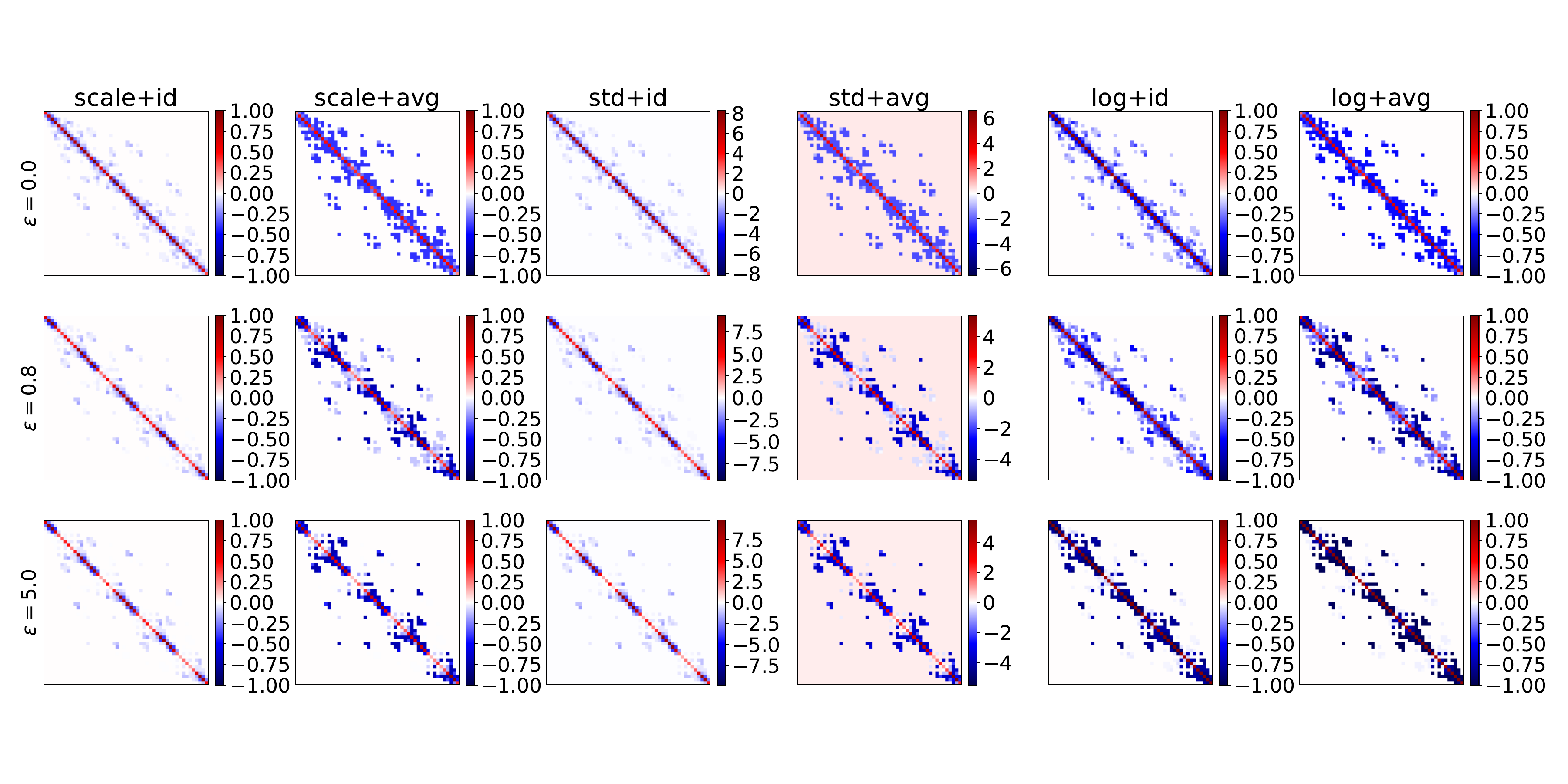}}
    \caption{\textcolor{ultramarine}{Graphical representation of \textcolor{reviewV2}{six} different normalization of the view matrix $\hat{\mathrm{V}} = \texttt{normalize}(\texttt{pooling}(\mathrm{A}_h, 50\textcolor{reviewV2}{, \texttt{sum}})) \in \mathbb{R}^{50\times 50}$ \textcolor{reviewV2}{for m}esh dimension $h=1.25$e-$1$, $\mu$ pattern Figure~\ref{fig:mu_pattern}(d) \textcolor{reviewV2}{and} $\varepsilon=0,0.8,5$ for the first, second and third rows, respectively.}}
        \label{fig:poolingpattern}

\end{center}
\end{figure}

We observe that the \textit{view} \textcolor{reviewV2}{V} defined in the previous section cannot be used as input of an ANN yet. In particular, it features very large values that might impact the stability of the gradient algorithm. For this reason, propose three normalization techniques:
\begin{equation}
\label{eq:norm_modes}
\begin{split}
  &\textcolor{reviewV2}{\texttt{std}(\mathrm{V})_{ij}} = \frac{(\mathrm{V})_{ij}-\bar{\mathrm{v}}}{\sigma}, \ \bar{\mathrm{v}} = \sum_{i,j}\frac{(\mathrm{V})_{ij}}{ m^2}, \ \sigma = \sqrt{\frac{1}{m^2}\sum_{i,j} [(\mathrm{V})_{ij} - \bar{\mathrm{v}}]^2},\\
  &\textcolor{reviewV2}{\texttt{scale}(\mathrm{V})_{ij}} = (\mathrm{V})_{ij} / \max_{i,j} \abs{(\mathrm{V})_{ij}},\\
  &\textcolor{reviewV2}{\texttt{log}(\mathrm{V})_{ij} = \texttt{scale}(\log(\abs{(V)_{ij}} + 1) (V)_{ij}/\abs{(V)_{ij}})}.
\end{split}
\end{equation}
\textcolor{reviewV2}{The first approach is the most employed in the field of deep learning. The argument behind the second and third definition of Eq.~\bref{eq:norm_modes}} is that we would like to preserve the sparsity pattern of the matrix. \textcolor{reviewV2}{In particular, the \texttt{log} normalization yields linearly distributed values, since the exponent $\varepsilon$ of diffusion coefficient $\mu = 10^\varepsilon$ is linearly distributed. Another possibility is to apply} these normalization to  the element-wise division of $\mathrm{V}$ and $\mathrm{C}$ 
\begin{equation}
    \textcolor{reviewV2}{\texttt{avg}(\mathrm{V}, \mathrm{C})_{ij}} = \frac{(\mathrm{V})_{ij}}{ (\mathrm{C})_{ij}},
\label{eq:mean46}
\end{equation}
\textcolor{ultramarine}{with the exception that if $(\mathrm C)_{ij} = 0$ then \textcolor{reviewV2}{$\texttt{avg}(\mathrm{V}, \mathrm{C})_{ij} = 0$}. Hence, the sparsity pattern of $\mathrm V$ is preserved, indeed $(\mathrm C)_{ij} = 0$ only if $(\mathrm V)_{ij} = 0$.}
\textcolor{reviewV2}{Table~\ref{tab:norm_modes} summarizes all the normalization techniques we propose. Figure~\ref{fig:sparsity_pattern} shows and example of the sparsity pattern of the matrix $\mathrm A_h$ ($h=1.25e$-$1$) and} \textcolor{ultramarine}{ Figure~\ref{fig:poolingpattern} shows three \textcolor{reviewV2}{examples} of a side-to-side view of the four normalizations $\hat{\mathrm{V}}$ of the view matrix $\mathrm{V}$.}
At this stage, we have obtained a matrix $\hat{\mathrm V} \in \mathbb R^{m \times m}$, \textcolor{reviewV2}{$\hat{\mathrm V} = \texttt{normalize}(\mathrm V, \mathrm C)$} with $m$ chosen a priori. We will discuss our choice in the next section. 

\begin{figure}[t]
    \makebox[\textwidth][c]{\includegraphics[trim={0 0cm 0 0cm},clip,width=13cm]{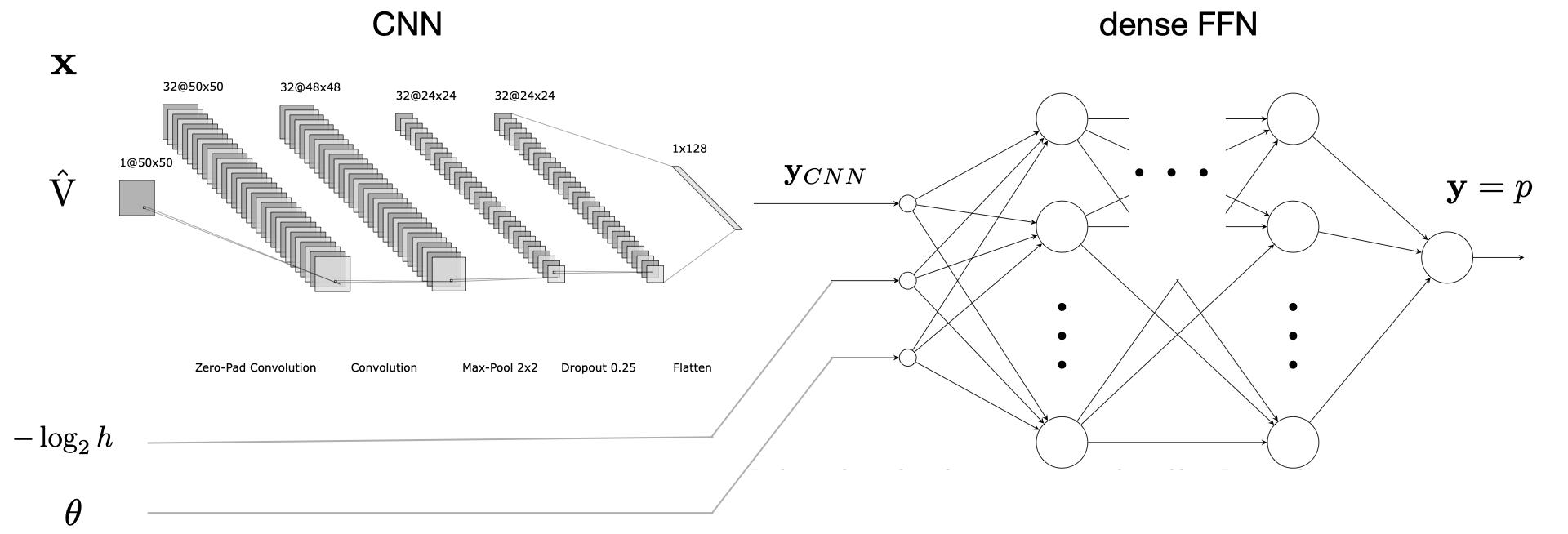}}
    \caption{Architecture \textcolor{ultramarine}{of our model represented by the ANN $\mathscr F(\hat{\mathrm{V}}, - \log_2(h), \theta; \boldsymbol{\gamma}) = p = \rho$ where $\hat{\mathrm{V}} = \texttt{normalize}(\texttt{pooling}(\mathrm{A}_h, 50\textcolor{reviewV2}{, \texttt{sum}}))$. In particular, it is comprised by the composition of two ANNs:} a CNN such that $\mathbf y_{CNN} =\mathscr F_{CNN}(\hat{\mathrm{V}};\boldsymbol{\gamma}_{CNN})$ and a dense FFN $\mathscr F_{FFN}(\mathbf y_{CNN},-\log_2(h), \theta ;\boldsymbol{\gamma}_{FNN}) = \rho$. \textcolor{ultramarine}{Given \textcolor{reviewV2}{$\mathrm{\hat{V}}$}, the mesh dimension $h$ and $\theta$ it predicts $\rho$ and thus also the optimal $\theta^*$ that minimizes $\rho$.}}
    \label{fig:ann}
\end{figure}

\begin{algorithm}[t]
{
    $\mathrm{V}$, $\mathrm{C}$  $\leftarrow$ \texttt{pooling}($\mathrm{A}_h$, $m$ \textcolor{reviewV2}{, \texttt{op}}) (Algorithm~\ref{a:view})\;
    \textcolor{reviewV2}{$\hat{\mathrm{V}} \leftarrow \texttt{normalize}(\mathrm V, \mathrm C)$ (Table \ref{tab:norm_modes})\;} 
    $\theta^*$  $\leftarrow$ $\textnormal{argmin}_{\theta} \, \mathscr F(\hat{\mathrm{V}}, -\log_2(h), \theta;\boldsymbol{\gamma})$ \; 
    $\mathbf{u}_h^{(k+1)}  \leftarrow \texttt{AMG}(\mathbf{u}_h^{(0)}, \mathrm{A}_h, \mathbf{f}_h, \theta^*, \nu_1, \nu_2, N_{max}, tol)$ (Algorithm 2)
    }

    \caption{ANN-enhanced AMG \newline {$\mathbf{u}_h^{(k+1)} = \texttt{ANN\_AMG}(\mathbf{u}_h^{(0)}, \mathrm{A}_h, h, \mathbf{f}_h, \nu_1, \nu_2, N_{max}, tol, \gamma$)}}
    \label{a:ann_amg}
\end{algorithm}

\subsubsection{The ANN-enhanced AMG Algorithm}

We show in Algorithm~\ref{a:ann_amg} how we intend to use the prediction of the optimal strong threshold parameter $\theta^*$ realized by ANN within the AMG solver, which we call ANN-enhanced AMG algorithm. In particular, our approach determines $\theta^*$ to be used in the {\tt AMG} algorithm starting from the matrix $\mathrm A_h$ and the mesh size $h$. This leverages on a map from a manipulation of $\mathrm A_h$ ($\hat{\mathrm V}$), $h$ and $\theta$ to a suitable performance index $p$ of the AMG solver \textcolor{ultramarine}{($\rho$)}. Specifically, this map is realized by an ANN $\mathscr F(\mathbf x; \boldsymbol{\gamma})$ such that its inputs are $\mathbf x =(\hat{\mathrm{V}}, -\log_2(h), \theta)$, while the output $\mathbf y$ (the predicted value of the regression) coincides with a suitable performance index, say $\mathbf y=p(\mathrm A_h,h,\theta$), of the linear solver with AMG preconditioner, \textcolor{ultramarine}{which we select as the approximated convergence factor $p=\rho(\mathrm{A}_h, \theta)$}. The steps in the ANN-enhanced AMG Algorithm~\ref{a:ann_amg} are the following:
\begin{itemize}
    \item (1--\textcolor{reviewV2}{2}) as the matrix $\mathrm A_h$ can not be directly used as input of an ANN $\mathscr F$, suitable pooling and normalization steps are performed to assemble $\hat{\mathrm V}=\texttt{normalize}(\texttt{pooling}(\mathrm{A}_h, m, \textcolor{reviewV2}{\texttt{op}})$ from $\mathrm A_h$, where \texttt{pooling} is defined in Algorithm~\ref{a:view}, \textcolor{reviewV2}{and \texttt{op}, \texttt{normalize} and $m$ are hyperparameters of our model defined inside Section~\ref{sec:pooling} and Section~\ref{sec:normalization}, respectively};
    \item (\textcolor{reviewV2}{3}) the ANN $\mathscr F$ built for the model problem is used to determine $\theta^*$ in order to minimize \textcolor{ultramarine}{the approximate convergence factor $\rho$} of the AMG;
    \item (\textcolor{reviewV2}{4}) the {\tt AMG} Algorithm~\ref{a:mgcycle} is used with $\theta^*$.
\end{itemize}
    
\subsubsection{ANN architecture (building the ANN of step \textcolor{reviewV2}{3} in Algorithm~\ref{a:ann_amg})}\label{sec:aaarchit}
We now build the ANN $\mathscr F(\mathbf x; \boldsymbol{\gamma}) \textcolor{ultramarine}{\; = \mathscr F(\hat{\mathrm{V}}, -\log_2(h), \theta; \boldsymbol{\gamma})}$. We recall that, following the former pooling and normalization steps, we have 
$$\hat{\mathrm V}=\texttt{normalize}(\texttt{pooling}(\mathrm{A}_h, m, \textcolor{reviewV2}{\texttt{op}})),$$
\textcolor{ultramarine}{where} \texttt{pooling} is defined in Algorithm~\ref{a:view} and \textcolor{reviewV2}{\texttt{op},} \texttt{normalize} \textcolor{reviewV2}{and $m$ are the hyperparameters defined in Section~\ref{sec:pooling} and} Section~\ref{sec:normalization}\textcolor{reviewV2}{, respectively}.
We use a model that is the composition of two networks as depicted in Figure~\ref{fig:ann}. Since the matrix view $\hat{\mathrm V}$ is a structured input, we first employ a CNN such that 
$$
\mathbf y_{CNN} =\mathscr F_{CNN}(\hat{\mathrm{V}};\boldsymbol{\gamma}_{CNN}),
$$
depending on the parameters $\boldsymbol{\gamma}_{CNN}$. Then, the output of the first CNN altogether with the remaining inputs $-\log_2(h)$ and $\theta$, constitute the inputs of a second network, which we select as a dense Feed Forward Network (FFN). This dense FFN is such that 
$$\textcolor{reviewV2}{\rho(A_h,\theta)} \; =\mathscr F_{FFN}(\mathbf y_{CNN},-\log_2(h), \theta ;\boldsymbol{\gamma}_{FNN})
$$ 
and depends on the parameters $\boldsymbol{\gamma}_{CNN}$.

\subsubsection{Evaluating the performance of the model}\label{sec:model_eval_def}
Since the a priori choice of the strong threshold parameter $\theta$ is based on the map \textcolor{reviewV2}{$\mathrm A_h \rightarrow \theta^*$ defined by step (3) of Algorithm~\ref{a:ann_amg}, it is not enough to have a small loss to verify that the model is accurate.} \textcolor{ultramarine}{With this aim, we introduce some quantities of interest. \textcolor{reviewV2}{Let $\mathrm A_h$ be fixed, and let}
\begin{itemize}
    \item $\rho_{\textnormal{ANN}}$ be the convergence factor of the AMG-ANN algorithm
    \item $\rho_{0.25}$ be the convergence factor of the AMG method for $\theta=0.25$
    \item $\rho_{\textnormal{MIN}}$ be the convergence factor of the AMG method with 
\textcolor{reviewV2}{$$\theta^*=\argmin_{\theta \in \textnormal{dataset for this A}_h} \rho(\theta\textcolor{reviewV2}{; \mathrm A_h}).$$}
\end{itemize}
\vspace{-0.5cm}
Moreover, we define
\begin{equation}
    P = 1 - \frac{\rho_\textnormal{ANN}}{ \rho_{0.25}} \quad \textnormal{  and  } \quad P_\textnormal{MAX} = 1 - \frac{\rho_{MIN}}{ \rho_{0.25}}.
    \label{eq:perf}
\end{equation}
Finally, we define $PB$ as the percentage of cases where $P \geq 0$, and we define $P_{<0}$ as the performance $P$ of the cases where $P < 0$.}

\subsection{Test Case 1}\label{sec:case1}

\textcolor{reviewV2}{We test our algorithm fixing the hyperparameters of the view, we discuss their tuning in Section~\ref{sec:view_hyper}. Namely, we employ \texttt{op=sum}, $m = 50$} (our choice is motivated by interpreting $\hat{\mathrm V}$ as a color image in input to the first CNN network; experience indicates that this kind of CNN network is able to excellently process color images of similar size) \textcolor{reviewV2}{ and \texttt{normalize=std+id} (see Table~\ref{tab:norm_modes}).}

Table~\ref{tab:general_arch} shows the results of changing the architecture of the model. In particular, we consider architectures with two convolutional layers each composed by a convolution with zero-padding, $3\times 3$ kernel and ReLU activation and $D_i-1$ other convolutions with $3\times 3$ kernel and ReLU activation (without padding). The last elements of the convolutional layers are a $2 \times 2$ max-pooling and Dropout with rate $P_i$, each layer has $W_i$ hidden units (for $i=1,2$). The output of the convolutional part has $O$ hidden units; the dense part is composed by $D_3$ dense layers with $W_3$ hidden units.

\begin{table}[t]
\centering
 \begin{tabular}{ lllllllll rr  }
  \hline
  $W_1$&$D_1$&$P_1$&$W_2$&$D_2$&$P_2$&$O$&$W_3$&$D_3$&loss&MAE \\ 
  \hline
32 & 2 & 0.25 & -  & - & -   & 128 & 64 & 2 & 7.36e-5 & 4.33e-3\\
32 & 2 & 0.25 & 32 & 2 & 0.5 & 128 & 64 & 2 & 9.28e-5 & 5.40e-3\\
32 & 2 & 0.0  & -  & - & -   & 128 & 64 & 3 & 7.85e-5 & 4.96e-3\\
32 & 2 & 0.25 & -  & - & -   & 128 & 64 & 3 & 7.72e-5 & 4.94e-3\\
32 & 2 & 0.5  & -  & - & -   & 128 & 64 & 3 & 7.86e-5 & 5.16e-3\\
32 & 2 & 0.25 & -  & - & -   & 256 & 64 & 3 & 8.19e-5 & 5.13e-3\\
32 & 2 & 0.25 & 64 & 2 & 0.5 & 128 & 64 & 4 & 1.88e-4 & 9.68e-3\\
  \hline
 \end{tabular}
\caption{Computed loss (MSE) and MAE for different ANNs architectures \textcolor{ultramarine}{trained with dataset 1. The quantities $W_1$, $D_1$, $P_1$, $W_2$, $D_2$, $P_2$, $O$, $W_3$ and $D_3$ are defined in Section~\ref{sec:aaarchit}.} The batch size is $32$, normalization \textcolor{reviewV2}{\texttt{std+id} (see Table~\ref{tab:norm_modes})}, training lasts $500$ epochs and the optimizer is the Adam algorithm \textcolor{ultramarine}{(with default Tensorflow learning rate)}. }
\label{tab:general_arch}
\end{table}

\begin{figure}[t]
\begin{center}
	\makebox[\textwidth][c]{\includegraphics[width=1.2\textwidth,trim={0 0.4cm 0 0},clip]{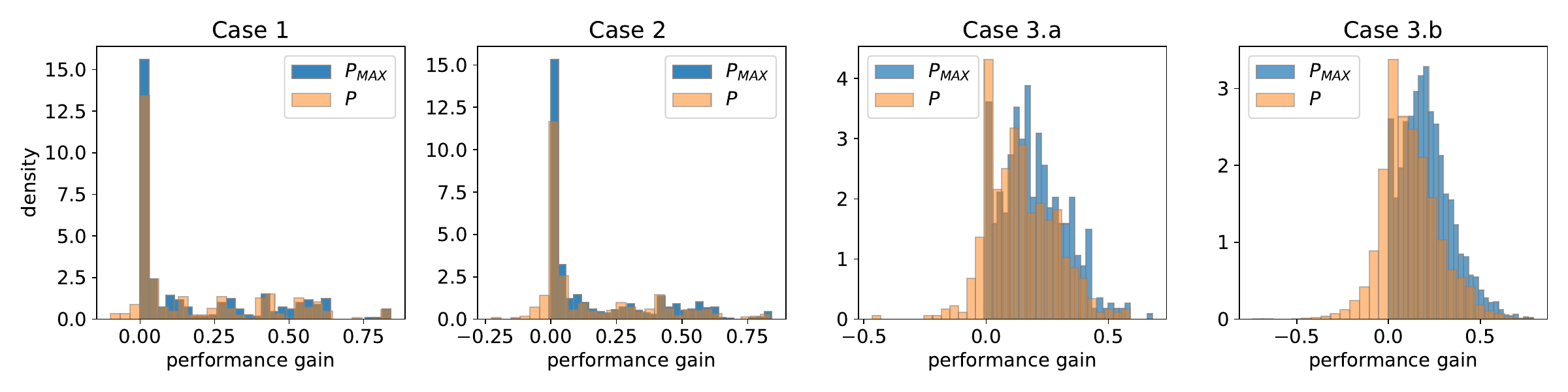}}
\end{center}
\caption{\textcolor{ultramarine}{Histogram of the performance of the AMG-ANN method: for each test case we evaluate the performance of the model $P$ in orange and the best performance $P_{\textnormal{MAX}}$ in blue (as defined in Eq.~\bref{eq:perf}). \textcolor{reviewV2}{\textit{From left to right.}} Performance for Case 1 of Section~\ref{sec:case_1}. Performance for the Case 2 of Section~\ref{sec:case_2}. Performance for the Case \textcolor{reviewV2}{3.a and 3.b} of Section~\ref{sec:case_3}.}}
\label{fig:perf_hist}
\end{figure}

\textcolor{ultramarine}{
In Table~\ref{tab:perf1} we report the performance indexes for the first six models of Table~\ref{tab:general_arch} (see Section~\ref{sec:model_eval_def}). We have chosen as architecture for our model the one reported in the first row of Table~\ref{tab:general_arch}. We trained this model for up to $1000$ epochs and employed early stopping. As result, it has a loss of $6.36 \cdot 10^{-5}$. We show the computed performance in Table~\ref{tab:perf_summary} (first row). Figure~\ref{fig:perf_hist} (left) shows an histogram of the performance gain $P$. We observe that in $20\%$ of the cases we have a performance gain $P \geq 43\%$.}

\begin{table}[t]
\centering
{\color{ultramarine}\begin{tabular}{ ccccccc  }
  \hline
  $PB$ & \multicolumn{2}{c}{$P$ (avg/median)} & \multicolumn{2}{c}{$P / P_\textnormal{MAX}$ (avg/median)} & \multicolumn{2}{c}{$P_{<0}$ (avg/median)}  \\ 
  \hline
  92.96\% & 16.06\% & 24.69\% & 81.31\% & 97.51\% & -3.928\% & -3.090\% \\
  90.36\% & 16.63\% & 33.74\% & 85.95\% & 97.14\% & -4.102\% & -2.123\% \\
  92.70\% & 16.64\% & 20.72\% & 81.18\% & 98.80\% & -2.363\% & -0.805\% \\
  92.96\% & 16.56\% & 17.97\% & 80.09\% & 97.56\% & -1.855\% & -0.826\% \\
  91.66\% & 16.31\% & 22.28\% & 82.39\% & 99.08\% & -3.462\% & -1.324\% \\
  91.14\% & 15.59\% & 22.32\% & 80.87\% & 96.80\% & -3.812\% & -1.094\% \\
  \hline
 \end{tabular}
 }
\caption{Evaluation of the performance of the first six models of Table~\ref{tab:general_arch}. The quantities $PB$, $P$, $P_{MAX}$ and $P_{<0}$ are defined at the end of Section~\ref{sec:aaarchit}.}
\label{tab:perf1}
\end{table}

\subsection{Test Case 2: an enhanced dataset}\label{sec:case_1}
In order to further test the robustness of the model to unseen data (i.e. test cases that are not in the training set), we test the prediction capabilities of the ANN on a new dataset. We call the latter dataset, ``dataset 2", while the one employed so far is called ``dataset 1". In particular, we solve the same model problem \bref{eq:2prob1} but with a different diffusion coefficient, defined as
\begin{equation}
    \mu(x,y) = 
\begin{dcases}
    10^{\varepsilon_2} \quad &\text{ if } (x,y) \in \Omega_{gray},\\
    10^{\varepsilon_1} & \text{ if } (x,y) \in \Omega_{white},
\end{dcases}
\label{eq:mu_patt2}
\end{equation}
where $\varepsilon_1$ and $\varepsilon_2$ are parameters to be chosen and $\Omega_{gray},\Omega_{white}$ is a partition of $\Omega$ as shown in Figure~\ref{fig:mu_pattern}.

\begin{figure}[t]
        \centering
        \includegraphics[trim={1.5cm 2cm 1.5cm 2cm},clip,width=0.9\textwidth]{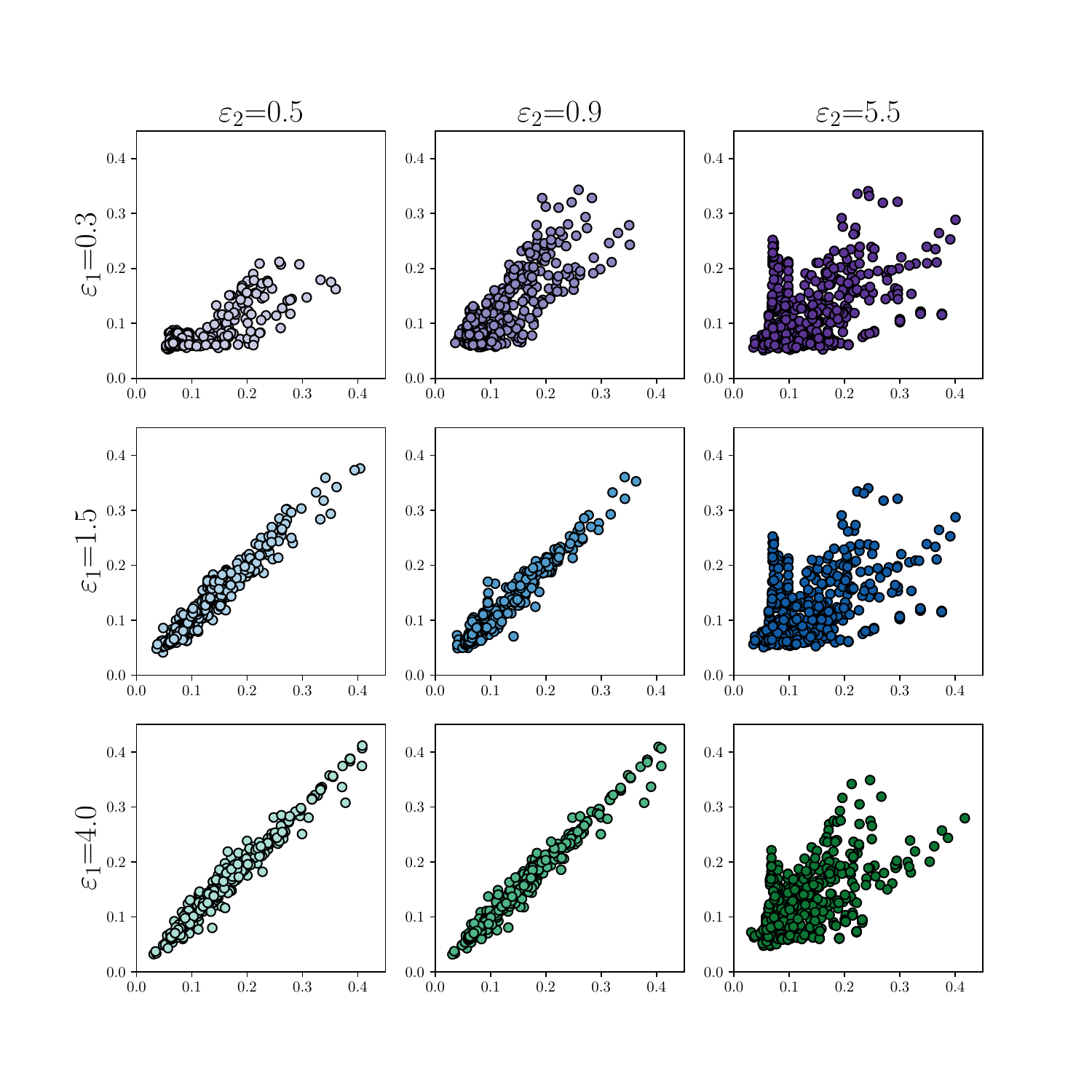}
\caption{Predictions on dataset 2 made by the model designed in Section~\ref{sec:aaarchit} and trained with dataset 1. Each plot \textcolor{reviewV2}{show data obtained} from a fixed combination of $\varepsilon_1$ and $\varepsilon_2$. Namely, $\varepsilon_1$ is constant for subplots on the same row and $\varepsilon_2$ is constant for subplots on the same column of the plot grid. On the x-axis there is the true value of $\rho$, on the y-axis the predicted value.}
    \label{fig:adv_pred_nothing_and_model3_test2}
\end{figure}

In Figure~\ref{fig:adv_pred_nothing_and_model3_test2}, we show the performance of the model that we designed in the previous sections, trained \textcolor{ultramarine}{with only} dataset 1. In particular, we choose randomly three values for $\varepsilon_1$ and three values for $\varepsilon_2$. We can see how the predictions maintain accuracy in some scenarios, but fail in other ones.

\begin{table}[t]
\centering

 \begin{tabular}{ llllll rr  }
  \hline
  $W_1$&$D_1$&$P_1$&$O$&$W_3$&$D_3$&loss&MAE \\ 
  \hline
32 & 2 & 0.05 & 128 & 64 & 3 & 1.75e-4 & 8.60e-3\\
32 & 2 & 0.00 & 128 & 64 & 4 & 1.63e-4 & 8.03e-3\\
32 & 2 & 0.25 & 128 & 64 & 3 & 1.58e-4 & 8.14e-3\\
32 & 2 & 0.25 & 128 & 64 & 4 & 1.35e-4 & 7.34e-3\\
32 & 2 & 0.50 & 128 & 64 & 3 & 1.55e-4 & 7.98e-3\\
  \hline
 \end{tabular}

\caption{\textcolor{ultramarine}{Computed loss (MSE) and MAE for different ANNs architectures (with one convolutional layer) trained with dataset 3. The quantities $W_1$, $D_1$, $P_1$, $O$, $W_3$ and $D_3$ are defined in Section~\ref{sec:aaarchit}. The batch size is $32$, \textcolor{reviewV2}{\texttt{normalize}=\texttt{std+id} (see Table~\ref{tab:norm_modes})}, training lasts $200$ epochs and the optimizer is the Adam algorithm (with default Tensorflow learning rate)}.}
\label{tab:adv_1_lay_dropout}
\end{table}

\begin{table}[t]
\centering
 \begin{tabular}{ lllllllll rr  }
  \hline
  $W_1$&$D_1$&$P_1$&$W_2$&$D_2$&$P_2$&$O$&$W_3$&$D_3$&loss&MAE \\ 
  \hline
16 & 2 & 0.25 & 12 & 2 & 0.25 & 128 & 64 &  4 & 1.72e-4 & 8.34e-3\\
16 & 2 & 0.25 & 12 & 2 & 0.50 & 256 & 256 & 3 & 1.56e-4 & 8.10e-3\\
16 & 2 & 0.25 & 16 & 2 & 0.25 & 128 & 64 &  4 & 1.68e-4 & 8.28e-3\\
16 & 2 & 0.25 & 32 & 2 & 0.50 & 128 & 64 &  4 & 1.61e-4 & 8.23e-3\\
32 & 2 & 0.25 & 16 & 2 & 0.25 & 128 & 64 &  3 & 1.55e-4 & 8.05e-3\\
32 & 2 & 0.25 & 16 & 2 & 0.25 & 128 & 64 &  4 & 1.55e-4 & 7.85e-3\\
32 & 2 & 0.25 & 32 & 2 & 0.50 & 128 & 64 &  4 & 1.75e-4 & 8.64e-3\\
  \hline
 \end{tabular}

\caption{\textcolor{ultramarine}{Computed loss (MSE) and MAE for different ANNs architectures (with two convolutional layers) trained with dataset 3. The quantities $W_1$, $D_1$, $P_1$, $W_2$, $D_2$, $P_2$, $O$, $W_3$ and $D_3$ are defined in Section~\ref{sec:aaarchit}. \textcolor{reviewV2}{The} hyperparameters are the same of Table~\ref{tab:adv_1_lay_dropout}}.}
\label{tab:adv_2_lay}
\end{table}

\begin{table}[t]
\centering
 \begin{tabular}{ llllll rr  }
  \hline
  $W_1$&$D_1$&$P_1$&$O$&$W_3$&$D_3$&loss&MAE \\ 
  \hline
32 & 2 & 0.25 & 128 &  64 & 3 & 1.71e-5 & 8.19e-3\\
32 & 2 & 0.25 & 128 &  64 & 4 & 1.35e-5 & 7.34e-3\\
32 & 2 & 0.25 & 128 &  64 & 5 & 1.48e-5 & 7.88e-3\\
32 & 2 & 0.25 & 128 & 128 & 3 & 1.51e-5 & 7.60e-3\\
32 & 2 & 0.25 & 128 & 128 & 4 & 1.43e-5 & 7.62e-3\\
32 & 2 & 0.25 & 256 &  64 & 5 & 1.50e-5 & 7.86e-3\\
32 & 2 & 0.25 & 256 & 256 & 3 & 1.52e-5 & 7.80e-3\\
32 & 2 & 0.25 & 512 & 128 & 3 & 1.60e-5 & 8.01e-3\\
32 & 3 & 0.25 & 128 &  64 & 4 & 1.51e-5 & 7.64e-3\\
32 & 3 & 0.25 & 128 & 128 & 3 & 1.48e-5 & 7.79e-3\\
32 & 3 & 0.25 & 128 & 128 & 4 & 1.34e-5 & 7.18e-3\\
32 & 3 & 0.25 & 256 & 256 & 3 & 1.56e-5 & 8.10e-3\\
  \hline
 \end{tabular}
\caption{\textcolor{ultramarine}{Computed loss (MSE) and MAE for different ANNs architectures (with one convolutional layer) trained with dataset 3. The quantities $W_1$, $D_1$, $P_1$, $O$, $W_3$ and $D_3$ are defined in Section~\ref{sec:aaarchit}. We change only the hyperparameters of the dense layers \textcolor{reviewV2}{as specified in columns $W_3$ and $D_3$}. \textcolor{reviewV2}{All the other} hyperparameters are the same \textcolor{reviewV2}{as those reported in} Table~\ref{tab:adv_1_lay_dropout}}.}
\label{tab:adv_1_lay}
\end{table}

\begin{table}[t]
\centering
 \begin{tabular}{ llllll rr  }
  \hline
  $W_1$&$D_1$&$P_1$&$O$&$W_3$&$D_3$&loss&MAE \\ 
  \hline
16 & 4 & 0.25 & 128 & 128 & 4 & 1.32e-4 & 7.29e-3\\
16 & 4 & 0.50 & 128 & 128 & 4 & 1.56e-4 & 7.98e-3\\
16 & 5 & 0.25 & 128 & 128 & 4 & 1.54e-4 & 7.77e-3\\
16 & 3 & 0.50 & 128 & 128 & 4 & 1.51e-4 & 7.86e-3\\
24 & 2 & 0.25 & 128 & 128 & 4 & 1.53e-4 & 7.71e-3\\
24 & 2 & 0.50 & 128 & 128 & 4 & 1.52e-4 & 7.76e-3\\
24 & 3 & 0.25 & 128 & 128 & 4 & 1.40e-4 & 7.34e-3\\
24 & 4 & 0.50 & 128 & 128 & 4 & 1.60e-4 & 7.98e-3\\
32 & 3 & 0.25 & 128 & 128 & 4 & 1.34e-5 & 7.18e-3\\
32 & 3 & 0.50 & 128 & 128 & 4 & 1.47e-5 & 7.75e-3\\
40 & 2 & 0.25 & 128 & 128 & 4 & 1.27e-4 & 7.30e-3\\
40 & 3 & 0.25 & 128 & 128 & 4 & 1.32e-4 & 7.18e-3\\
  \hline
 \end{tabular}
\caption{\textcolor{ultramarine}{Computed loss (MSE) and MAE for different ANNs architectures (with one convolutional layer) trained with dataset 3. The quantities $W_1$, $D_1$, $P_1$, $O$, $W_3$ and $D_3$ are defined in Section~\ref{sec:aaarchit}. Here, we change only the hyperparameters of the convolutional layer. Not specified hyperparameters are the same of Table~\ref{tab:adv_1_lay_dropout}}}
\label{tab:adv_1_lay_basis}
\end{table}

We proceed to show how the model behaves when the training is instead done with training samples from both datasets. Dataset~2 contains $5184$ entries, we define the test set to be the union of the $20\%$ of dataset~2 and the $50\%$ of the dataset~1. In this way, the union of the training and validation set contains $4800$ datapoints from the dataset~1 and $4147$ from dataset~2. The ratio between the number of entries of the validation set and the training set is defined to be $1:3$. We call this combination dataset~3. The aim is to have a balanced training dataset in which each definition of $\mu$ is equally represented.

If not otherwise stated, we stop the training at $200$ epochs. As shown in Table~\ref{tab:adv_1_lay_dropout}, dropout improves the training, thus it will be employed in all the models.
We have also tried employing batch normalization as a regularization technique on some of these models and a deeper model with three convolutional layers but it did not lead to any significant improvements. This can be explained by the fact that batch normalization effectiveness is most evident in very deep models; see \cite{he2016deep}.

\begin{figure}[t]
        \centering
        \includegraphics[trim={1.5cm 2cm 1.5cm 2cm},clip,width=0.9\textwidth]{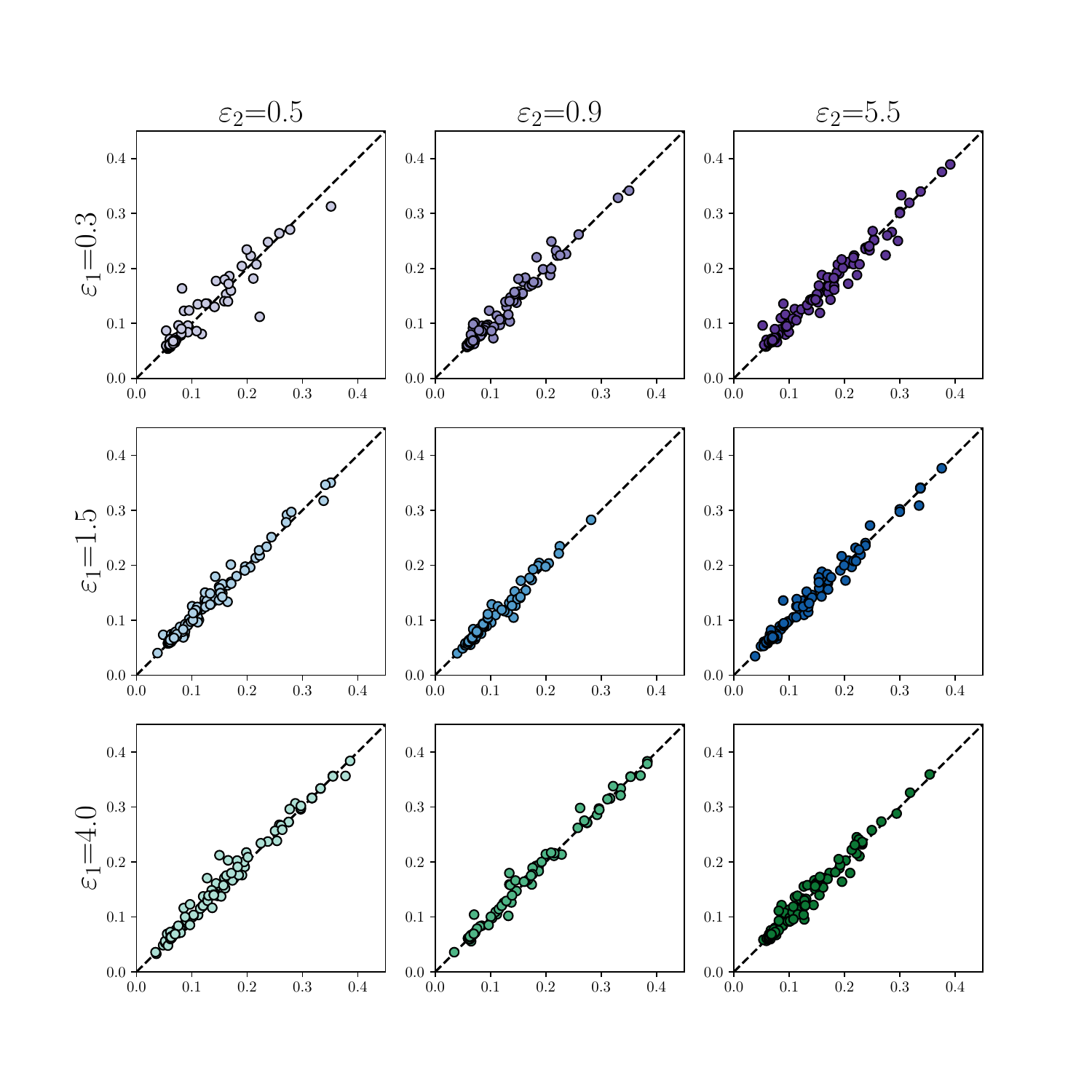}
\caption{Prediction of the model trained with dataset 3 on the test dataset 2. Each plot \textcolor{reviewV2}{show data obtained} from a fixed combination of $\varepsilon_1$ and $\varepsilon_2$. Namely, $\varepsilon_1$ is constant for subplots on the same row and $\varepsilon_2$ is constant for subplots on the same column of the plot grid. On the x-axis there is the true value of $\rho$, on the y-axis the predicted value.}
    \label{fig:adv_pred_nothing_and_model3_test2bis}
\end{figure}

Table~\ref{tab:adv_2_lay} shows training of models with two convolutional layers. By comparing it with Table~\ref{tab:adv_1_lay}, where the MSE and MAE are reported for different ANN architectures with one layer, it is possible to notice that models with only one layer achieve lower loss. From Table~\ref{tab:adv_1_lay}, it is also possible to appreciate that the model that in the previous section achieved the lowest loss is not the same in this case. In particular, a deeper model performs better. This is not surprising since this means that we need a more complex model to explain the data, and indeed we are using a more diversified dataset. In Table~\ref{tab:adv_1_lay_basis}, we repeat the same test case for different architectures of the convolutional layer. Two applications of convolution with 40 hidden units seems to be the best choice. The architecture that we choose for the model is the second to last of Table~\ref{tab:adv_1_lay_basis}. Employing training with batch size 32, the Adam optimizer and early stopping (up to 1000 epochs), we obtain a loss on the test ``dataset~3" of $8.83 \cdot 10^{-5}$ and MAE $4.84 \cdot 10^{-3}$. On the test ``dataset 1" we achieve a loss of $7.73 \cdot 10^{-5}$ and a MAE of $4.48 \cdot 10^{-3}$ and on test ``dataset 2" we obtained a loss of $1.39 \cdot 10^{-4}$ and a MAE of $6.51 \cdot 10^{-3}$. The predictions are reported in Figure~\ref{fig:adv_pred_nothing_and_model3_test2bis}.
\textcolor{ultramarine}{Figure~\ref{fig:perf_hist} (center) shows an histogram of the performance gain $P$. In particular, in $20\%$ of the cases $P\geq37\%$. Table~\ref{tab:perf_summary} (second row) summarizes the results.}

\begin{remark}
\textcolor{ultramarine}{
We observe that the difference in predictions showed in
Figures~\ref{fig:adv_pred_nothing_and_model3_test2} and~\ref{fig:adv_pred_nothing_and_model3_test2bis} are significant. This begs the question of how assessing the quality of a training dataset. Unfortunately, it is very hard to assess the quality of the dataset a-priori. Indeed, this question is equivalent to predict the neural network generalization, which is still an open question in the field of ML. However, in the case of very large datasets, it could be useful to train a small ANN on a small sample of the dataset and validate the results.
}
\end{remark}

\color{ultramarine}
\subsection{Test Case 3: diffusion coefficient with \textcolor{reviewV2}{different values on each tile}}\label{sec:case_2}
In this case, we test the prediction capabilities of AMG-ANN whenever the diffusion coefficient shows a more complicated pattern. Namely, we generalize the definition of the diffusion coefficient $\mu$ in the following way. Let $\texttt{size} \in \mathbb N$ be a positive integer indicating the size of the pattern (e.g. in Figure~\ref{fig:mu_pattern}, from left to right, \texttt{size} $= 2, 2, 4, 4$ since the patterns are two stripes, $2\times 2$ checkerboard, four stripes and $4 \times 4$ checkerboard, respectively). Let \texttt{mode} $= 1, 2$ be an integer that indicates if the pattern is at stripes or checkerboard like, respectively. These two parameters determine a partition $\{\Omega_i\}_{i=1, ..., {\tt{size}}^{\tt{mode}} }$ of the domain $\Omega$. On each element of the partition $\Omega_i$ we set $\mu(x,y) = 10^{(\boldsymbol{\varepsilon})_i}$ (constant), where $\boldsymbol{\varepsilon} \in \mathbb{R}^{\texttt{size}^\texttt{mode}}$ is a given vector of parameters. What just described is synthesized in Algorithm~\ref{algo:mu}. Notice that if \texttt{size} is not a power of two, the pattern of $\mu$ is not aligned with the mesh. 

We choose $f$ of Eq.~\bref{eq:2prob1} such that $u(x, y) = \sin(\texttt{size} \pi x / 2) \sin(\texttt{size} \pi y / 2)$ if \texttt{size} is odd and $u(x, y) = \cos(\texttt{size} \pi x / 2) \cos(\texttt{size} \pi y / 2)$ otherwise. Varying \texttt{size} $=2,...,10$, \texttt{mode} $= 1, 2$ and $\boldsymbol{\varepsilon}$ by sampling its component at random with uniform distribution in the interval $(-20, 20)$, we create a dataset with about 200~000 samples ($\theta$ and $h$ vary in the same range as defined in Section~\ref{sec:numeric_res}).
{\color{reviewV2}
\begin{table}[t]
\centering
{\color{reviewV2}\begin{tabularx}{1.0\textwidth}{ llX  }
  \hline
   Name & \#channels $c$ & list of \texttt{op}  \\ 
  \hline
  \texttt{sum} & 1  & \texttt{op}$(v_1, v_2) = v_1 + v_2$ \\
  \hline
  \texttt{max} & 1  & \texttt{op}$(v_1, v_2) = \max\{v_1, \abs{v_2}\}$ \\
  \hline
  \texttt{pp+np} & 2 &$\texttt{op}_1(v_1, v_2) = \max\{\max\{0, v_2\}, v_1\}$\hfill\textcolor{white}{,} $\texttt{op}_2(v_1, v_2) = \max\{\max\{0, -v_2\}, v_1\}$ \\
  \hline
  \texttt{pp+np+sum} & 3 &$\texttt{op}_1(v_1, v_2) = \max\{\max\{0, v_2\}, v_1\}$\hfill\textcolor{white}{,}\hfill$\texttt{op}_2(v_1, v_2) = \max\{\max\{0, -v_2\}, v_1\}$\hfill\textcolor{white}{,} $\texttt{op}_3(v_1, v_2) = v_1 + v_2$ \\
  \hline
 \end{tabularx}
 }
\caption{\color{reviewV2}Summary of the possible ways to obtain $\hat{\mathrm V}' \in \mathbb R^{m\times m \times c}$ by stacking together $\{ \mathrm{\hat{V}}_i= \texttt{normalize}(\texttt{pooling}(\mathrm A_h, m, \texttt{op}_i))\}_{i = 1, ..., c}$, where $c$ and $\{\texttt{op}_i\}_{i = 1, ..., c}$ are defined in the second and third columns, respectively.}
\label{tab:view_modes}
\end{table}
\subsubsection{Tuning of the hyperparameters of the pooling operation}\label{sec:view_hyper}

To tune the hyperparameters \texttt{op}, \texttt{normalize} and $m$ we consider a subset of the dataset consisting of about 15~000 samples. By stacking together $c$ views $\hat{\mathrm V}$ obtained with different \texttt{op} we can obtain a tensor $\hat{\mathrm V}' \in \mathbb R^{m\times m \times c}$ that can be interpreted as a multi-channel image. Namely, we consider four possible approaches summarized in Table~\ref{tab:view_modes}. The first two rows represent the most used approaches in the computer vision field, whereas the argument for the approach in the third row comes from Eq.\bref{eq:amg_interpolation} where you can see the relevance of splitting the values into the positive and negative part. The fourth approach just combines the information of the first three together. We also consider all the six possible choice of \texttt{normalize} reported in Table~\ref{tab:norm_modes} and eight values of $m = 30, 40, ..., 100$. Thus, there are 192 combinations of hyperparameter $m$, \texttt{op} and \texttt{normalize}. For each one of this choices we train a network with learning rate 0.001, batch size 32, 1 convolutional layer with no dropout and $W_1=32, D_1=2, W_3=64, D_3 = 3$. Figure~\ref{fig:hyper_view} shows the boxplot of the loss of the model trained for 100 epochs (top row) and 200 epochs (bottom row), grouped by hyperparameter category.
From these results it seems that the choices \texttt{normalize}=\texttt{log+id} and \texttt{normalize}=\texttt{log+avg} provide the best results. Concerning the choice of \texttt{op}$(\cdot, \cdot)$, it seems that \texttt{op}=\texttt{pp+np+sum} leads to better results. The view size $m$ seems to be inversely proportional to the loss. However, comparing $m=60$ and $m=100$ when the model is trained for 200 epochs, we notice that, on average, the loss is only about $10\%$ larger even if the view is $64\%$ smaller. 

\begin{figure}[t]
\begin{center}
	\makebox[\textwidth][c]{\includegraphics[width=1.2\textwidth,trim={0 1cm 0 0},clip]{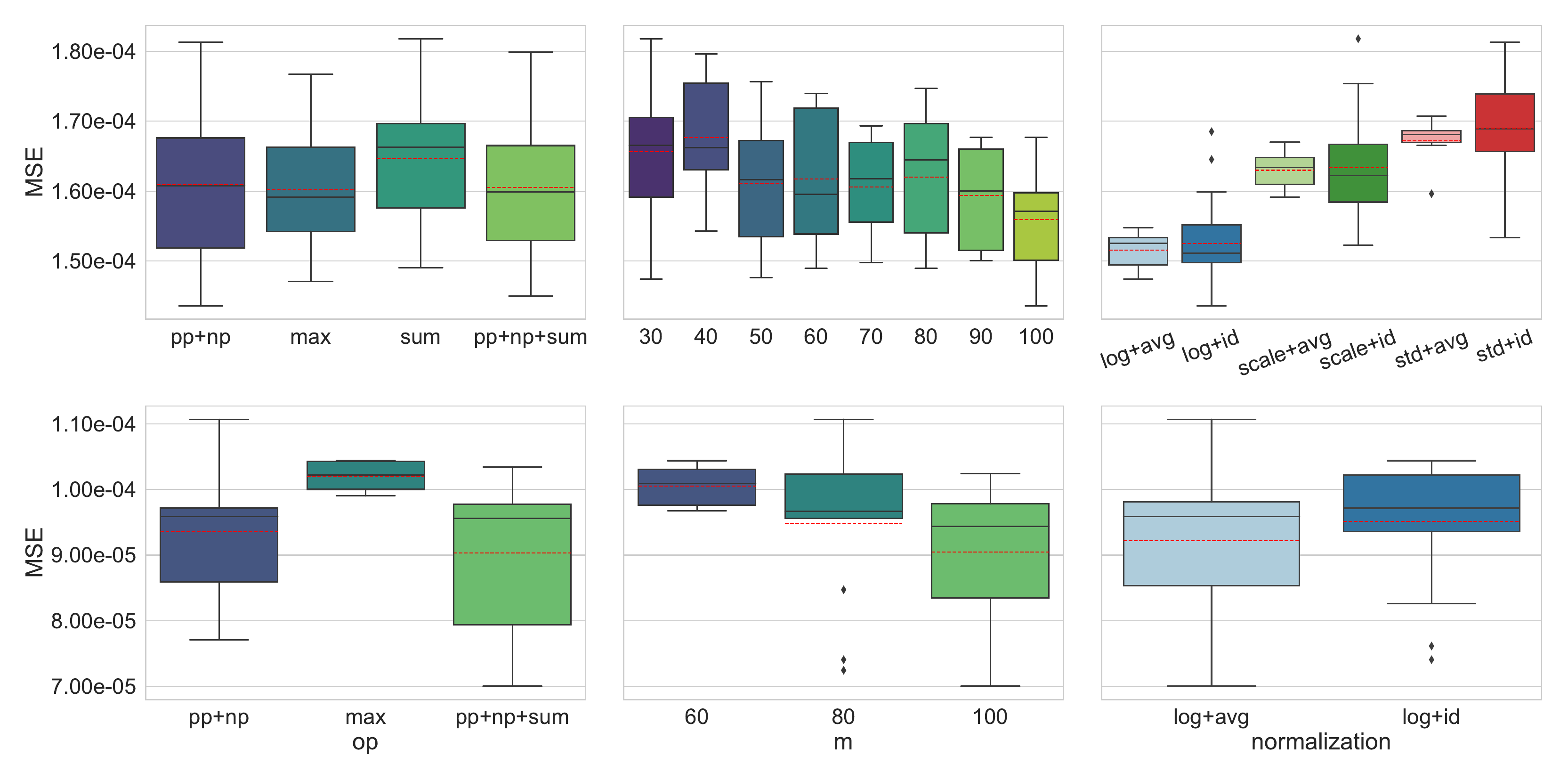}}
\end{center}
\caption{\textcolor{reviewV2}{Boxplot (with mean in red) representing the effect of different choices of the hyperparaments \texttt{op} (\textit{left}), $m$ (\textit{center}) and \texttt{normalize} (\textit{right}) on the test loss (MSE). The training is stopped after 100 epochs (\textit{top row}), 200 epochs (\textit{bottom row}).}}
\label{fig:hyper_view}
\end{figure}

\begin{table}[t]
\centering
\color{reviewV2}
 	\makebox[\textwidth][c]{\begin{tabular}{ lllllllll rr  }
  \hline
  \texttt{op}&$m$&\texttt{normalize}&$W_1$&$D_1$&$P_1$&$O$&$W_3$&$D_3$&loss&MAE \\ 
  \hline
{\texttt{pp+np+sum}} & 100 & \texttt{log+avg} & 40 & 3 & 0 & 128 & 128 & 5 & 4.31e-5 & 3.65e-3\\
{\texttt{pp+np+sum}} & 100 & \texttt{log+id}  & 40 & 3 & 0 & 128 & 128 & 5 & 3.87e-5 & 3.47e-3\\
{\texttt{pp+np+sum}} & 80  & \texttt{log+avg} & 40 & 3 & 0 & 128 & 128 & 5 & 4.07e-5 & 3.57e-3\\
  \hline
 \end{tabular}}
\caption{\textcolor{reviewV2}{Case 3.a: computed loss (MSE) and MAE for different choices of the hyperparameters \texttt{op}, $m$ and \texttt{normalize} hyperparameters. All the other hyperparameters are the same as those reported in Table~\ref{tab:adv_1_lay_dropout}}}
\label{tab:hyper_view}
\end{table}
\subsubsection{Test Case 3.a: Predictions on 15~000 samples}
We train our model on the reduced dataset containing ``only'' 15~000 samples. After tuning the architecture of the neural network, we train three different models for 500 epochs. Table~\ref{tab:hyper_view} shows the results. The choices \texttt{normalize}=\texttt{log+id} and \texttt{normalize}=\texttt{log+avg} produce similar results as before. We can conclude that $m=100$ seems better, but the improvement on the loss is only about 10\% w.r.t.\ $m=80$. Thus $m = 80$ could be a reasonable choice to save time both in the offline and online phase. The second model results into the lowest loss and MAE, its performance is reported in Table~\ref{tab:perf_summary}. Moreover, we observe that in 20\% of the cases it has a performance gain $P \geq 27\%$.
}

\textcolor{reviewV2}{
\subsubsection{Test Case 3.b: Predictions on 200~000 samples}
In this section we train our model with the full dataset containing 200~000 samples. We use the following choices as hyperparameters $m=50$, $\texttt{normalize}=\textnormal{\texttt{scale+avg}}$ and $\texttt{op} = \textnormal{sum}$ since they provide a good compromise between accuracy and efficiency for tuning the architecture of the ANN.} 

First, we check if using as loss function MSE or MAE makes any difference on the final performance: we fix two models architectures and train the ANNs changing only the loss function. Results are reported in Table~\ref{tab:mae_vs_mse}: there seems not to be any significant difference between the two.

We then perform an hyperparameter optimization on the model architecture: we train different models for 40 epochs and check their performance in terms of loss. We test different dropout probabilities, different sizes and depth for the convolution filters and different architectures for the dense part. We report some results in Table~\ref{tab:case2_models}, where we can see that larger models, w.r.t.\ the ones in Table~\ref{tab:adv_1_lay_basis}, perform better. This is expected since we have a rather larger dataset. 

We choose as architecture of our model the one reported in the seventh row of Table~\ref{tab:case2_models}. \textcolor{reviewV2}{Using as hyperparameters of the view \texttt{op}=\textnormal{\texttt{pp+np+sum}}, $m=100$, \texttt{normalize}=\texttt{log+avg}, the model is trained for 500 epochs (validation loss was still decreasing when training ended). It reaches a loss of $6.95 \cdot 10^{-5}$ MSE ($4.69 \cdot 10^{-3}$ MAE). We show its performance in Table~\ref{tab:perf_summary} (fourth row).} The histogram of the performance is reported in Figure~\ref{fig:perf_hist} (right). Even if its performance is not as good as the one of the model obtained in Section~\ref{sec:case_1}, we would like to point out that there is still a margin of improvement for the hyperparameters optimization step \textcolor{reviewV2}{and more epochs could be used}.

\begin{algorithm}[t]
\color{ultramarine}
    $j$ $\leftarrow$ $1$\;
    \For{$i\leftarrow 1$ \KwTo \textnormal{\texttt{mode}}}{
        $j$ $\leftarrow$ $j + \lfloor ((\mathbf x)_i + 1) \texttt{size} / 2 \rfloor \, \texttt{size}^{i-1}$    \;
    }
    $\mu$ $\leftarrow$ $10^{(\boldsymbol{\varepsilon})_j}$\;
    \caption{Diffusion coefficient $\mu$ in a given point $\mathbf{x} \in \Omega = (-1, 1)^2$, fixed the parameters \texttt{mode}, \texttt{size} and $\boldsymbol{\varepsilon} \in \mathbb{R}^{\texttt{size}^\texttt{mode}}$. \newline $\mu = \mu(\mathbf{x}; \texttt{mode}, \texttt{size}, \boldsymbol \varepsilon)$}
    \label{algo:mu}
\end{algorithm}

\begin{table}[t]
\centering
\addtolength{\leftskip} {-2cm}
\addtolength{\rightskip}{-2cm}
\color{ultramarine}\begin{tabular}{ l cccccccc  }
\cline{2-9}
& loss & MAE & MSE & $PB$ & \multicolumn{2}{c}{$P$ (avg/median)} & \multicolumn{2}{c}{$P_{<0}$ (avg/median)}  \\
\hline
Model 1 & MSE & 1.00e-3 & 2.95e-4 & 62.3\% & 3.14\% & 2.67\% & -13.1\% & -7.91\%\\
Model 2 & MSE & 9.67e-4 & 2.71e-4 & 62.5\% & 3.55\% & 2.19\% & -12.3\% & -7.74\%\\
Model 1 & MAE & 9.35e-4 & 2.83e-4 & 63.7\% & 3.69\% & 3.56\% & -13.4\% & -7.47\%\\
Model 2 & MAE & 1.00e-3 & 3.18e-4 & 61.8\% & 2.85\% & 2.86\% & -14.1\% & -8.69\%\\
\hline
\end{tabular}
\caption{\color{ultramarine}Comparison of the performance of two fixed models when using MSE or MAE as a loss. The two models are the second and third of Table~\ref{tab:case2_models}. They were trained for 20 epochs. batch size 32, default Tensorflow learning rate of Adam. The quantities $PB$, $P$ and $P_{<0}$ are defined in Section~\ref{sec:aaarchit}.}
\label{tab:mae_vs_mse}
\end{table}

\begin{table}[htbp]
\centering
 \color{ultramarine}\begin{tabular}{ llllll rr  }
  \hline
  $W_1$&$D_1$&$P_1$&$O$&$W_3$&$D_3$&loss&MAE \\ 
  \hline
32&2&0.25&128&128&4&3.64e-3&1.24e-2\\
32&3&0.25&128&128&4&2.22e-3&9.35e-3\\
32&4&0.25&128&128&4&2.16e-3&9.16e-3\\
32&5&0.25&128&128&4&2.26e-3&9.48e-3\\
40&2&0.25&128&128&4&2.45e-3&9.90e-3\\
40&3&0.25&128&128&4&2.12e-3&9.22e-3\\
40&3&0.25&128&128&5&1.98e-3&8.75e-3\\
40&4&0.25&128&128&4&2.22e-3&9.31e-3\\
64&3&0.25&128&128&4&2.19e-3&9.31e-3\\
64&4&0.25&128&128&4&3.16e-3&1.01e-2\\
  \hline
 \end{tabular}
\caption{\color{ultramarine}Computed loss (MSE) and MAE on for different architectures of the ANN trained on the dataset of Section~\ref{sec:case_2} for 40 epochs. The quantities $W_1$, $D_1$, $P_1$, $O$, $W_3$ and $D_3$ are defined in Section~\ref{sec:aaarchit}. Not specified hyperparameters are the same of Table~\ref{tab:mae_vs_mse}}
\label{tab:case2_models}
\end{table}

\subsection{Test Case 4: the stationary Stokes problem}\label{sec:case_3}
\begin{figure}[t]
    \begin{center}
    \includegraphics[trim={10cm 1cm 10cm 1cm},clip,width=1.0\textwidth]{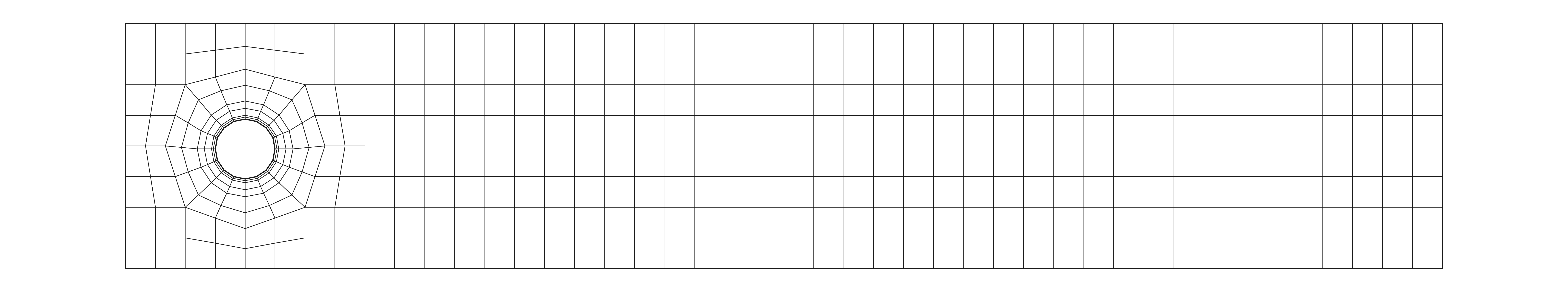}
    \end{center}
    \caption{\color{reviewV2}Representation of the domain of the Stokes problem~\bref{eq:stokes} with a quadrilateral mesh $\mathcal{T}_h$.}
    \label{fig:stokes_mesh}
\end{figure}
Finally, we show how we can extend the prediction capabilities of the ANN trained in the previous section in the case we consider a different model problem. 
\textcolor{reviewV2}{Namely, we consider the model problem defined in Eq.~\bref{eq:stokes}.}
Namely, we consider the channel flow around a cylinder with parabolic inflow profile in $\Omega =(0,2.2)\times (0,0.41) \setminus B_r(0.2,0.2) \subset \mathbb{R}^2$ with $r=0.05$ where $B_r(x, y)$ is the ball of radius $r$ centered in $(x, y) \in \mathbb{R}^2$. \textcolor{reviewV2}{We define $\Gamma_0 =[0,2.2]\times \{0, 0.41\} \cup\partial B_r(0.2,0.2), \Gamma_{in} = 0\times[0,0.41], \Gamma_{out} = 2.2\times[0,0.41]$.} We report in Figure~\ref{fig:stokes_mesh} a graphical representation of the computational domain together with a computational mesh. \textcolor{reviewV2}{We solve its discrete formulation \bref{eq:stokes_algebraic}} with MINRES preconditioned with
\begin{eqnarray*} \mathrm{P}^{-1} = \left(\begin{array}{cc} \mathrm{A}_h & 0 \\ 0 & \mathrm{M}_h \end{array}\right) ^{-1} = \left(\begin{array}{cc} \mathrm{A}_h^{-1} & 0 \\ 0 & \mathrm{M}_h^{-1} \end{array}\right), \end{eqnarray*}
where $\mathrm{M}_h$ is the mass matrix in the pressure space. For the approximation of the velocity block A${}_h$ we will perform a single AMG V-cycle: our AMG-ANN will be applied to this block.

\begin{figure}[t]  
\begin{centering}
\begin{subfigure}{0.45\textwidth}
\includegraphics[width=\linewidth]{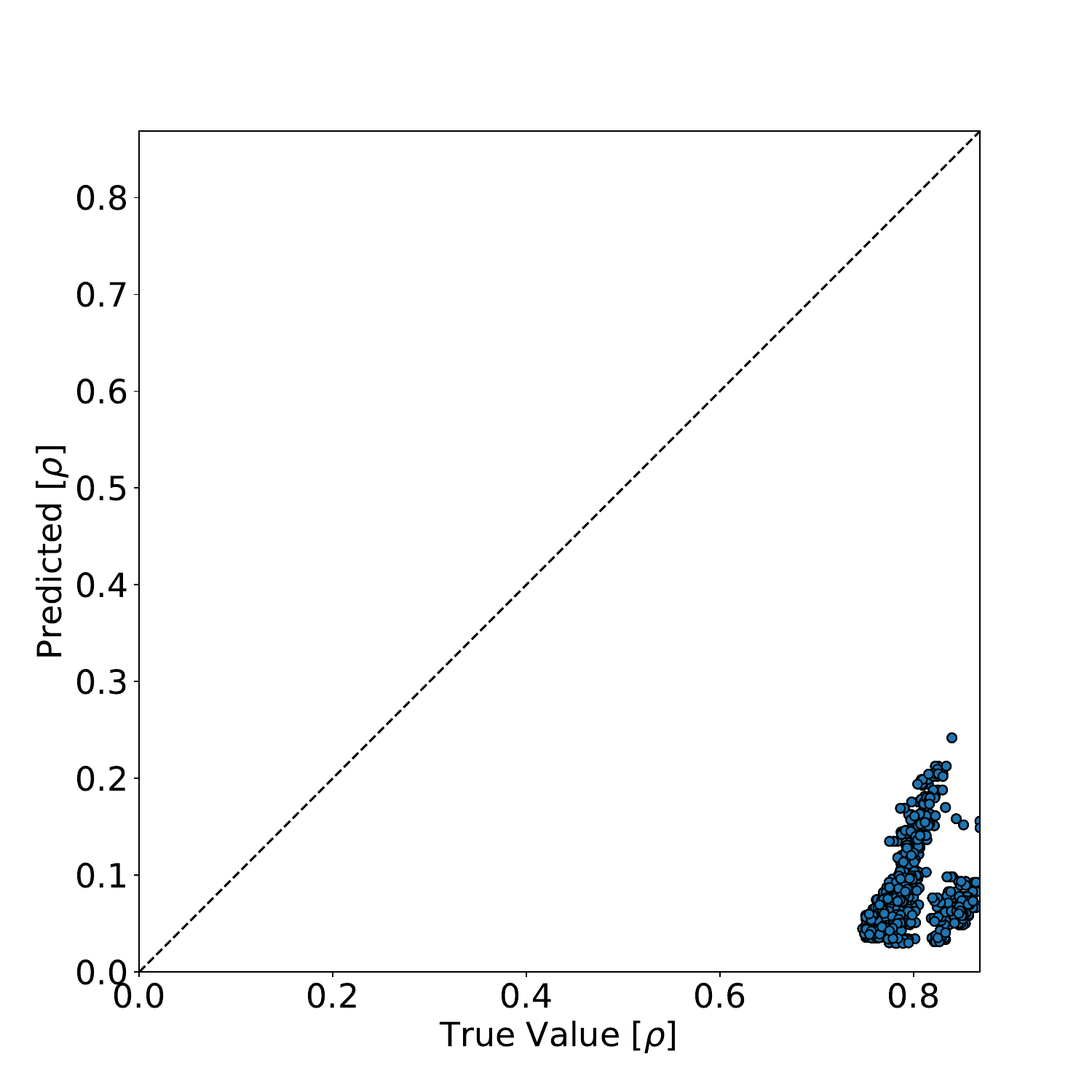}
\end{subfigure}
\begin{subfigure}{0.45\textwidth}
\includegraphics[width=\linewidth]{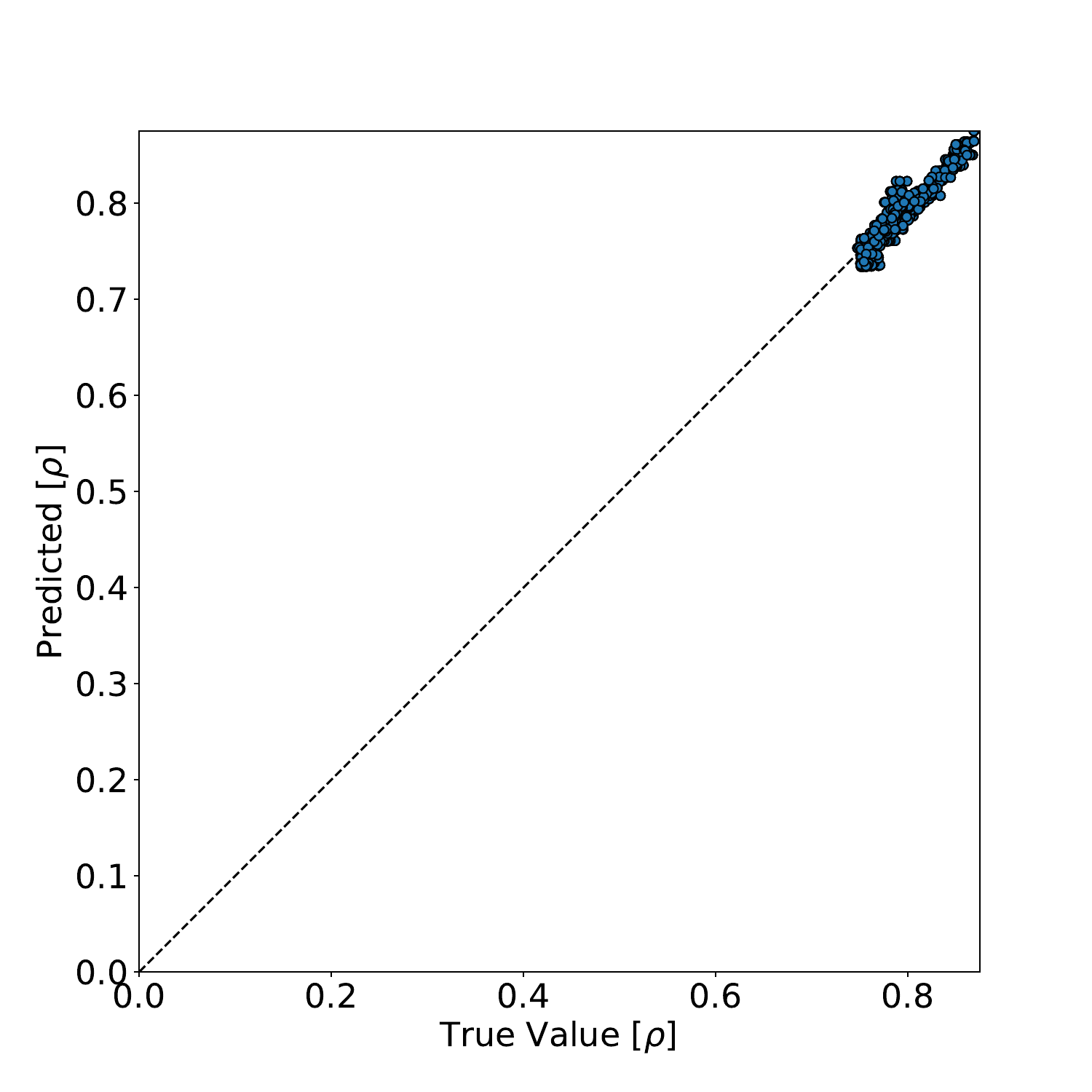}
\end{subfigure}
\caption{\color{ultramarine}\textit{Left.} \textcolor{reviewV2}{Predictions on the Stokes problem dataset defined in Section~\ref{sec:case_3} of the model trained in Section~\ref{sec:case_2} (i.e.\ with only Laplace problems)}. \textit{Right.} The predictions of the same model after training for an extra 20 epochs on a dataset composed at 50\% by Stokes problems and 50\% by samples of the dataset defined in Section~\ref{sec:case_2}. Not specified hyperparameters are the same of Table~\ref{tab:mae_vs_mse}. } 
\label{fig:stokestl}
\end{centering}
\end{figure}

Then, we build a dataset with 3600 samples by varying $\nu = 0.001, 0.1, 10$, $U = 0.00001, 0.001, 0.1, 10$, $h = 0.14 / 2^i$ with $i = 0, 1, 2, 3, 4, 5$. Our baseline is the model trained in the Section~\ref{sec:case_2}. Using transfer learning, we would like to extend its applicability range to the Stokes problem. To this end, using the weights (and architecture) of the pre-trained model obtained in Section~\ref{sec:case_2}, we train it with a dataset composed at 50\% by Stokes problem samples and 50\% by data from the dataset we build in Section~\ref{sec:case_2} (i.e.\ the dataset used to pre-trained the model). Indeed, we keep some of the old data to avoid catastrophic forgetting \cite{catastrophicforgetting}, i.e.\ the tendency of ANNs to forget how to perform a task upon learning new information. We train the model only for 20 epochs, we employ early stopping and a 80\%-20\% training-validation split. In Figure~\ref{fig:stokestl} we report the prediction of the model before and after this step of training. At the end of the training the model has a loss of $1.52\cdot 10^{-4}$ MSE (MAE $9.09\cdot 10^{-3}$). We report its performance in Table~\ref{tab:perf_summary}. We remark that, in this case, a performance gain $P$ of 3\% (in terms of $\rho$) results into about 30\% performance gain in terms of elapsed CPU time. Meanwhile on the full original dataset we still maintain a $PB > 50\%$ and $P > 0$ both in mean and median. Thus, with a small computational effort, we succeeded in extending the knowledge of our model.

\begin{table}[t]
  \centering
  \addtolength{\leftskip} {-2cm}
 \addtolength{\rightskip}{-2cm}
 \color{ultramarine}\begin{tabular}{ lccccccc  }
    \cline{2-8}
     & $PB$ & \multicolumn{2}{c}{$P$ (avg/median)} & \multicolumn{2}{c}{$P / P_\textnormal{MAX}$ (avg/median)} & \multicolumn{2}{c}{$P_{<0}$ (avg/median)}  \\ 
    \hline
    Case 1 & 96.61\% & 17.27\% & 4.11\% & 86.96\% & 99.17\% &  -1.65\% & -0.63\% \\
    Case 2    & 89.13\% & 14.92\% & 2.81\% & 80.09\% & 96.01\% &  -3.72\% & -2.04\% \\
    \textcolor{reviewV2}{Case 3.a}    & \textcolor{reviewV2}{88.60\% }&  \textcolor{reviewV2}{15.03\%} & \textcolor{reviewV2}{13.60\%} & \textcolor{reviewV2}{72.06\%} & \textcolor{reviewV2}{82.24\%} & \textcolor{reviewV2}{-5.47\%} & \textcolor{reviewV2}{-2.25\%} \\
    \textcolor{reviewV2}{Case 3.b}    & \textcolor{reviewV2}{80.65\% }&  \textcolor{reviewV2}{12.60\%} & \textcolor{reviewV2}{10.94\%} & \textcolor{reviewV2}{65.82\%} & \textcolor{reviewV2}{69.79\%} & \textcolor{reviewV2}{-6.00\%} & \textcolor{reviewV2}{-3.91\%} \\
    Case 4    & 90.16\% &  2.45\% & 1.24\% & 79.62\% & 90.76\% &  -0.36\% & -0.39\% \\
    \hline
   \end{tabular}
  \caption{\color{ultramarine}Evaluation of the performance of the best model for each case. The models and cases are defined in Section~\ref{sec:aaarchit}, Section~\ref{sec:case_1} and Section~\ref{sec:case_2}, Section~\ref{sec:case_3} respectively. The quantities $PB$, $P$, $P_{MAX}$ and $P_{<0}$ are defined at the end of Section~\ref{sec:aaarchit}.}
  \label{tab:perf_summary}
  \end{table}

\color{black}
\section{Conclusions}\label{sec:conc}
In this work, we developed an ANN-based approach to enhance the computational efficiency of the AMG methods, i.e. to accelerate their performance. In particular, we accurately predicted the value of the strong threshold parameter $\theta$ that maximizes the performance with respect to the matrix $\mathrm{A}_h$ of the linear system to be solved. In order to be able \textcolor{ultramarine}{to} apply the model independently of the matrix of the linear system, we introduced a pooling operator. We measured the efficiency of the AMG method using the approximate convergence factor and we designed a model that predicts its value. In this way, we are able to choose the strong threshold parameter that minimizes the predicted approximated convergence factor.
Moreover, we have shown that, as expected, the approximated convergence factor of the AMG method is strictly correlated to the elapsed CPU time during the application of the AMG method to the linear system solution, thus demonstrating that it provides a good measure of the performance of the solver. This a priori, optimal selection of the strong threshold parameter allows us to efficiently \textcolor{reviewV2}{choose} a value of $\theta$ that significantly decreases the elapsed CPU time with respect to the ``classical" value. \textcolor{reviewV2}{We introduced a set of indicators to measure the performance of our model: we show that if the dataset is smaller than 15~000 samples our model is better than using the literature value of $\theta$ in about 90\% of the cases with on average a 15\% gain in performance. On the other hand, from our computations, it seems that in case of large datasets more work on the tuning of the hyperparameter is needed. Finally, we have addressed a generalization test case moving from the elliptic scalar differential problem to the Stokes system exploiting the so called transfer learning. The preliminary results are encouraging and further investigation will be the subject of further research.} 

Possible further developments \textcolor{reviewV2}{also} include: using the ANN to optimize the value of other AMG parameters, \textcolor{ultramarine}{such as the maximum row sum parameter, the choice of whether using W-cycles or V-cycles or the number of levels of aggressive coarsening. Possible other improvements include} \textcolor{reviewV2}{tuning of the hyperparameters of the model};
testing \textcolor{ultramarine}{a wider range of differential models}, in particular in three-dimensional \textcolor{reviewV2}{configuration}.

\bibliographystyle{abbrv}
\bibliography{references}

\begin{thebibliography}{10}

\bibitem{abhyankar2018petsc}
S.~Abhyankar, J.~Brown, E.~M. Constantinescu, D.~Ghosh, B.~F. Smith, and
  H.~Zhang.
\newblock {PETSc/TS}: A modern scalable {ODE/DAE} solver library.
\newblock {\em arXiv preprint arXiv:1806.01437}, 2018.

\bibitem{antoniettimanuzzi2}
P.~F. Antonietti, F.~Dassi, and E.~Manuzzi.
\newblock Machine learning based refinement strategies for polyhedral grids
  with applications to {V}irtual {E}lement and polyhedral discontinuous
  {Galerkin} methods.
\newblock {\em Journal of Computational Physics}, 2022, in press.

\bibitem{antoniettimanuzzi}
P.~F. Antonietti and E.~Manuzzi.
\newblock Refinement of polygonal grids using convolutional neural networks
  with applications to polygonal {D}iscontinuous {Galerkin} and {V}irtual
  {E}lement methods.
\newblock {\em Journal of Computational Physics}, 452:110900, 2022.

\bibitem{AntoniettiMelas_2020}
P.~F. Antonietti and L.~Melas.
\newblock Algebraic multigrid schemes for high-order nodal discontinuous
  {G}alerkin methods.
\newblock {\em SIAM Journal on Scientific Computing}, 42(2):A1147--A1173, 2020.

\bibitem{dealII91}
D.~Arndt, W.~Bangerth, T.~C. Clevenger, D.~Davydov, M.~Fehling,
  D.~Garcia-Sanchez, G.~Harper, T.~Heister, L.~Heltai, M.~Kronbichler, R.~M.
  Kynch, M.~Maier, J.-P. Pelteret, B.~Turcksin, and D.~Wells.
\newblock The \texttt{deal.II} library, version 9.1.
\newblock {\em Journal of Numerical Mathematics}, 27(4):203--213, 2019.

\bibitem{falgout96mgcg}
S.~F. Ashby and R.~D. Falgout.
\newblock A parallel multigrid preconditioned conjugate gradient algorithm for
  groundwater flow simulations.
\newblock {\em Nuclear Science and Engineering}, 124(1):145--159, 1996.

\bibitem{baker2012scaling}
A.~H. Baker, R.~D. Falgout, T.~V. Kolev, and U.~M. Yang.
\newblock Scaling hypre's multigrid solvers to 100,000 cores.
\newblock In {\em High-Performance Scientific Computing}, pages 261--279.
  Springer, 2012.

\bibitem{baker2010improving}
A.~H. Baker, T.~V. Kolev, and U.~M. Yang.
\newblock Improving algebraic multigrid interpolation operators for linear
  elasticity problems.
\newblock {\em Numerical Linear Algebra with Applications}, 17(2-3):495--517,
  2010.

\bibitem{bank2015algebraic}
R.~Bank, R.~Falgout, T.~Jones, T.~A. Manteuffel, S.~F. McCormick, and J.~W.
  Ruge.
\newblock Algebraic multigrid domain and range decomposition
  ({AMG}-{DD}/{AMG}-{RD}).
\newblock {\em SIAM Journal on Scientific Computing}, 37(5):S113--S136, 2015.

\bibitem{BastianBlattScheichl_2012}
P.~Bastian, M.~Blatt, and R.~Scheichl.
\newblock Algebraic multigrid for discontinuous {G}alerkin discretizations of
  heterogeneous elliptic problems.
\newblock {\em Numerical Linear Algebra with Applications}, 19(2):367--388,
  2012.

\bibitem{boffi2013mixed}
D.~Boffi, F.~Brezzi, M.~Fortin, et~al.
\newblock {\em Mixed finite element methods and applications}, volume~44.
\newblock Springer, 2013.

\bibitem{Brandt}
A.~Brandt and O.~E. Livne.
\newblock {\em Multigrid Techniques}.
\newblock Society for Industrial and Applied Mathematics, 2011.

\bibitem{Brannick}
J.~Brannick, M.~Brezina, S.~MacLachlan, T.~Manteuffel, S.~McCormick, and
  J.~Ruge.
\newblock An energy-based {AMG} coarsening strategy.
\newblock {\em Numerical Linear Algebra with Applications}, 13(2-3):133--148,
  2006.

\bibitem{brezinafalgout2002AMGe}
M.~Brezina, A.~Cleary, R.~Falgout, V.~Henson, J.~Jones, T.~Manteuffel,
  S.~McCormick, and J.~Ruge.
\newblock Algebraic multigrid based on element interpolation ({AMGe}).
\newblock {\em SIAM Journal on Scientific Computing}, 22, 2002.

\bibitem{brezis2010functional}
H.~Brezis.
\newblock {\em Functional Analysis, Sobolev Spaces and Partial Differential
  Equations}.
\newblock Springer Science \& Business Media, 2010.

\bibitem{buiwnagoseikuffor2018}
Q.~M. Bui, L.~Wang, and D.~Osei-Kuffuor.
\newblock Algebraic multigrid preconditioners for two-phase flow in porous
  media with phase transitions.
\newblock {\em Advances in Water Resources}, 114:19--28, 2018.

\bibitem{cleary2000robustness}
A.~J. Cleary, R.~D. Falgout, V.~E. Henson, J.~E. Jones, T.~A. Manteuffel, S.~F.
  McCormick, G.~N. Miranda, and J.~W. Ruge.
\newblock Robustness and scalability of algebraic multigrid.
\newblock {\em SIAM Journal on Scientific Computing}, 21(5):1886--1908, 2000.

\bibitem{discacciati2020controlling}
N.~Discacciati, J.~S. Hesthaven, and D.~Ray.
\newblock Controlling oscillations in high-order discontinuous {G}alerkin
  schemes using artificial viscosity tuned by neural networks.
\newblock {\em Journal of Computational Physics}, 409:109304, 2020.

\bibitem{falgoutvassilevski2004}
R.~Falgout and S.~Vassilevski.
\newblock On generalizing the amg framework.
\newblock {\em SIAM Journal on Scientific Computing}, 42(4):1669--1693, 2004.

\bibitem{falgout02hypre}
R.~Falgout and U.~Yang.
\newblock hypre: A library of high performance preconditioners.
\newblock {\em Computational Science-ICCS 2002, Pt Iii, Proceedings},
  2331:632--641, 04 2002.

\bibitem{fresca2021comprehensive}
S.~Fresca, L.~Dede', and A.~Manzoni.
\newblock A comprehensive deep learning-based approach to reduced order
  modeling of nonlinear time-dependent parametrized {PDEs}.
\newblock {\em Journal of Scientific Computing}, 87(2):1--36, 2021.

\bibitem{goodfellow2016deep}
I.~Goodfellow, Y.~Bengio, and A.~Courville.
\newblock {\em Deep Learning}.
\newblock MIT Press, 2016.

\bibitem{Gottschalk2021}
H.~Gottschalk and K.~Kahl.
\newblock Coarsening in algebraic multigrid using gaussian processes.
\newblock {\em ETNA - Electronic Transactions on Numerical Analysis},
  54:514–533, 2021.

\bibitem{greenfeld2019learning}
D.~Greenfeld, M.~Galun, R.~Basri, I.~Yavneh, and R.~Kimmel.
\newblock Learning to optimize multigrid {PDE} solvers.
\newblock In {\em International Conference on Machine Learning}, pages
  2415--2423. PMLR, 2019.

\bibitem{he2015delving}
K.~He, X.~Zhang, S.~Ren, and J.~Sun.
\newblock Delving deep into rectifiers: Surpassing human-level performance on
  imagenet classification.
\newblock In {\em Proceedings of the IEEE international conference on computer
  vision}, pages 1026--1034, 2015.

\bibitem{he2016deep}
K.~He, X.~Zhang, S.~Ren, and J.~Sun.
\newblock Deep residual learning for image recognition.
\newblock In {\em Proceedings of the IEEE conference on computer vision and
  pattern recognition}, pages 770--778, 2016.

\bibitem{Heinlein_et_al_2021}
A.~Heinlein, A.~Klawonn, M.~Lanser, and J.~Weber.
\newblock Combining machine learning and adaptive coarse spaces-a hybrid
  approach for robust {FETI-DP} methods in three dimensions.
\newblock {\em SIAM Journal on Scientific Computing}, 0(0):S816--S838, 2021.

\bibitem{hensonvassilevski2001}
V.~E. Henson and P.~S. Vassilevski.
\newblock Algebraic multigrid preconditioners for two-phase flow in porous
  media with phase transitions.
\newblock {\em SIAM Journal on Scientific Computing}, 23(2):629--650, 2001.

\bibitem{hesthaven2018non}
J.~S. Hesthaven and S.~Ubbiali.
\newblock Non-intrusive reduced order modeling of nonlinear problems using
  neural networks.
\newblock {\em Journal of Computational Physics}, 363:55--78, 2018.

\bibitem{hu2018squeeze}
J.~Hu, L.~Shen, and G.~Sun.
\newblock Squeeze-and-excitation networks.
\newblock In {\em Proceedings of the IEEE conference on computer vision and
  pattern recognition}, pages 7132--7141, 2018.

\bibitem{hughes2012finite}
T.~J.~R. Hughes.
\newblock {\em The Finite Element Method: Linear Static and Dynamic Finite
  Element Analysis}.
\newblock Courier Corporation, 2012.

\bibitem{ioffe2015batch}
S.~Ioffe and C.~Szegedy.
\newblock Batch normalization: Accelerating deep network training by reducing
  internal covariate shift.
\newblock {\em arXiv preprint arXiv:1502.03167}, 2015.

\bibitem{janssensadvancing}
M.~Janssens and S.~Hulshoff.
\newblock Advancing artificial neural network parameterisation for atmospheric
  turbulence using a variational multiscale model.
\newblock {\em Journal of Advances in Modeling Earth Systems}, page
  e2021MS002490, 2021.

\bibitem{joneslee2006}
J.~Jones and B.~Lee.
\newblock A multigrid method for variable coefficient maxwell's equations.
\newblock {\em SIAM Journal on Scientific Computing}, 27(5):1689--1708, 2006.

\bibitem{KATRUTSA2020112524}
A.~Katrutsa, T.~Daulbaev, and I.~Oseledets.
\newblock Black-box learning of multigrid parameters.
\newblock {\em Journal of Computational and Applied Mathematics}, 368:112524,
  2020.

\bibitem{kingma2014adam}
D.~P. Kingma and J.~Ba.
\newblock Adam: A method for stochastic optimization.
\newblock {\em arXiv preprint arXiv:1412.6980}, 2014.

\bibitem{catastrophicforgetting}
J.~Kirkpatrick, R.~Pascanu, N.~Rabinowitz, J.~Veness, G.~Desjardins, A.~A.
  Rusu, K.~Milan, J.~Quan, T.~Ramalho, A.~Grabska-Barwinska, D.~Hassabis,
  C.~Clopath, D.~Kumaran, and R.~Hadsell.
\newblock Overcoming catastrophic forgetting in neural networks.
\newblock {\em Proceedings of the National Academy of Sciences},
  114(13):3521--3526, 2017.

\bibitem{krizhevsky2012imagenet}
A.~Krizhevsky, I.~Sutskever, and G.~E. Hinton.
\newblock Imagenet classification with deep convolutional neural networks.
\newblock In {\em Advances in Neural Information Processing Systems}, pages
  1097--1105, 2012.

\bibitem{lecun1998gradient}
Y.~LeCun, L.~Bottou, Y.~Bengio, and P.~Haffner.
\newblock Gradient-based learning applied to document recognition.
\newblock {\em Proceedings of the IEEE}, 86(11):2278--2324, 1998.

\bibitem{li2020efficient}
R.~Li and C.~Zhang.
\newblock Efficient parallel implementations of sparse triangular solves for
  gpu architectures.
\newblock In {\em Proceedings of the 2020 SIAM Conference on Parallel
  Processing for Scientific Computing}, pages 106--117. SIAM, 2020.

\bibitem{lions72}
J.~L. Lions and E.~Magenes.
\newblock {\em Non-Homogeneous Boundary Value Problems and Applications}.
\newblock Springer-Verlag Berlin Heidelberg, 1972.

\bibitem{Mishra}
S.~Mishra.
\newblock A machine learning framework for data driven acceleration of
  computations of differential equations.
\newblock {\em Mathematics in Engineering}, 1(1):118--146, 2019.

\bibitem{Neittaanmaki202227}
P.~Neittaanm\"{a}ki and S.~Repin.
\newblock Artificial intelligence and computational science.
\newblock {\em Intelligent Systems, Control and Automation: Science and
  Engineering}, 76:27--35, 2022.

\bibitem{Quarteroni_2017}
A.~Quarteroni.
\newblock {\em Numerical Models for Differential Problems}, volume~16.
\newblock Springer International Publishing, 2017.

\bibitem{QV94}
A.~Quarteroni and A.~Valli.
\newblock {\em Numerical Approximation of Partial Differential Equations}.
\newblock Springer--Verlag, Berlin and Heidelberg, 1994.

\bibitem{raissi2017machine}
M.~Raissi, P.~Perdikaris, and G.~E. Karniadakis.
\newblock Machine learning of linear differential equations using {G}aussian
  processes.
\newblock {\em Journal of Computational Physics}, 348:683--693, 2017.

\bibitem{raissi2019physics}
M.~Raissi, P.~Perdikaris, and G.~E. Karniadakis.
\newblock Physics-informed neural networks: A deep learning framework for
  solving forward and inverse problems involving nonlinear partial differential
  equations.
\newblock {\em Journal of Computational Physics}, 378:686--707, 2019.

\bibitem{regazzoni2019machine}
F.~Regazzoni, L.~Dede', and A.~Quarteroni.
\newblock Machine learning for fast and reliable solution of time-dependent
  differential equations.
\newblock {\em Journal of Computational Physics}, 397:108852, 2019.

\bibitem{ruge}
J.~W. Ruge and K.~Stüben.
\newblock {\em 4. Algebraic Multigrid}, pages 73--130.
\newblock SIAM, 1987.

\bibitem{russakovsky2015imagenet}
O.~Russakovsky, J.~Deng, H.~Su, J.~Krause, S.~Satheesh, S.~Ma, Z.~Huang,
  A.~Karpathy, A.~Khosla, M.~Bernstein, et~al.
\newblock Imagenet large scale visual recognition challenge.
\newblock {\em International Journal of Computer Vision}, 115(3):211--252,
  2015.

\bibitem{seabold2010statsmodels}
S.~Seabold and J.~Perktold.
\newblock statsmodels: Econometric and statistical modeling with python.
\newblock In {\em 9th Python in Science Conference}, 2010.

\bibitem{SiefertTuminaro_et_al_2014}
C.~Siefert, R.~Tuminaro, A.~Gerstenberger, G.~Scovazzi, and S.~S. Collis.
\newblock Algebraic multigrid techniques for discontinuous {G}alerkin methods
  with varying polynomial order.
\newblock {\em Computational Geosciences}, 18(5):597--612, 2014.

\bibitem{srivastava2014dropout}
N.~Srivastava, G.~Hinton, A.~Krizhevsky, I.~Sutskever, and R.~Salakhutdinov.
\newblock Dropout: a simple way to prevent neural networks from overfitting.
\newblock {\em The Journal of Machine Learning Research}, 15(1):1929--1958,
  2014.

\bibitem{trottenberg2001multigrid}
K.~St{\"u}ben.
\newblock {\em An Introduction to Algebraic Multigrid}, chapter~A, pages
  413--532.
\newblock Elsevier Science, 2001.

\bibitem{STUBEN2001281}
K.~Stüben.
\newblock A review of algebraic multigrid.
\newblock {\em Journal of Computational and Applied Mathematics},
  128(1):281--309, 2001.
\newblock Numerical Analysis 2000. Vol. VII: Partial Differential Equations.

\bibitem{sutskever2013importance}
I.~Sutskever, J.~Martens, G.~Dahl, and G.~Hinton.
\newblock On the importance of initialization and momentum in deep learning.
\newblock In {\em International conference on machine learning}, pages
  1139--1147, 2013.

\bibitem{tassi2021amachine}
T.~Tassi, A.~Zingaro, and L.~Dede'.
\newblock A {M}achine {L}earning approach to enhance the {S}{U}{P}{G}
  stabilization method for advection-dominated differential problems.
\newblock {\em MOX Report, Politecnico di Milano}, 58, 2021.

\bibitem{tieleman2012lecture}
T.~Tieleman and G.~Hinton.
\newblock Lecture 6.5-rmsprop: Divide the gradient by a running average of its
  recent magnitude.
\newblock {\em COURSERA: Neural networks for machine learning}, 4(2):26--31,
  2012.

\bibitem{Waskom2021}
M.~L. Waskom.
\newblock seaborn: statistical data visualization.
\newblock {\em Journal of Open Source Software}, 6(60):3021, 2021.

\bibitem{webster1994algebraic}
R.~Webster.
\newblock An algebraic multigrid solver for {Navier-Stokes} problems.
\newblock {\em International Journal for Numerical Methods in Fluids},
  18(8):761--780, 1994.

\bibitem{xuzikatanov2017}
J.~Xu and L.~Zikatanov.
\newblock Algebraic multigrid methods.
\newblock {\em Acta Numerica}, 26:591--721, 2017.

\end{thebibliography}
\end{document}